\documentclass[10pt]{article}
\usepackage{geometry}
\usepackage[utf8]{inputenc} 
\usepackage[T1]{fontenc}    
\usepackage{microtype}
\usepackage{url}            
\usepackage{booktabs} 
\usepackage{hyperref}
\usepackage{authblk}
\usepackage{xcolor}         
\usepackage{amsfonts}       
\usepackage{nicefrac}       
\usepackage{soul}
\usepackage{lipsum}

\makeatletter
\def\blfootnote{\xdef\@thefnmark{}\@footnotetext}
\makeatother

\usepackage{amsmath,amssymb,amsthm}
\usepackage{graphicx,subcaption}
\usepackage{mathdots}
\usepackage{multirow}
\usepackage{blkarray}
\usepackage{circuitikz}
\usepackage{mathtools}
\usepackage{csquotes}
\usepackage[capitalize,noabbrev]{cleveref}
\usepackage{enumitem}
\usepackage{empheq}
\usepackage{booktabs}
\graphicspath{ {./figure/} }

\usepackage[ruled,vlined]{algorithm2e}
\SetAlFnt{\small}

\numberwithin{equation}{section}

\theoremstyle{remark}

\crefname{assumption}{Assumption}{Assumptions}

\newcommand\norm[1]{\left\Vert#1\right\Vert}
\newcommand\R{\mathbb{R}}

\newcommand\cL{\mathcal{L}}
\newcommand\cB{\mathcal{B}}
\newcommand\cC{\mathcal{C}}
\newcommand\cF{\mathcal{F}}

\title{A Neumann-Neumann Acceleration with Coarse Space for Domain Decomposition of Extreme Learning Machines
\blfootnote{This work was supported by the National Research Foundation~(NRF) of Korea grant funded by the Korea government~(MSIT) (No. RS-2023-00208914).
}}
\author{Chang-Ock Lee and Byungeun Ryoo
}
\affil{Department of Mathematical Sciences, KAIST, Daejeon 34141, Korea}
\date{}

\begin{document}
\maketitle

\begin{abstract}
  Extreme learning machines~(ELMs), which preset hidden layer parameters and solve for last layer coefficients via a least squares method, can typically solve partial differential equations faster and more accurately than Physics Informed Neural Networks.
  However, they remain computationally expensive when high accuracy requires large least squares problems to be solved.
  Domain decomposition methods~(DDMs) for ELMs have allowed parallel computation to reduce training times of large systems.
  This paper constructs a coarse space for ELMs, which enables further acceleration of their training.
  By partitioning interface variables into coarse and non-coarse variables, selective elimination introduces a Schur complement system on the non-coarse variables with the coarse problem embedded.
  Key to the performance of the proposed method is a Neumann-Neumann acceleration that utilizes the coarse space.
  Numerical experiments demonstrate significant speedup compared to a previous DDM method for ELMs.

\end{abstract}


{\small \textbf{AMS subject classifications~(2020)}
65N55, 65Y05, 68T07
}

\section{Introduction}
\label{Sec:Int}
The Universal Approximation Theorem~\cite{cybenko1989approximation,hornik1989multilayer,huang2006universal,pinkus1999approx} has facilitated the study of neural networks as solutions to partial differential equations~(PDEs).
The most widely used method among them is the Physics Informed Neural Network~(PINN)~\cite{raissi2019physics}, which typically employ gradient based optimization methods.
Extreme learning machines~(ELMs)~\cite{huang2011extreme,huang2006extreme} instead use least squares methods to solve for the last layer coefficients by fixing the hidden layer parameters.
Applied to PDEs, ELMs often outperform PINNs in terms of both accuracy and speed~\cite{dong2021local,dwivedi2020physics}.
Nevertheless, the time complexity of the least squares problem limits ELMs' applicability to large systems.

Domain decomposition methods~(DDMs) are a class of fast solvers for large systems suitable for parallel computation.
Methods mainly fall into two categories: overlapping and non-overlapping methods.
In substructuring methods, which are non-overlapping methods, a reduced problem on the interface of non-overlapping subdomains is solved after which the local subproblems can be solved independently.
Choosing which reduced problem to solve and how to accelerate it leads to different methods, Neumann-Neumann~\cite{bourgat1988variational} and FETI~\cite{farhat1991method} methods being two prominent examples.
For many problems, scalability can only be achieved by adding a coarse space~\cite{toselli2005domain}.
Substructuring methods construct coarse spaces by partitioning interface variables into coarse and non-coarse components.
Eliminating all but the non-coarse interface variables embeds a global problem coupled only in the coarse variables into the reduced system.
BDDC~\cite{dohrmann2003preconditioner} and FETI-DP~\cite{farhat2001feti} are two such methods.

Application of ideas from DDM to PINNs has resulted in many accelerated training regimes.
These include overlapping methods such as DeepDDM~\cite{li2020deep}, FBPINN~\cite{moseley2023finite}.
Several works have incorporated two-level~\cite{valentin2021coarse,yang2022additive} and multilevel schemes~\cite{dolean2024multilevel} to these methods.
Nonoverlapping methods have also been developed~\cite{basir2023generalized,jagtap2020extended,jagtap2020conservative}.

Domain decomposition techniques for ELMs have also been developed, among which are Local ELM~\cite{dong2021local} and Distributed PIELM~\cite{dwivedi2020physics}.
However in both methods, the local network parameters cannot be optimized in parallel.
Our previous work, Domain Decomposition for ELM~(DDELM)~\cite{lee2024nonoverlapping}, overcomes this limitation by applying techniques motivated by substructuring methods.
In DDELM an auxiliary variable on the interface is introduced to couple the local networks.
A Schur complement reduces the system to one on the interface variables.
DDELM showed that parallel computation of the interface problem can reduce computation time.

In this paper, we extend DDELM by introducing a coarse space.
We first partition the auxiliary interface variables into coarse and non-coarse components.
After eliminating the flux terms on the coarse variables, we obtain a Schur complement system on the non-coarse variables with an embedded coarse problem.

A key term appearing in the resulting systems serves as a Dirichlet to Neumann map.
Driven by this observation, we introduce a Neumann-Neumann acceleration by constructing a corresponding Neumann to Dirichlet map.
Although this does not yield a preconditioner for the entire system, it can be understood to precondition some internal steps of the method.

Experimentally, choosing corner variables as the coarse component did not result in significant acceleration.
Furthermore, the Neumann-Neumann acceleration even degraded convergence in this setting.
We additionally introduce a change of variables that makes the Neumann-Neumann acceleration faster and more stable.
Numerical experiments show that the combination of the change of variables and Neumann-Neumann acceleration greatly improves convergence.

The paper is organized as follows.
In \cref{Sec:Prelim}, we give a brief description of the classical DDMs related to this research and our previous work on DDELM.
In \cref{Sec:Coarse}, we present the coarse space and Neumann-Neumann acceleration for DDELM.
We also detail the change of variables and a reweighting scheme for the loss terms.
In \cref{Sec:Num}, we provide an ablation study to verify the effectiveness of various components of our acceleration scheme and present the results of numerical experiments.
Performance is compared to the original DDELM algorithm for Poisson and biharmonic equations.
We conclude with remarks in \cref{Sec:Conc}.

\section{Preliminaries}\label{Sec:Prelim}
\label{Sec:DDM}
We consider a general PDE with boundary conditions:
\begin{equation}
  \label{eqn:pde}
  \left\{\begin{aligned}
  \cL u &= f && \text{ in } \Omega,\\
  \cB u &= g && \text{ on } \partial\Omega,
  \end{aligned}\right.
\end{equation}
where $\Omega \subset \R^{d}$ is a bounded domain, $\partial \Omega$ denotes the boundary of $\Omega$, $f$ and $g$ are functions, $\cL$ is a differential operator,
and $\mathcal{B}$ represents boundary conditions, e.g. $\mathcal{B}u=\left. u\right|_\Omega=g$ gives Dirichlet boundary conditions.
We assume that problem~\eqref{eqn:pde} is well-posed so that the solution exists uniquely.

\subsection{Nonoverlapping domain decomposition methods}
Let $\Omega$ be divided into $N$ nonoverlapping subdomains $\Omega^i$.
For certain PDEs, a coupled system equivalent to the problem~\eqref{eqn:pde} can be constructed using appropriate continuity conditions $\cC$ and flux operator $\cF$:
\begin{equation}
\label{eqn:nonoverlap}
\left\{\begin{aligned}
  \cL u^i &= f && \text{ in } \Omega^i,\\
  \cB u^i &= g && \text{ on } \partial\Omega^i \cap \partial\Omega,\\
  \cC u^i &= \cC u^j && \text{ on } \partial\Omega^i \cap \partial\Omega^j,\\
  \cF u^i &= -\cF u^j && \text{ on } \partial\Omega^i \cap \partial\Omega^j.
\end{aligned}\right.
\end{equation}
Taking the biharmonic equation as an example, the equation
\begin{equation*}
  \left\{\begin{aligned}
  \Delta^2 u &= f && \text{ in } \Omega,\\
  u &= g_1 && \text{ on } \partial\Omega,\\
  \Delta u &= g_2 && \text{ on } \partial\Omega
  \end{aligned}\right.
\end{equation*}
is equivalent to the coupled system
\begin{equation}
  \label{eqn:biharmonic_coupled}
  \left\{\begin{aligned}
  \Delta^2 u^i &= f && \text{ in } \Omega^i,\\
  u^i &= g_1 && \text{ on } \partial\Omega^i \cap \partial\Omega,\\
  \Delta u^i &= g_2 && \text{ on } \partial\Omega^i \cap \partial\Omega,\\
  u^i &= u^j && \text{ on } \partial\Omega^i \cap \partial\Omega^j,\\
  \Delta u^i &= \Delta u^j && \text{ on } \partial\Omega^i \cap \partial\Omega^j,\\
  \frac{\partial u^i}{\partial n^i} &= -\frac{\partial u^j}{\partial n^j} && \text{ on } \partial\Omega^i \cap \partial\Omega^j,\\
  \frac{\partial\Delta u^i}{\partial n^i} &= -\frac{\partial\Delta u^j}{\partial n^j} && \text{ on } \partial\Omega^i \cap \partial\Omega^j,
  \end{aligned}\right.
\end{equation}
where $n^i$ is the outward normal vector of $\partial\Omega^i$.

\subsubsection{Dirichlet-Dirichlet methods}
This class of methods aims to solve for the flux on the interface, after which the local subproblems can be solved in parallel.
They are called Dirichlet-Dirichlet methods because the system for the flux is essentially a Neumann to Dirichlet map and local Dirichlet to Neumann maps are used to form a preconditioner.
FETI with Dirichlet preconditioner belongs to this category.
In its most basic form, FETI suffers in two aspects.
One is that domain decomposition methods require a coarse problem to be scalable, but neither the system for the flux nor the accompanying preconditioner incorporates one.
Another is that when a subdomain does not share boundaries with the entire domain, the Neumann to Dirichlet map has a nontrivial null space.

FETI-DP is a variant of FETI that introduces a coarse space into its construction.
It chooses certain degrees of freedom on the interface to be primal, which is to say at each iteration they hold the same value across subdomains.
This creates a system coupled only on these primal variables, which is embedded into the system for the flux and acts as a coarse problem.
Appropriate choice of the primal degrees of freedom confers FETI-DP scalability~\cite{mandel2001convergence}, while the primal variables removes the singularity of the Neumann to Dirichlet map. 
In simple 2D elliptic equations, a coarse space consisting only of corner variables is enough to achieve scalability.

\subsubsection{Neumann-Neumann method}
Neumann-Neumann methods are dual to the Dirichlet-Dirichlet methods, solving for the interface degrees of freedom with local Neumann to Dirichlet maps as preconditioner.
The system for the interface degrees of freedom is a Dirichlet to Neumann map so there is no nontrivial null space even with floating subdomains; however, the preconditioner has a nontrivial null space.
In continuous form, the basic form of the method can be written as
\begin{align*}
  &(D_i)\left\{
         \begin{aligned}
                \cL u_{n+1/2}^i &= f && \text{ in } \Omega_i,\\
                \cB u_{n+1/2}^i &= g && \text{ on } \partial\Omega_i \cap \partial \Omega\\
                \cC u_{n+1/2}^i &= \mu_n && \text{ on } \partial\Omega_i \backslash \partial \Omega
         \end{aligned}
  \right.\\
  &(N_i)\left\{
         \begin{aligned}
                \cL v_{n+1}^i &= 0 && \text{ in } \Omega_i,\\
                \cB v_{n+1}^i &= 0 && \text{ on } \partial\Omega_i \cap \partial \Omega\\
                \cF v_{n+1}^i &= \sum_i \cF u_{n+1/2}^i && \text{ on } \partial\Omega_i \backslash \partial \Omega
         \end{aligned}
  \right.\\
  &\mu_{n+1} = \mu_n - \theta \sum_i \cC v_{n+1}^i.
  \end{align*}
Here $\theta$ is the step size for the Richardson iteration that the continuous form is written in.

\subsection{Extreme learning machines}
ELMs optimize neural networks, typically single hidden layer multi layer perceptrons, by fixing the weights and biases of the hidden layers and solving for the last layer coefficients via a least squares algorithm.
In this paper, we will consider single hidden layer multi layer perceptrons of the form
\begin{align*}
N\!N(x) &= \sum_{j=1}^{M} c_j \sigma(w_j^T x + b_j)\\
& \eqqcolon \sum_{j=1}^{M} c_j \phi_j(x),
\end{align*}
where $M$ is the number of neurons in the hidden layer.

Let $x_I$, $x_B$ be points in the interior and boundary of $\Omega$, respectively.
An ELM solution to problem~\eqref{eqn:pde} is obtained by randomly initializing $\{w_j\}$ and $\{b_j\}$, then solving for $\{c_j\}$ in the least squares problem
\[
\min_{\{c_j\}} \norm{\cL N\!N(x_I) - f(x_I)}^2 + \norm{\cB N\!N(x_B) - g(x_B)}^2.
\]
For a linear operator $\cL$, we solve the least squares system
\[
\begin{bmatrix}
  \cL \phi_1(x_I) & \cL \phi_2(x_I) & \cdots & \cL \phi_M(x_I)\\
  \cB \phi_1(x_B) & \cB \phi_2(x_B) & \cdots & \cB \phi_M(x_B)
\end{bmatrix}
\begin{bmatrix}
  c_1\\
  c_2\\
  \vdots\\
  c_M
\end{bmatrix}
=
\begin{bmatrix}
  f(x_I)\\
  g(x_B)
\end{bmatrix}.
\]

\subsubsection{Domain decomposition for ELM~(DDELM)}
\label{subsubsec:ddelm}
We shall denote by superscript the local networks,
\begin{align*}
N\!N^i(x) &=\sum_{j=1}^{M^i} c_j^i \sigma(w_j^{i}\cdot x+b_j^i) \\
&\eqqcolon \sum_{j=1}^{M^i} c_j^i \phi_j^i(x),
\end{align*}
with $M^i$ the number of neurons in the hidden layer of $N\!N^i$.
Let $x_I^i$, $x_B^i$, $x_\Gamma^i$ be points in the interior, boundary, and interface of $\Omega^i$, respectively.
We use $x_\Gamma$ to refer to the union $\bigcup_{i=1}^{N} x_\Gamma^i$.
By introducing an auxiliary variable $\mu$ on the interface, the condition
\[
  \cC u^i = \cC u^j \text{ on } \partial\Omega^i \cap \partial\Omega^j
\]
in system~\eqref{eqn:nonoverlap} can be replaced with the equivalent condition
\[
  \cC u^i = \mu \text{ on } \partial\Omega^i \backslash \partial\Omega.
\]
Discretization of \eqref{eqn:nonoverlap} then yields
\begin{equation*}
\begin{bmatrix}
  \begin{bmatrix}
    \cL \phi^1(x_I^1)\\
    \cB \phi^1(x_B^1)\\
    \cC \phi^1(x_\Gamma^1)
  \end{bmatrix}
  & 0 & \cdots & 0 & B^1 \vspace{3pt}\\
  0 & \begin{bmatrix}
    \cL \phi^2(x_I^2)\\
    \cB \phi^2(x_B^2)\\
    \cC \phi^2(x_\Gamma^2)
\end{bmatrix} & \cdots & 0 & B^2\\
\vdots & \vdots & \ddots & \vdots & \vdots \vspace{3pt}\\
0 & 0 & \cdots & \begin{bmatrix}
  \cL \phi^N(x_I^N)\\
  \cB \phi^N(x_B^N)\\
  \cC \phi^N(x_\Gamma^N)
\end{bmatrix} & B^N \vspace{3pt}\\
\cF\phi^1(x_\Gamma) & \cF\phi^2(x_\Gamma) & \cdots & \cF\phi^N(x_\Gamma) & 0
\end{bmatrix}
\begin{bmatrix}
  c^1\\
  c^2\\
  \vdots\\
  c^N\\
  \mu(x_\Gamma)
\end{bmatrix}=
\begin{bmatrix}
  \begin{bmatrix}
    f(x_I^1)\\
    g(x_B^1)\\
    0
  \end{bmatrix}\vspace{3pt}\\
  \begin{bmatrix}
    f(x_I^2)\\
    g(x_B^2)\\
    0
  \end{bmatrix}\\
  \vdots \vspace{3pt}\\
  \begin{bmatrix}
    f(x_I^N)\\
    g(x_B^N)\\
    0
  \end{bmatrix}\vspace{3pt}\\
  0
\end{bmatrix},
\end{equation*}
where $B^i=\begin{bmatrix}0&0&-(R^i)^{T}\end{bmatrix}^T$ and $R^i$ is the restriction operator taking $\mu(x_\Gamma)$ to $\mu(x_\Gamma^i)$.
Note that, in the last row, $\cF\phi^i(x_\Gamma)$ is understood to be nonzero only where $x_\Gamma$ intersects with $x_\Gamma^i$.
Let
\begin{equation*}
  \begin{gathered}
    K^i=\begin{bmatrix}\cL \phi^i(x_I^i)^T & \cB \phi^i(x_B^i)^T & \cC \phi^i(x_\Gamma^i)^T\end{bmatrix}^T,\quad
    K = \mathrm{diag}(K^1, K^2, \ldots, K^N),\\
    B = \begin{bmatrix}(B^1)^{T} & (B^2)^{T} & \cdots & (B^N)^{T}\end{bmatrix}^T,\\
    A^i = \cF \phi^i(x_\Gamma),\quad
    A = \begin{bmatrix}A^{1} & A^{2} & \cdots & A^{N}\end{bmatrix},\\
    f^i = \begin{bmatrix}f(x_I^i)^T & g(x_B^i)^T & 0\end{bmatrix}^T,\quad f = \begin{bmatrix}(f^1)^T & (f^2)^T & \cdots & (f^N)^T\end{bmatrix}^T,\\
    c = \begin{bmatrix}(c^1)^T & (c^2)^T & \cdots & (c^N)^T\end{bmatrix}^T
  \end{gathered}
\end{equation*}
to rewrite the system as
\begin{equation}
  \label{eqn:ddelm}
  \begin{bmatrix}
    K & B\\
    A & 0
  \end{bmatrix}
  \begin{bmatrix}
    c\\
    \mu
  \end{bmatrix}
  =
  \begin{bmatrix}
    f\\
    0
  \end{bmatrix}.
\end{equation}
Let $D=-AK^+B$ and $d=-AK^+f$, where $K^+$ is the pseudo-inverse of $K$.
Instead of solving~\eqref{eqn:ddelm} directly, we first eliminate $A$ from the second row to obtain
\begin{equation}
  \label{eqn:schur}
  \begin{bmatrix}
    K & B\\ 0 & D
  \end{bmatrix}
  \begin{bmatrix}
    c\\
    \mu
  \end{bmatrix}
  =
  \begin{bmatrix}
    f \\ d
  \end{bmatrix}
\end{equation}
and form the normal equation
\begin{equation*}
  \begin{bmatrix}
    K^TK & K^TB\\
    B^TK & B^TB+D^TD
  \end{bmatrix}
  \begin{bmatrix}
    c\\
    \mu
  \end{bmatrix}
  =
  \begin{bmatrix}
    K^Tf\\
    B^Tf+D^Td
  \end{bmatrix}.
\end{equation*}
Taking the Schur complement gives the solution
\begin{equation*}
  \left\{
    \begin{aligned}
      (B^TB+D^TD-B^TKK^+B)\mu &= (B^T-B^TKK^+)f+D^Td,\\
      c&=K^+(f-B\mu).
    \end{aligned}
  \right.
\end{equation*}
We may first find $\mu$ via the conjugate gradient~(CG) method from which $c$ can be obtained through a local least squares system.

\section{Coarse Space Formulation for DDELM}
\label{Sec:Coarse}
Partition and reorder $\Gamma$ into corner points, $\Pi$, and the rest, $\Delta$.
The system~\eqref{eqn:ddelm} can then be written as
\begin{equation*}
  \begin{bmatrix}
    K & B_\Pi & B_\Delta\\
    A_\Pi & 0 & 0\\
    A_\Delta & 0 & 0
  \end{bmatrix}
  \begin{bmatrix}
    c\\
    \mu_\Pi\\
    \mu_\Delta
  \end{bmatrix}
  =
  \begin{bmatrix}
    f\\
    0\\
    0
  \end{bmatrix}.
\end{equation*}
Eliminating $A_\Pi$ from the second row, the augmented system becomes
\begin{equation}
  \label{eqn:intermediate}
  \left[\begin{array}{@{}ccc|c@{}}
    K & B_\Pi & B_\Delta & f\\
    0 & -A_\Pi K^+B_\Pi & -A_\Pi K^+B_\Delta & -A_\Pi K^+ f\\
    A_\Delta & 0 & 0 & 0
  \end{array}\right].
\end{equation}
Let
\begin{equation*}
  \begin{gathered}
    D_{\alpha\beta} = -A_\alpha K^+B_\beta = -\sum_i A_\alpha^i (K^i)^+ B_\beta^i \text{ for } \alpha,\beta\in \left\{\Pi, \Delta\right\},\\
    d_\alpha = -A_\alpha K^+ f = -\sum_i A_\alpha^i (K^i)^+ f^i \text{ for } \alpha\in \left\{\Pi, \Delta\right\},\\
    \tilde{K} = \begin{bmatrix}K & B_\Pi\\ 0 & D_{\Pi\Pi}\end{bmatrix},\quad \tilde{B} = \begin{bmatrix}B_\Delta\\ D_{\Pi\Delta} \end{bmatrix},\quad \tilde{A} = \begin{bmatrix}A_\Delta & 0\end{bmatrix},\quad \tilde{f} = \begin{bmatrix}f\\ d_\Pi\end{bmatrix},\\
    \tilde{D} = -\tilde{A}\tilde{K}^+\tilde{B},\quad \tilde{d}= -\tilde{A}\tilde{K}^+\tilde{f}.
  \end{gathered}
\end{equation*}
Then the system~\eqref{eqn:intermediate} can be written as
\begin{equation*}
  \left[\begin{array}{@{}cc|c@{}}
    \tilde{K} & \tilde{B} & \tilde{f}\\
    \tilde{A} & 0 & 0
  \end{array}\right].
\end{equation*}
We now eliminate $\tilde{A}$ from the second row to obtain
\begin{equation}
  \label{eqn:coarse}
  \left[\begin{array}{@{}cc|c@{}}
    \tilde{K} & \tilde{B} & \tilde{f}\\
    0 & \tilde{D} & \tilde{d}
  \end{array}\right].
\end{equation}
Note the similarity to system~\eqref{eqn:schur}.
We thus form the normal equation as in vanilla DDELM and solve for $\mu_\Delta$ then $\mu_\Pi$ and $c$:
\begin{subequations}
\label{eqn:coarse_sol}
\begin{empheq}[left = \empheqlbrace]{align}
  \label{eqn:coarse_sol_1}
  (\tilde{B}^T\tilde{B}+\tilde{D}^T\tilde{D}-\tilde{B}^T\tilde{K}\tilde{K}^+\tilde{B})\mu_\Delta &= (\tilde{B}^T-\tilde{B}^T\tilde{K}\tilde{K}^+)\tilde{f}+\tilde{D}^T\tilde{d},\\
  \label{eqn:coarse_sol_2}
  \begin{bmatrix}
    c \\ \mu_\Pi
  \end{bmatrix} &= \tilde{K}^+(\tilde{f}-\tilde{B}\mu_\Delta) .
\end{empheq}
\end{subequations}
Each of these equations requires the application of $\tilde{K}^+$, the pseudo-inverse of a system coupled on the coarse variables.
Incidentally, $\tilde{K}$ also looks much like \eqref{eqn:schur} so that $
\begin{bmatrix}
  c^T & \mu_\Pi^T
\end{bmatrix}^T = \tilde{K}^+ \begin{bmatrix}
  \varphi^T & \psi^T
\end{bmatrix}^T$ is solved by
\begin{equation*}
  \left\{\begin{aligned}
  (B_\Pi^TB_\Pi+D_{\Pi\Pi}^TD_{\Pi\Pi}-B_\Pi^TKK^+B_\Pi)\mu_\Pi &= (B_\Pi^T-B_\Pi^TKK^+)\varphi+D_{\Pi\Pi}^T\psi\\
  c &= K^+(\varphi-B_\Pi\mu_\Pi)
  \end{aligned}\right. .
\end{equation*}
Let $S_\Pi$ be the system involving $\mu_\Pi$, which is a Schur complement of a normal equation.
We assemble $S_\Pi$ in full and perform a Cholesky factorization for repeated application of its inverse.

Since \eqref{eqn:coarse_sol_1} involves $\tilde{D}^T=-\tilde{B}^T(\tilde{K}^+)^T\tilde{A}^T$, we also need to know how to apply $(\tilde{K}^+)^T$ to $\tilde{A}^T=\begin{bmatrix}
  A_\Delta & 0
\end{bmatrix}^T$.
Hence, we only need to be able to compute
\begin{align*}
(\tilde{K}^+)^{T}\begin{bmatrix}
    \varphi \\ 0
  \end{bmatrix}
&= \tilde{K}
\begin{bmatrix}
  K^TK & K^TB_\Pi\\
  B_\Pi^TK & B_\Pi^TB_\Pi+D_{\Pi\Pi}^TD_{\Pi\Pi}
\end{bmatrix}^{-1}
\begin{bmatrix}
  \varphi\\ 0
\end{bmatrix}\\
&= \tilde{K}
\begin{bmatrix}
  (K^TK)^{-1}\varphi - K^+ B_\Pi \nu\\
  \nu
\end{bmatrix} \rlap{\quad $\left(\nu=S_\Pi^{-1}(-B_\Pi^T(K^+)^T\varphi)\right)$}\\
&= \begin{bmatrix}
  (K^+)^{T} \varphi - KK^+ B_\Pi \nu + B_\Pi \nu\\
  D_{\Pi\Pi} \nu
\end{bmatrix}.
\end{align*}

\subsection{Neumann-Neumann Acceleration}
\label{subsec:nn-accel}
Note how $AK^+B$ takes Dirichlet boundary information, solves a homogeneous problem, then adds the fluxes of neighboring subdomains.
This implies that $AK^+B$ is a Dirichlet to Neumann map.
A similar argument can be made for $\tilde{A}\tilde{K}^+\tilde{B}$.

Let $R_f^i$ be the restriction operator taking the flux on $\Gamma$ to that on $\Gamma^i$ and
\begin{equation*}
  \begin{gathered}
    \bar{K}^i = \begin{bmatrix}
    \cL \phi^i(x_I^i)\\
    \cB \phi^i(x_B^i)\\
    \cC \phi^i(x_\Pi^i)\\
    \cF \phi^i(x_\Delta^i)
    \end{bmatrix},\quad
    \bar{K} = \begin{bmatrix}
      \bar{K}^1 & & &\\
      & \bar{K}^2 & &\\
      & & \ddots &\\
      & & & \bar{K}^N
    \end{bmatrix},\quad
    \bar{B}^i = \begin{bmatrix}
      0\\0\\0\\-R_f^i
    \end{bmatrix},\quad
    \bar{B} = \begin{bmatrix}\bar{B}^{1}\\\bar{B}^{2}\\\vdots\\\bar{B}^{N}\end{bmatrix},\\
    \bar{A}^i = \cC \phi^i(x_\Delta),\quad
    \bar{A} = \begin{bmatrix}\bar{A}^{1}&\bar{A}^{2}&\cdots&\bar{A}^{N}\end{bmatrix}.
  \end{gathered}
\end{equation*}
Then, $\bar{A}\bar{K}^+\bar{B}$ is a Neumann to Dirichlet map.
Note that because the boundary type on the corners in $\bar{K}$ remains Dirichlet, $\bar{K}^+$ corresponds to a well-posed continuous problem.
Multiplying $P\coloneqq\bar{A}\bar{K}^+\bar{B}$ to the second row of \eqref{eqn:coarse} alters the solution in \eqref{eqn:coarse_sol} only by changing \eqref{eqn:coarse_sol_1} to
\begin{equation}
  \label{eqn:nn}
  (\tilde{B}^T\tilde{B}+\tilde{D}^TP^TP\tilde{D}-\tilde{B}^T\tilde{K}\tilde{K}^+\tilde{B})\mu_\Delta = (\tilde{B}^T-\tilde{B}^T\tilde{K}\tilde{K}^+)\tilde{f}+\tilde{D}^TP^TP\tilde{d}
.\end{equation}
Empirically, mixing \eqref{eqn:coarse_sol_1} and \eqref{eqn:nn} into 
\begin{equation}
  \label{eqn:nn-mixed}
  (\tilde{B}^T\tilde{B}+\tilde{D}^T(\theta P^TP + (1-\theta)I)\tilde{D}-\tilde{B}^T\tilde{K}\tilde{K}^+\tilde{B})\mu_\Delta = (\tilde{B}^T-\tilde{B}^T\tilde{K}\tilde{K}^+)\tilde{f}+\tilde{D}^T(\theta P^TP + (1-\theta)I)\tilde{d}
\end{equation}
has improved convergence further for $\theta$ values close to 1.

\subsection{Change of Variables}
\label{subsec:chov}
The coarse space constructed above consists of corner variables.
However, as with classical domain decomposition methods, the coarse space can incorporate more variables or employ a change of variables.
We can simply enrich the coarse space by adding more variables to $\Pi$.

\begin{figure}
  \centering
  \begin{subfigure}{0.495\textwidth}
  \begin{circuitikz}
  \tikzstyle{every node}=[font=\normalsize]
  \draw (6,6) to[short] (11,6);
  \draw (6,9) to[short] (7,6);
  \draw (6,6) to[short] (7,9);
  \draw (8,6) to[short] (7,9);
  \draw (7,6) to[short] (8,9);
  \draw (9,6) to[short] (8,9);
  \draw (8,6) to[short] (9,9);
  \draw (10,6) to[short] (9,9);
  \draw (9,6) to[short] (10,9);
  \draw (11,6) to[short] (10,9);
  \draw (10,6) to[short] (11,9);
  \node at (6,9.25) {$\mu_\Pi$};
  \node at (7,9.25) {$\mu_\Delta$};
  \node at (8,9.25) {$\mu_\Delta$};
  \node at (9,9.25) {$\mu_\Delta$};
  \node at (10,9.25) {$\mu_\Delta$};
  \node at (11,9.25) {$\mu_\Pi$};
  \node at (6,5.75) {1};
  \node at (7,5.75) {2};
  \node at (8,5.75) {3};
  \node at (9,5.75) {4};
  \node at (10,5.75) {5};
  \node at (11,5.75) {6};
  \end{circuitikz}
\end{subfigure}
\hfill
\begin{subfigure}{0.495\textwidth}
  \begin{circuitikz}
  \tikzstyle{every node}=[font=\normalsize]
  \draw (6,6) to[short] (11,6);
  \draw (6,9) to[short] (11,6);
  \draw (6,6) to[short] (7,9);
  \draw (8,6) to[short] (7,9);
  \draw (7,6) to[short] (8,9);
  \draw (9,6) to[short] (8,9);
  \draw (8,6) to[short] (9,9);
  \draw (10,6) to[short] (9,9);
  \draw (9,6) to[short] (10,9);
  \draw (11,6) to[short] (10,9);
  \draw (6,6) to[short] (11,9);
  \node at (6,9.25) {$\mu_\Pi'$};
  \node at (7,9.25) {$\mu_\Delta'$};
  \node at (8,9.25) {$\mu_\Delta'$};
  \node at (9,9.25) {$\mu_\Delta'$};
  \node at (10,9.25) {$\mu_\Delta'$};
  \node at (11,9.25) {$\mu_\Pi'$};
  \node at (6,5.75) {1};
  \node at (7,5.75) {2};
  \node at (8,5.75) {3};
  \node at (9,5.75) {4};
  \node at (10,5.75) {5};
  \node at (11,5.75) {6};
  \end{circuitikz}
\end{subfigure}
\caption{Basis $W$ for auxiliary variable $\mu$ (left) and the new basis $W'$ (right) on an edge.}
\label{fig:change}
\end{figure}
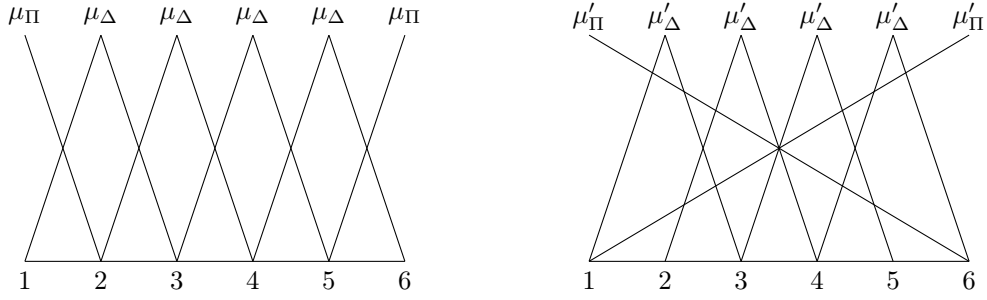

For a change of variables, let $T^i$ be the transition matrix from $W'$ to $W$ in \cref{fig:change}, where the auxiliary variable $\mu_\Pi$ represents point values of $\cC N\!N$, whereas $\mu_\Pi'$ represents a linear function present throughout entire edges of the subdomain.
Here,
\[
T^i=
\begin{bmatrix}
  1 & 0 & 0 & 0 & 0 & 0\\
  4/5 & 1 & 0 & 0 & 0 & 1/5\\
  3/5 & 0 & 1 & 0 & 0 & 2/5\\
  2/5 & 0 & 0 & 1 & 0 & 3/5\\
  1/5 & 0 & 0 & 0 & 1 & 4/5\\
  0 & 0 & 0 & 0 & 0 & 1\\
\end{bmatrix}.
\]
The change of variables is made by replacing $K^i$ with
\[
K_c^i = 
\begin{bmatrix}
  \cL \phi^i(x_I^i)\\ \cB \phi^i(x_B^i)\\ (T^i)^{-1} \cC \phi^i(x_\Gamma^i)
\end{bmatrix}.
\]
This change of variables in conjunction with the modification to $A$ described below has been experimentally shown to improve convergence, especially of the Neumann-Neumann acceleration proposed in \cref{subsec:nn-accel}.
We replace
\[
A^i=\begin{bmatrix}
  \cF \phi^i(x_\Pi) \\ \cF \phi^i(x_\Delta)
\end{bmatrix}
\text{ with }
A_m^i=\begin{bmatrix}
  \cF \phi^i(x_\mathcal{E}) \\ \cF \phi^i(x_\Delta)
\end{bmatrix},
\]
where $\cF \phi^i(x_\mathcal{E})$ represents the mean flux on edges, which is taken by excluding the opposing corner to preserve the degrees of freedom.
With \cref{fig:change}, $A_m^i$ becomes
\[
\begin{bmatrix}
  \frac{1}{5}\sum_{j=1}^5 \cF \phi^i(x_j)\\
  \cF \phi^i(x_2)\\
  \cF \phi^i(x_3)\\
  \cF \phi^i(x_4)\\
  \cF \phi^i(x_5)\\
  \frac{1}{5}\sum_{j=2}^6 \cF \phi^i(x_j)
\end{bmatrix}.
\]

\section{Numerical Experiments}
\label{Sec:Num}
We present the results of numerical experiments using the Neumann-Neumann acceleration on several PDEs.
Performance is assessed through a mixture of accuracy, CG iterations, and wall-clock time.
Accuracy is evaluated using relative $L^2$ and $H^1$ errors, defined as
\begin{align*}
  \text{Relative } L^2 \text{ error} &= \frac{\norm{u - u^*}_{L^2(\Omega)}}{\norm{u}_{L^2(\Omega)}}, \\
  \text{Relative } H^1 \text{ error} &= \frac{\norm{u - u^*}_{H^1(\Omega)}}{\norm{u}_{H^1(\Omega)}}.
\end{align*}
Here, $u$ is the reference solution and $u^*$ is the solution computed using the method being evaluated.

Each local neural network is a single hidden layer multi layer perceptron with $2^{16}/N$ neurons unless stated otherwise.
The activation function is the hyperbolic tangent function.
The training points on each subdomain are given by a $320/\sqrt{N}\times 320/\sqrt{N}$ uniform grid unless stated otherwise.
In all experiments, the auxiliary variables are computed using CG, stopping when the relative residual becomes less than $10^{-9}$.

When the exact solution is not known, reference solutions are computed with the finite element method using FEniCSx~\cite{baratta2023dolfinx}.
Order 16 Lagrangian elements on a $32\times32$ quadrilateral mesh are used.
A direct solver is used to solve the systems.

DDELM is implemented with Pytorch~\cite{paszke2019pytorch} and mpi4py~\cite{dalcin2021mpi}.
Except for the $2\times2$ case, all experiments were conducted on eight nodes with 192GB RAM:
six nodes equipped with two Intel Xeon Gold 6148 (2.4GHz 20c) processors; two nodes with two Intel Xeon Gold 6248 (2.5GHz 20c) processors.
Each subdomain was allocated one process.
In turn, each process was allocated $256/N$ cores and $N/8$ processes were spawned per node.
For the $2\times2$ case, four nodes with two Intel Xeon Gold 6148 (2.4GHz 20c) processors were used;
each process was allocated 20 cores and 1 process was spawned per node.
Note that although 40 cores were available per node, computation was faster when each process was limited to one CPU.
The cluster employs a 100Gbps Infiniband network.

\subsection{Initialization of parameters of the hidden layer}
We utilize the initialization strategy presented in our previous work for a single hidden layer~\cite{lee2024nonoverlapping}:
\begin{equation*}
  \label{eqn:init_be}
  \begin{split}
    w^{i}_{j} &\sim U([-l, l]^d), \text{ for } l=c\frac{(M^i)^{1/d}}{\mathrm{diam}(\bar{B}_r(\Omega^i))}\\
    b^{i}_{j} &= -w^{i}_{j} \cdot \xi^{i}_{j} \quad \text{ for } \xi^{i}_{j} \sim U(\bar{B}_{r}(\Omega^i)).
  \end{split}
\end{equation*}
Here, $U$ is the uniform distribution, $d$ is the dimension of the problem, $c$ is a proportionality constant, and $\bar{B}_{r}(\Omega^i)$ is the closed $r$-neighborhood about $\Omega^i$.
In the following experiments, $r$ will be $\mathrm{diam}(\Omega^i)/2$.
For the Poisson equation, $c=3/16$ is used, which makes $l=32$; for the biharmonic equation $c=3/32$ is used, which makes $l=64$.

\subsection{Ablation Study}
\label{subsec:ablation}
We investigate the effectiveness of the change of variables and modification of $A$ proposed in \cref{subsec:chov} and mixed equation form of the Neumann-Neumann acceleration proposed in \cref{subsec:nn-accel}.
Performance is evaluated using the 2-D Poisson problem
\begin{equation}
  \label{eqn:poisson}
  \left\{\begin{aligned}
    -\Delta u &= f &&\text{ in } \Omega=(0,1)^{2}, \\
    u &= g &&\text{ on } \partial \Omega,
  \end{aligned}\right.
\end{equation}
where we specify the exact solution as $u(x,y)=\sin(\pi x)\sin(\pi y)$.

\subsubsection{Change of Variables}
We first compare using a coarse space with the change of variables in \cref{subsec:chov} to using the standard coarse space consisting of corners points.
\begin{table}
  \caption{
    Relative $H^1$ error and number of CG iterations alternating between using $K$, $K_c$ and $A$, $A_m$ in solving the Poisson equation with exact solution $u(x,y)=\sin(\pi x)\sin(\pi y)$.
    A checkmark in the NN column indicates a solution obtained with Neumann-Neumann acceleration via \eqref{eqn:nn} whereas no checkmark indicates a solution obtained by \eqref{eqn:coarse_sol}.
    }
    \label{tab:chov}
  \centering
  \begin{tabular}{cccccccccc}
    \toprule
    Compo- & \multirow{2}{*}{NN} & \multicolumn{2}{c}{$2\times2$} & \multicolumn{2}{c}{$4\times4$} & \multicolumn{2}{c}{$8\times8$} & \multicolumn{2}{c}{$16\times16$}\\
    nents &&  $H^1$ error & Iter & $H^1$ error & Iter & $H^1$ error & Iter & $H^1$ error & Iter\\
    \midrule \midrule
    \multirow{2}{*}{$K$, $A$} &
    & 1.08e-09 & 826 & 1.78e-09 & 2,914 & 6.95e-09 & 5,345 & 4.37e-08 & 10,597\\
    & \checkmark
    & 5.44e-01 & 67 & 2.49e-08 & 592 & 1.55e-08 & 7,529 & 8.01e-01 & 27,360\\
    \midrule
    \multirow{2}{*}{$K_c$, $A$} &
    & 4.74e-08 & 599 & 1.92e-08 & 2,784 & 7.12e-09 & 5,198 & 1.40e-08 & 7,491\\
    & \checkmark
    & 7.62e-01 & 63 & 5.18e-06 & 2,349 & 1.40e-07 & 7,998 & 6.08e-08 & 15,763\\
    \midrule
    \multirow{2}{*}{$K$, $A_m$} &
    & 1.17e-07 & 588 & 1.13e-07 & 2,078 & 1.33e-07 & 3,877 & 1.27e-07 & 6,167\\
    & \checkmark
    & 7.81e-09 & 40 & 1.32e-07 & 155 & 5.41e-08 & 563 & 1.16e-07 & 4,219\\
    \midrule
    \multirow{2}{*}{$K_c$, $A_m$} &
    & 1.21e-09 & 682 & 7.96e-10 & 1,948 & 2.03e-09 & 2,802 & 1.41e-08 & 3,874\\
    & \checkmark
    & 2.82e-09 & 45 & 3.50e-08 & 284 & 3.14e-09 & 592 & 1.39e-08 & 1,138\\
  \bottomrule
  \end{tabular}
\end{table}
In \cref{tab:chov}, we see that adpoting $K_c$ with the change of variables gives better results than $K$.
Moreover, without the change of variables and modifcation of $A$, the Neumann-Neumann acceleration shows very unstable results.
Using $A_m$ seems to actually allow Neumann-Neumann acceleration to improve convergence.
Combining $K_c$ and $A_m$ gives the best results in terms of accuracy.
While using only $A_m$ can result in slightly fewer iterations up to the $8\times 8$ case, iteration increase is not consistent as the number of subdomains increases, resulting in considerably more iterations for the $16\times 16$ case.

In light of the superior accuracy and stable iteration counts, we will use both $K_c$ and $A_m$ for the coarse space throughout the following numerical experiment section.

\subsubsection{Mixing in Neumann-Neumann Acceleration}
We next consider the effect of $\theta$ in the mixed form \eqref{eqn:nn-mixed} of the proposed Neumann-Neumann acceleration.
At $\theta=0$, the method reduces to the original coarse space formulation.
The results are summarized in \cref{tab:nn-mixed}.
\begin{table}
  \caption{
    Relative $H^1$ error and number of CG iterations for different values of $\theta$ in \eqref{eqn:nn-mixed} solving the Poisson equation with exact solution $u(x,y)=\sin(\pi x)\sin(\pi y)$.
    }
    \label{tab:nn-mixed}
\centering
\begin{tabular}{ccccccccc}
  \toprule
  \multirow{2}{*}{$\theta$} & \multicolumn{2}{c}{$2\times2$} & \multicolumn{2}{c}{$4\times4$} & \multicolumn{2}{c}{$8\times8$} & \multicolumn{2}{c}{$16\times16$}\\
  &  $H^1$ error & Iter & $H^1$ error & Iter & $H^1$ error & Iter & $H^1$ error & Iter\\
  \midrule \midrule
  0 & 1.21e-09 & 682 & 7.96e-10 & 1,948 & 2.03e-09 & 2,802 & 1.41e-08 & 3,874 \\
  0.5 & 1.11e-09 & 564 & 6.52e-10 & 1,599 & 1.98e-09 & 2,167 & 1.40e-08 & 2,938 \\
  0.9 & 1.08e-09 & 375 & 6.12e-10 & 906 & 1.92e-09 & 1,176 & 1.39e-08 & 1,634 \\
  0.99 & 1.02e-09 & 201 & 5.95e-10 & 368 & 1.89e-09 & 503 & 1.39e-08 & 991 \\
  0.999 & 1.07e-09 & 92 & 7.23e-10 & 189 & 1.90e-09 & 328 & 1.38e-08 & 919 \\
  0.9999 & 2.48e-09 & 50 & 9.44e-10 & 209 & 2.10e-09 & 337 & 1.37e-08 & 906 \\
  1 & 2.82e-09 & 45 & 3.50e-08 & 284 & 3.14e-09 & 592 & 1.39e-08 & 1,138 \\
  \bottomrule
\end{tabular}
\end{table}

We can see that a suitable $\theta$ value close to 1 can reduce the number of iterations without sacrificing accuracy even compared to the $\theta=1$ case.
Optimal values of $\theta$ seem to increase with the number of subdomains.
Nevertheless, the method is robust to the choice of $\theta$ when $\theta$ is very close to 1.

In the remaining experiments, when the Neumann-Neumann acceleration is used, the mixed form with $\theta=.999$ will be used.

\subsection{Numerical Experiment Results}
In this section, we present the results of numerical experiments using the Neumann-Neumann acceleration with coarse space on several PDEs.
As mentioned in \cref{subsec:ablation}, we will use both $K_c$ and $A_m$ from \cref{subsec:chov} for the coarse space.
When using the Neumann-Neumann acceleration, we refer to the mixed form~\eqref{eqn:nn-mixed} with $\theta=.999$.
We denote the vanilla DDELM method from \cref{subsubsec:ddelm} by DDELM, DDELM with coarse space by DDELM-CS, and Neumann-Neumann accelerated DDELM-CS by DDELM-NN.
Reported computation times are averages of 10 runs.

\subsubsection{Poisson Equation}
Consider the 2-D Poisson problem~\eqref{eqn:poisson}.
First consider the case when the exact solution is
\[
u(x,y)=\sin(2\pi x)e^y.
\]
The results are shown in \cref{tab:poisson1}.
\begin{table}
  \caption{
    Relative $L^2$ and $H^1$ errors, number of CG iterations, and wall-clock time solving the Poisson equation with exact solution $u(x,y)=\sin(2\pi x)e^y$ using a total of 65,536 neurons.
    $N$ is the number of subdomains.
    }
  \label{tab:poisson1}
\centering
\begin{tabular}{ccccccccc}
  \toprule
  Method & $N$ & $L^2$ error & $H^1$ error & Iterations & Wall-clock time (sec)\\
  \midrule \midrule
  \multirow{4}{*}{DDELM}
  & $2\times2$ & 1.11e-10 & 1.31e-09 & 1,016 & 460.82 \\
  & $4\times4$ & 4.31e-09 & 5.81e-08 & 2,634 & 11.49 \\
  & $8\times8$ & 7.38e-08 & 8.65e-08 & 4,352 & 6.27 \\
  & $16\times16$ & 2.72e-07 & 2.13e-07 & 6,757 & 11.36 \\
  \midrule
  \multirow{4}{*}{DDELM-CS}
  & $2\times2$ & 1.09e-10 & 1.32e-09 & 615 & 461.70 \\
  & $4\times4$ & 3.36e-09 & 4.36e-08 & 1,640 & 10.61 \\
  & $8\times8$ & 7.92e-09 & 7.31e-08 & 2,369 & 4.73 \\
  & $16\times16$ & 5.11e-09 & 8.94e-08 & 2,962 & 5.59 \\
  \midrule
  \multirow{4}{*}{DDELM-NN}
  & $2\times2$ & 1.34e-09 & 8.24e-09 & 87 & 909.24 \\
  & $4\times4$ & 1.22e-09 & 4.56e-09 & 194 & 15.93 \\
  & $8\times8$ & 1.44e-09 & 1.03e-08 & 312 & 1.53 \\
  & $16\times16$ & 2.19e-09 & 3.00e-08 & 618 & 1.58 \\
  \bottomrule
\end{tabular}
\end{table}

\begin{figure}
  \begin{subfigure}{0.32\textwidth}
    \centering
    \includegraphics[width=\linewidth]{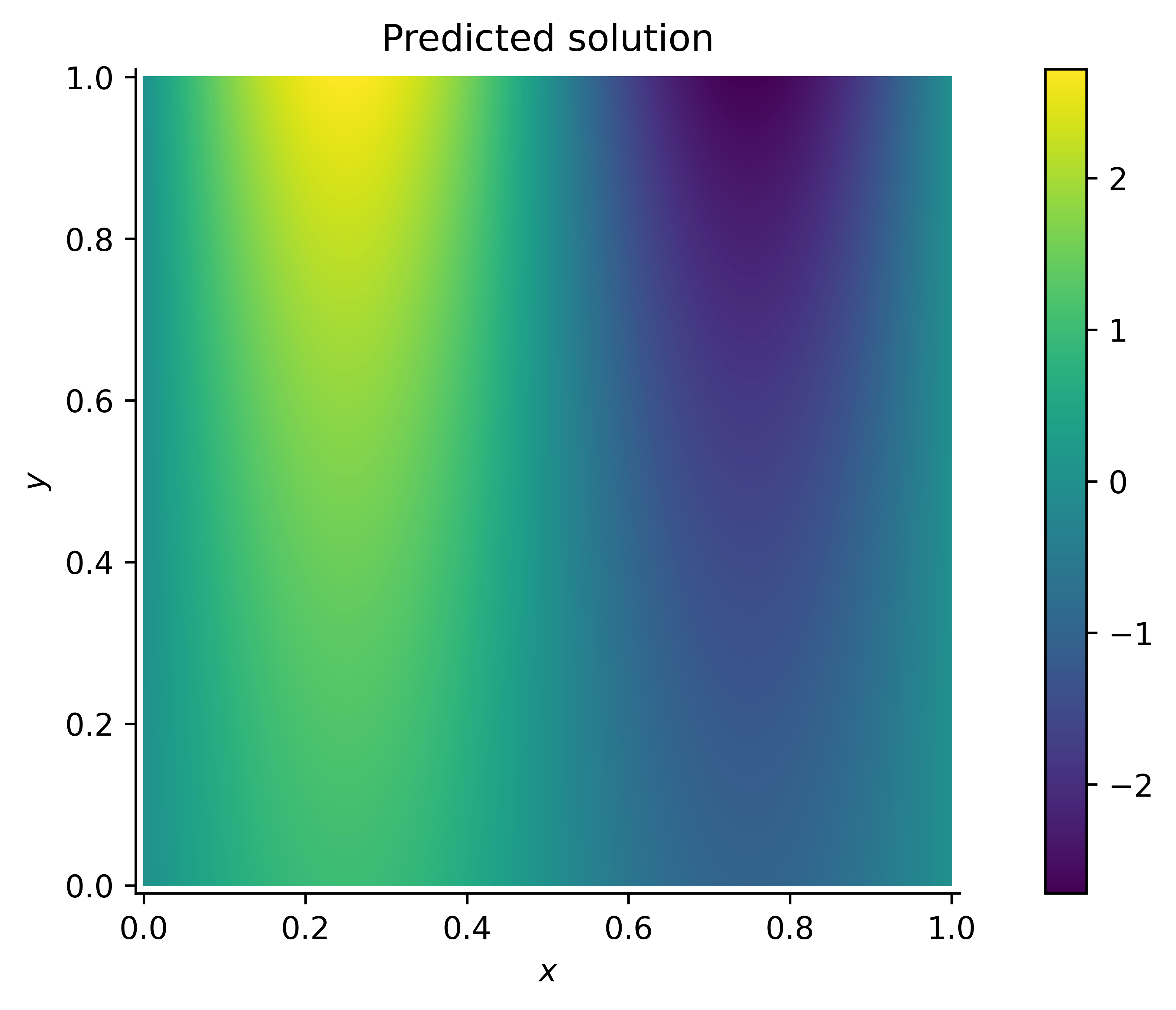}
    \caption*{(1a)}
  \end{subfigure}
  \hfill
  \begin{subfigure}{0.32\textwidth}
    \centering
    \includegraphics[width=.965\linewidth]{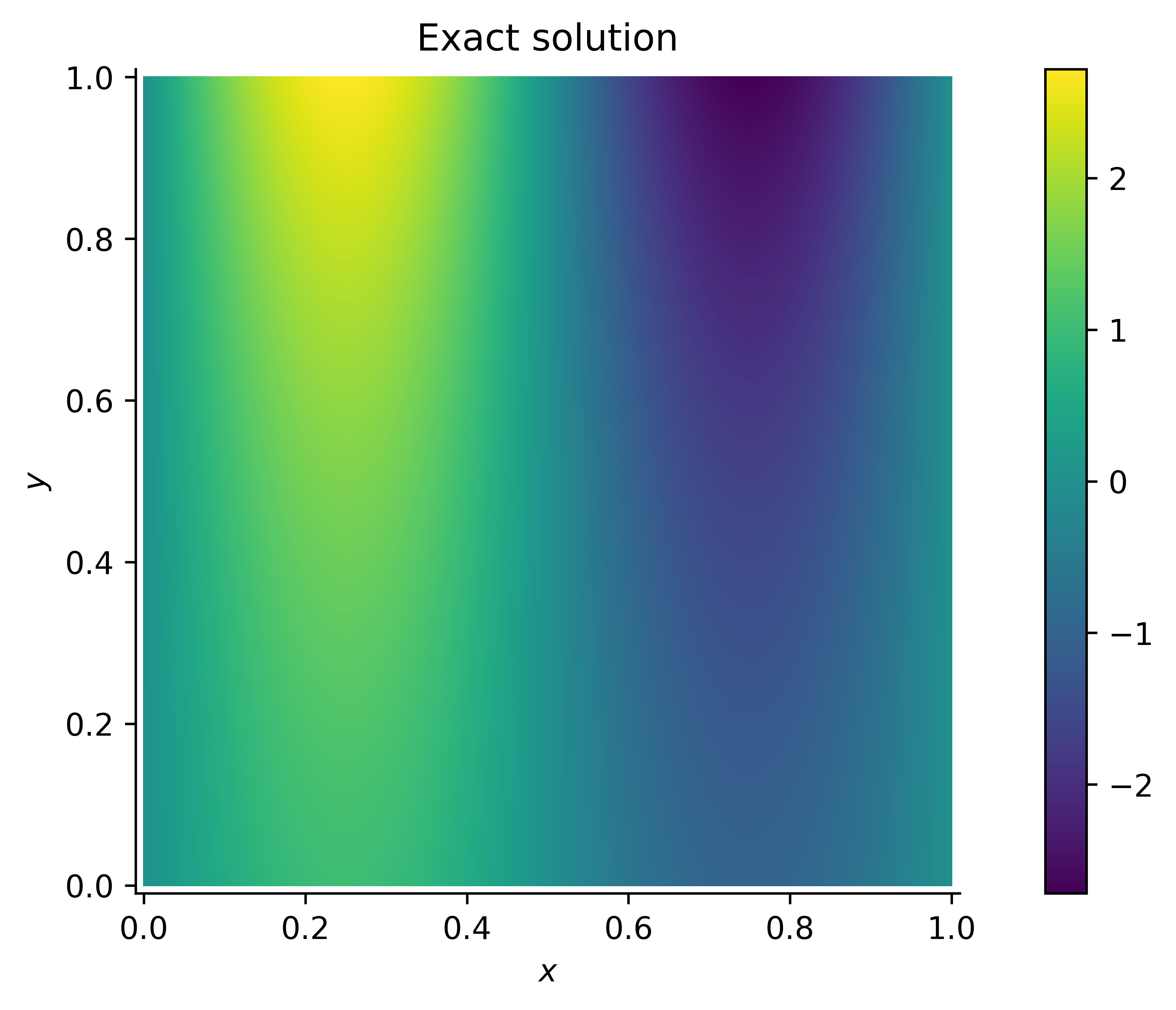}
    \caption*{(1b)}
  \end{subfigure}
  \begin{subfigure}{0.32\textwidth}
    \centering
    \includegraphics[width=.965\linewidth]{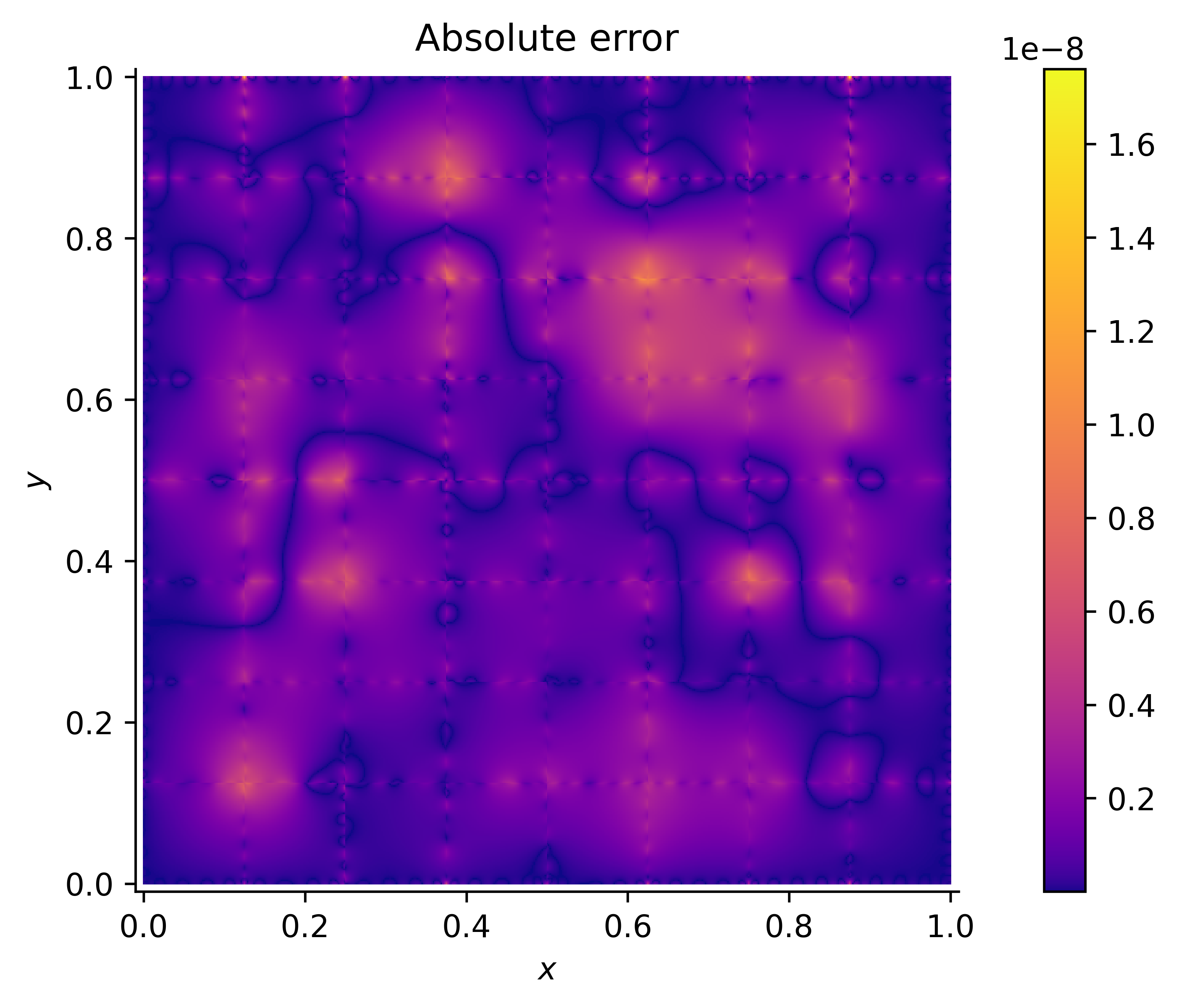}
    \caption*{(1c)}
  \end{subfigure}

  \begin{subfigure}{0.32\textwidth}
    \centering
    \includegraphics[width=\linewidth]{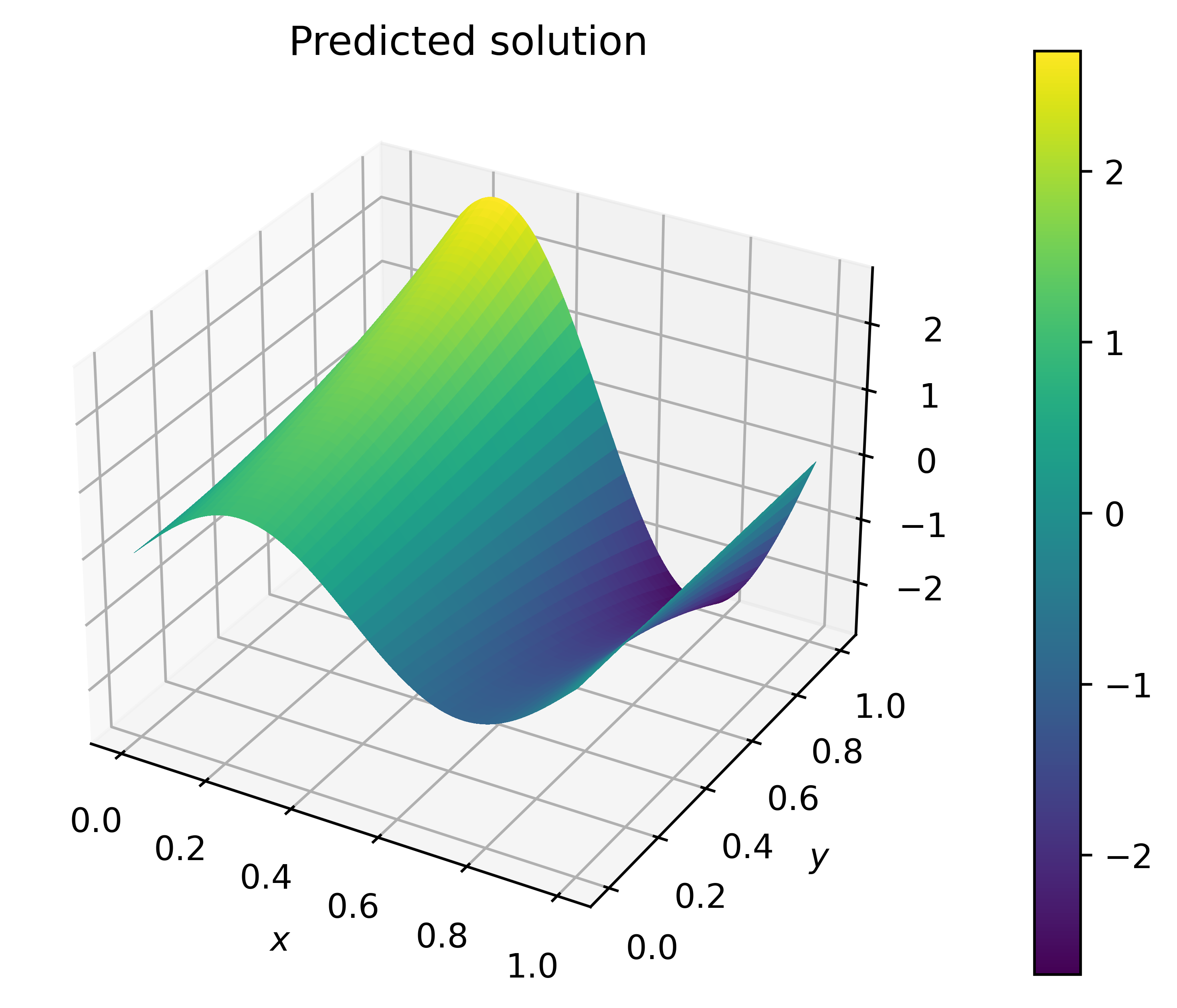}
    \caption*{(2a)}
  \end{subfigure}
  \hfill
  \begin{subfigure}{0.32\textwidth}
    \centering
    \includegraphics[width=.965\linewidth]{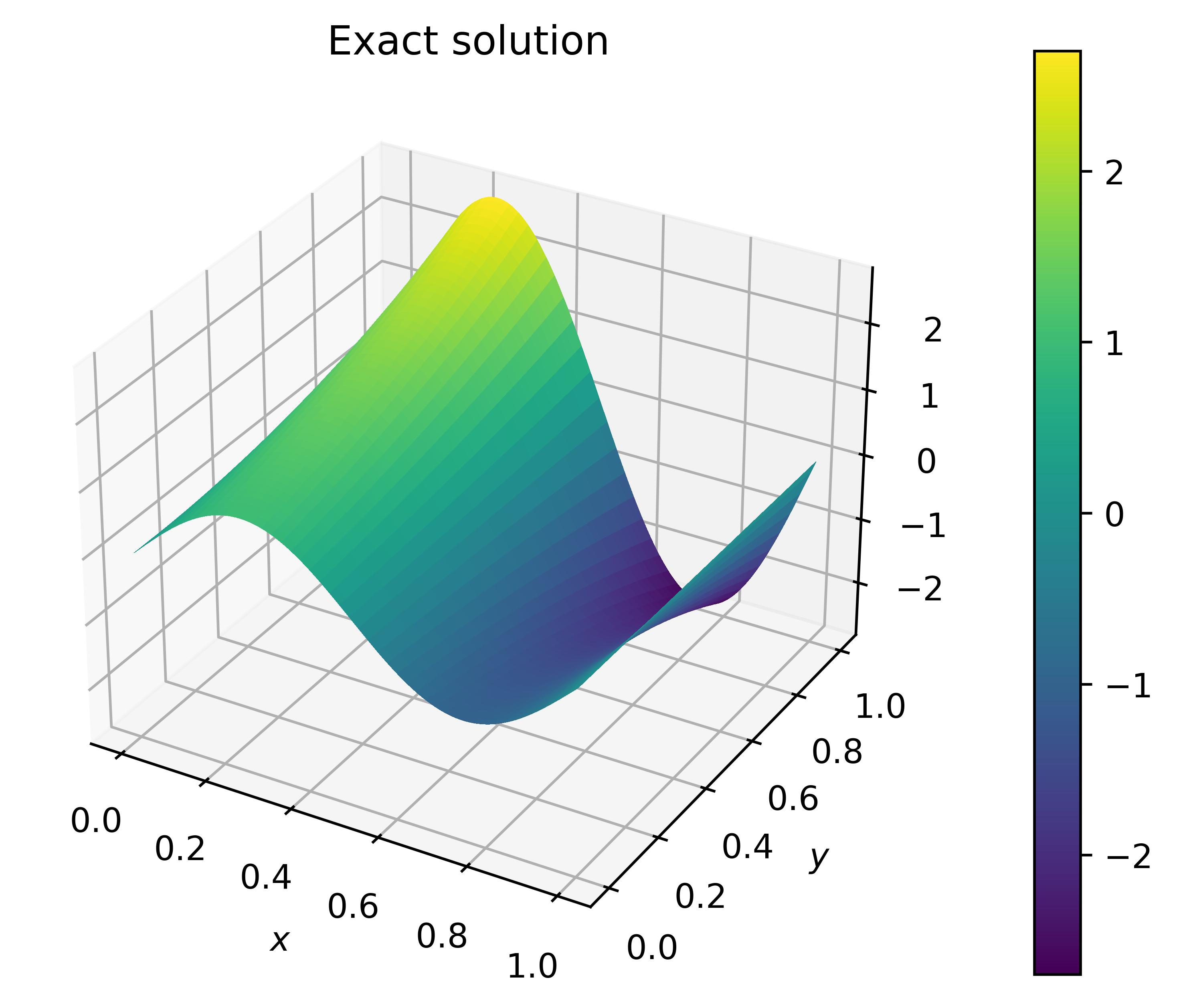}
    \caption*{(2b)}
  \end{subfigure}
  \hfill
  \begin{subfigure}{0.32\textwidth}
    \centering
    \includegraphics[width=.965\linewidth]{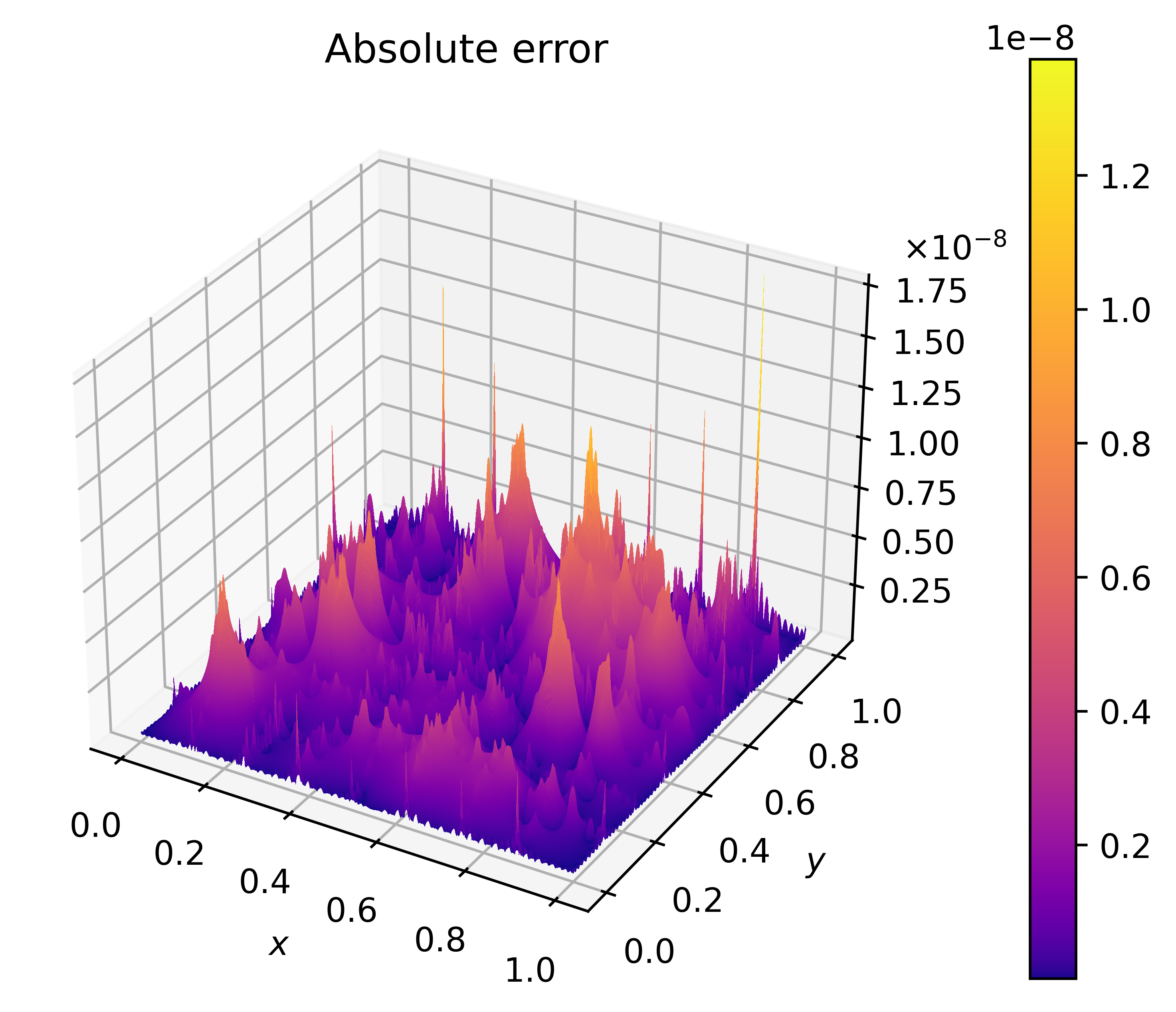}
    \caption*{(2c)}
  \end{subfigure}

  \begin{subfigure}{0.495\textwidth}
    \centering
    \includegraphics[width=\linewidth]{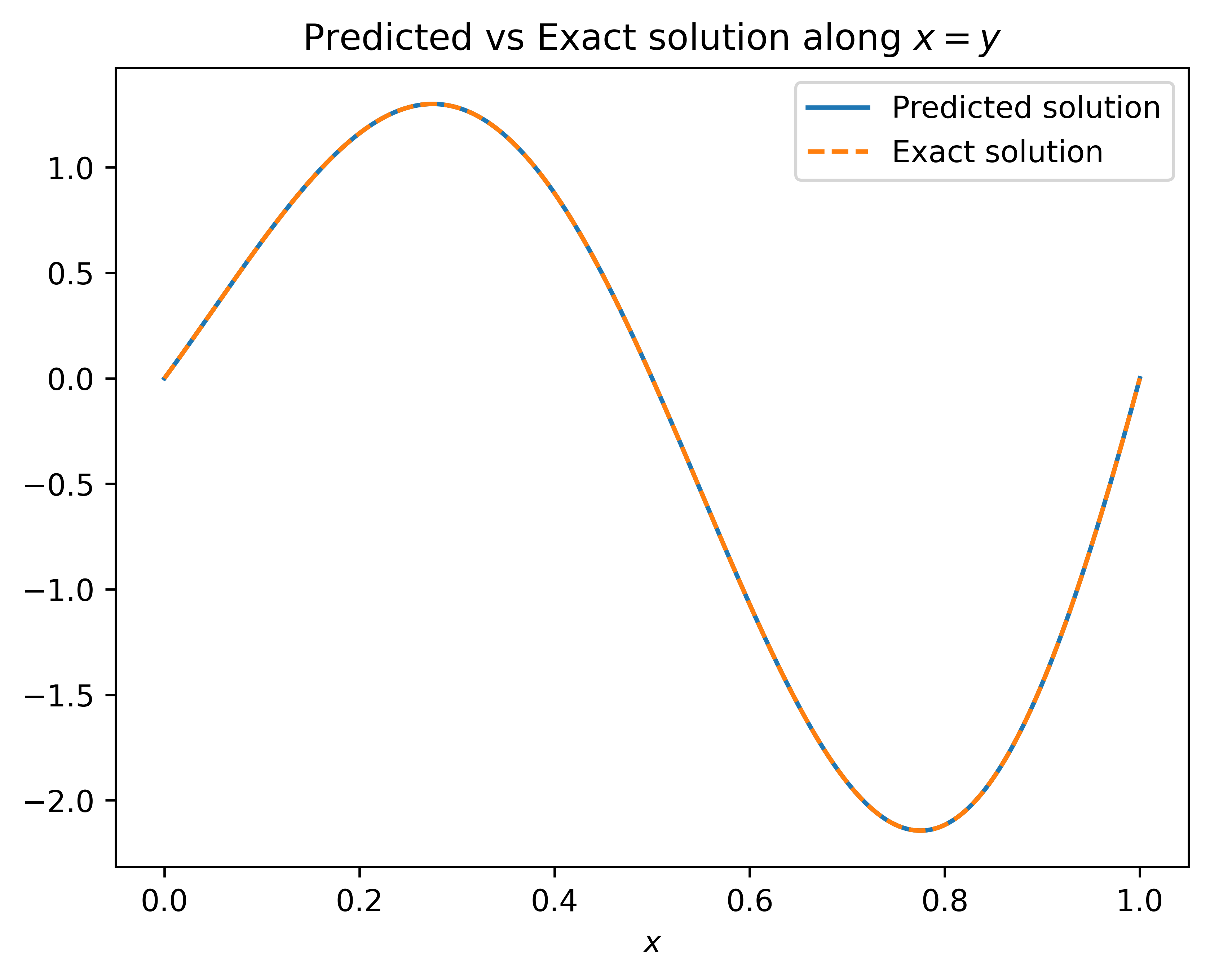}
    \caption*{(3a)}
  \end{subfigure}
  \hfill
  \begin{subfigure}{0.495\textwidth}
    \centering
    \includegraphics[width=.965\linewidth]{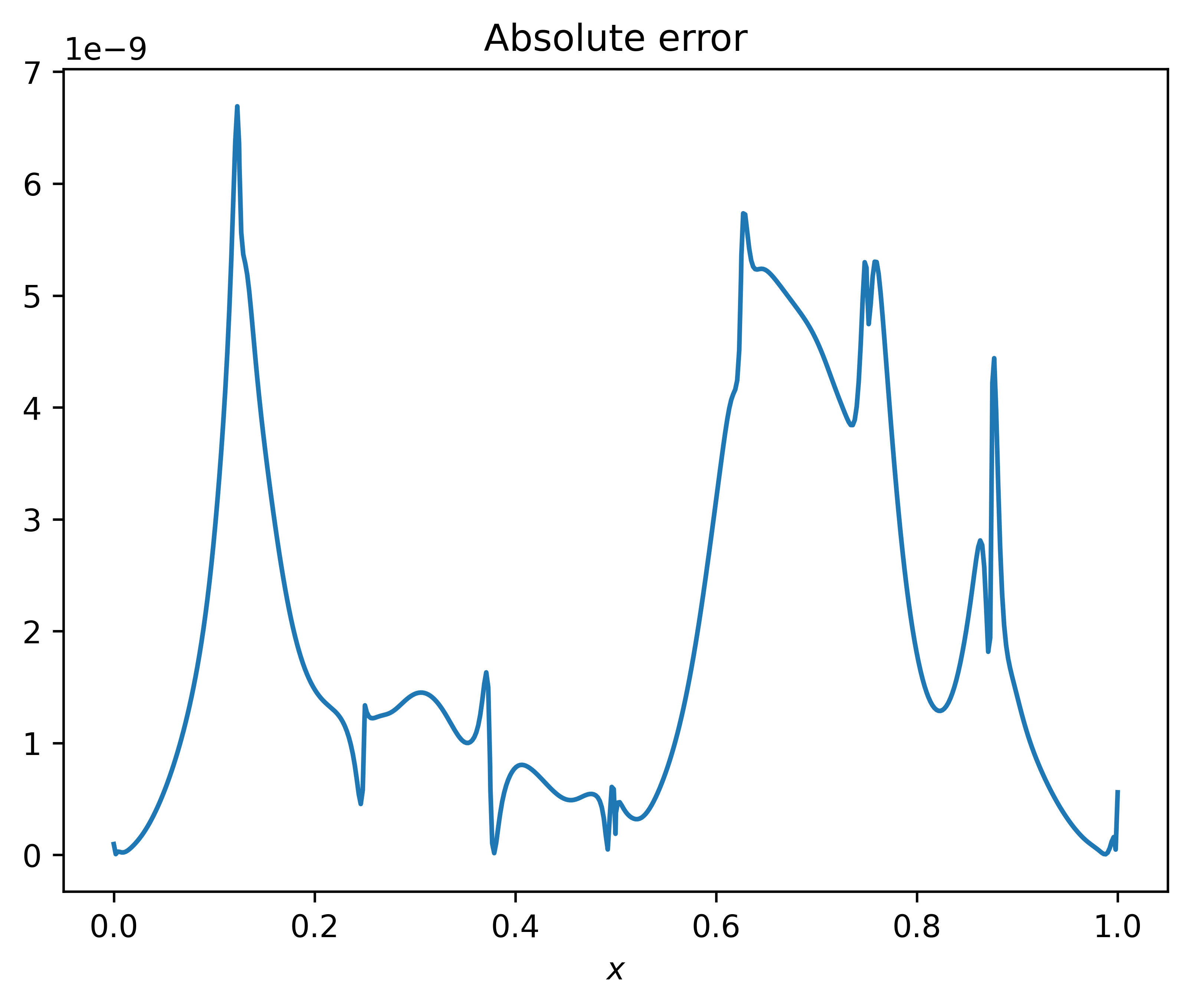}
    \caption*{(3b)}
  \end{subfigure}
  \caption{Solution and error of 2-D Poisson equation with exact solution $u(x, y)=\sin(2\pi x)e^y$.
  The DDELM-NN solution is computed on $8\times8$ subdomains.
  (1a) DDELM-NN solution 2D plot, (1b) Exact solution 2D plot, (1c) Absolute error of DDELM-NN solution 2D plot.
  (2a) DDELM-NN solution 3D plot, (2b) Exact solution 3D plot, (2c) Absolute error of DDELM-NN solution 3D plot.
  (3a) DDELM-NN and exact solution along the line $x=y$, (3b) Absolute error of DDELM-NN solution along the line $x=y$.}
  \label{fig:poisson_sin2piexp}
\end{figure}

We note that DDELM-NN is slower in the $2\times2$ and $4\times4$ cases.
This is because wall-clock time in these cases is dominated by factorization for least squares problems of which the Neumann-Neumann acceleration must perform twice: once for $K$ and once for $\tilde{K}$.
Compared to DDELM, DDELM-NN reduces the number of iterations fourteen-fold to achieve a four-fold speedup in the fastest $8\times8$ case.
While DDELM-CS reduces the number of iterations and wall-clock time, the reduction is not as significant as with DDELM-NN.
DDELM-NN also gives the lowest relative $L^2$ and $H^1$ errors except in the $2\times2$ case.
The $8\times8$ subdomain DDELM-NN solution and error is plotted in~\cref{fig:poisson_sin2piexp}.
We note excellent agreement between the DDELM-NN solution and the exact solution.
The magnitude of error is quite small even where it peaks at the subdomain corners.

We also investigate a more oscillatory case where $g\equiv 0$ and $f$ is an approximate Gaussian random field~(GRF),
\begin{equation}
  \label{eqn:grf}
\begin{gathered}
  f(\mathbf{x}) = \sum_{i=1}^{256} a_i \sin(w_i\cdot \mathbf{x} + b_i),\\
  a_i \sim \mathcal{N}(0, \alpha^4/256), \quad w_i \sim \mathcal{N}(0, \alpha^2I), \quad b_i \sim U((0, 2\pi)),
\end{gathered}
\end{equation}
where a larger $\alpha$ corresponds to a more oscillatory field.
In this case, the exact solution is not known and the FEM solution is used as reference.
The results for $\alpha=32$ are shown in \cref{tab:poisson2}.
\begin{table}
  \caption{
    Relative $L^2$ and $H^1$ errors, number of CG iterations, and wall-clock time solving the Poisson equation with $g\equiv 0$ and $f$ an approximate GRF as in \eqref{eqn:grf} with $\alpha=32$ using a total of 65,536 neurons.
    $N$ is the number of subdomains.
    }
  \label{tab:poisson2}
\centering
\begin{tabular}{ccccccccc}
  \toprule
  Method & $N$ & $L^2$ error & $H^1$ error & Iterations & Wall-clock time (sec)\\
  \midrule \midrule
  \multirow{4}{*}{DDELM}
  & $2\times2$ & 2.19e-06 & 3.16e-05 & 1,028 & 458.21 \\
  & $4\times4$ & 1.85e-06 & 2.73e-05 & 3,337 & 12.44 \\
  & $8\times8$ & 1.88e-06 & 2.82e-05 & 5,558 & 7.99 \\
  & $16\times16$ & 1.95e-06 & 3.08e-05 & 7,881 & 13.25 \\
  \midrule
  \multirow{4}{*}{DDELM-CS}
  & $2\times2$ & 2.19e-06 & 3.16e-05 & 595 & 460.01 \\
  & $4\times4$ & 1.85e-06 & 2.73e-05 & 1,866 & 10.97 \\
  & $8\times8$ & 1.88e-06 & 2.82e-05 & 2,658 & 5.24 \\
  & $16\times16$ & 1.92e-06 & 3.02e-05 & 2,977 & 5.66 \\
  \midrule
  \multirow{4}{*}{DDELM-NN}
  & $2\times2$ & 2.19e-06 & 3.16e-05 & 97 & 916.46 \\
  & $4\times4$ & 1.85e-06 & 2.72e-05 & 196 & 16.03 \\
  & $8\times8$ & 1.85e-06 & 2.78e-05 & 310 & 1.52 \\
  & $16\times16$ & 1.91e-06 & 2.97e-05 & 590 & 1.52 \\
  \bottomrule
\end{tabular}
\end{table}
\begin{figure}
  \begin{subfigure}{0.32\textwidth}
    \centering
    \includegraphics[width=\linewidth]{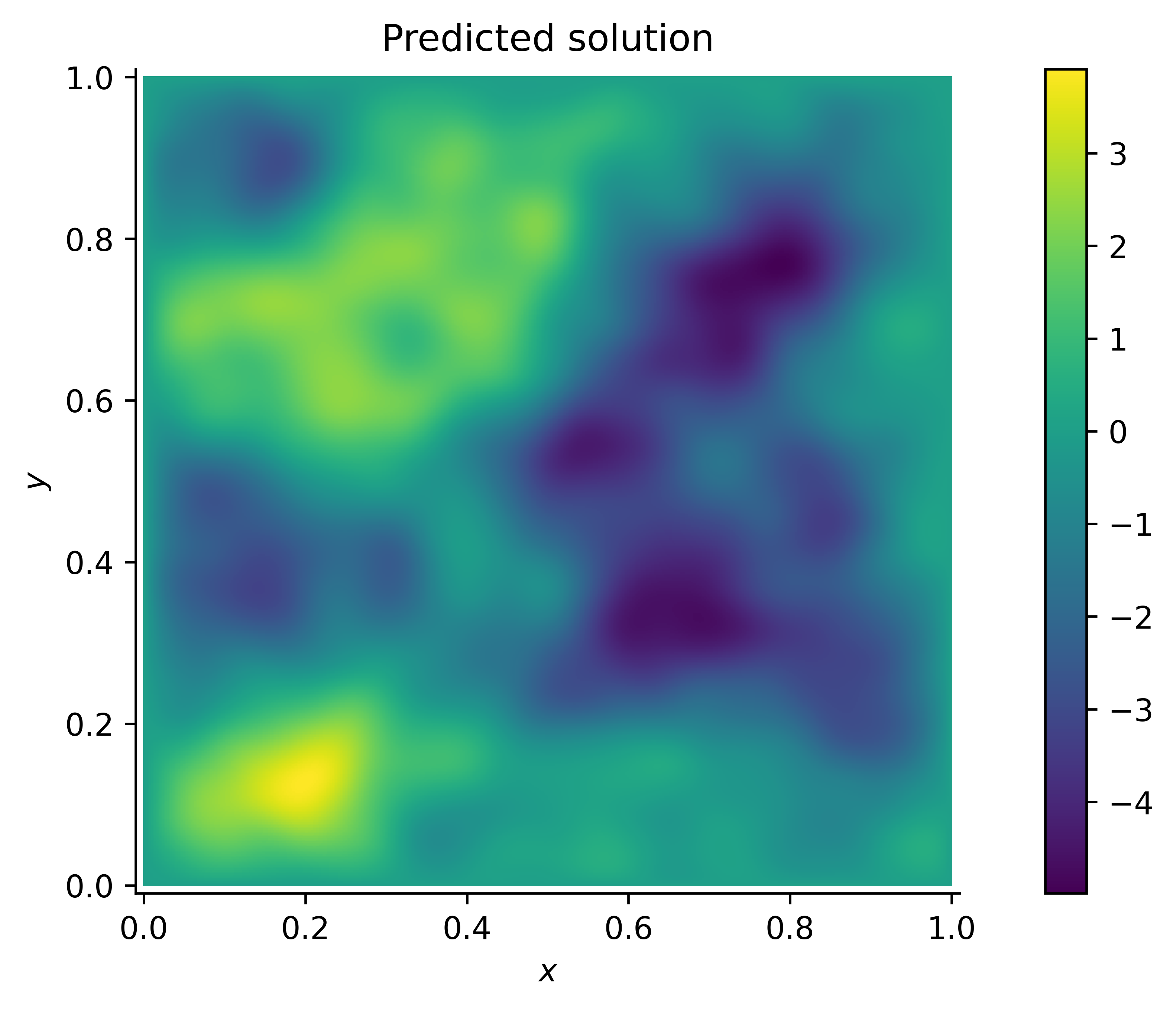}
    \caption*{(1a)}
  \end{subfigure}
  \hfill
  \begin{subfigure}{0.32\textwidth}
    \centering
    \includegraphics[width=.965\linewidth]{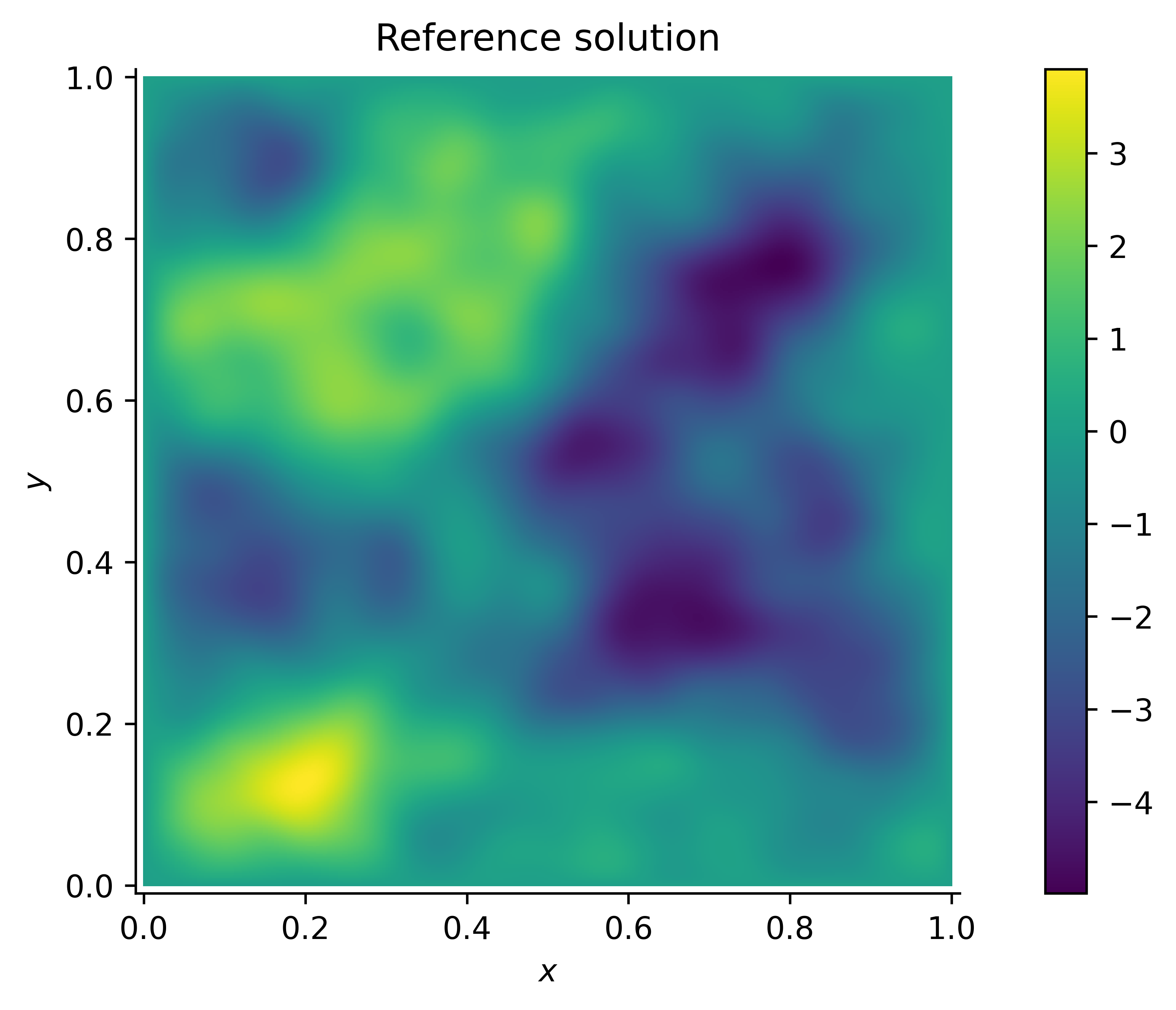}
    \caption*{(1b)}
  \end{subfigure}
  \begin{subfigure}{0.32\textwidth}
    \centering
    \includegraphics[width=.965\linewidth]{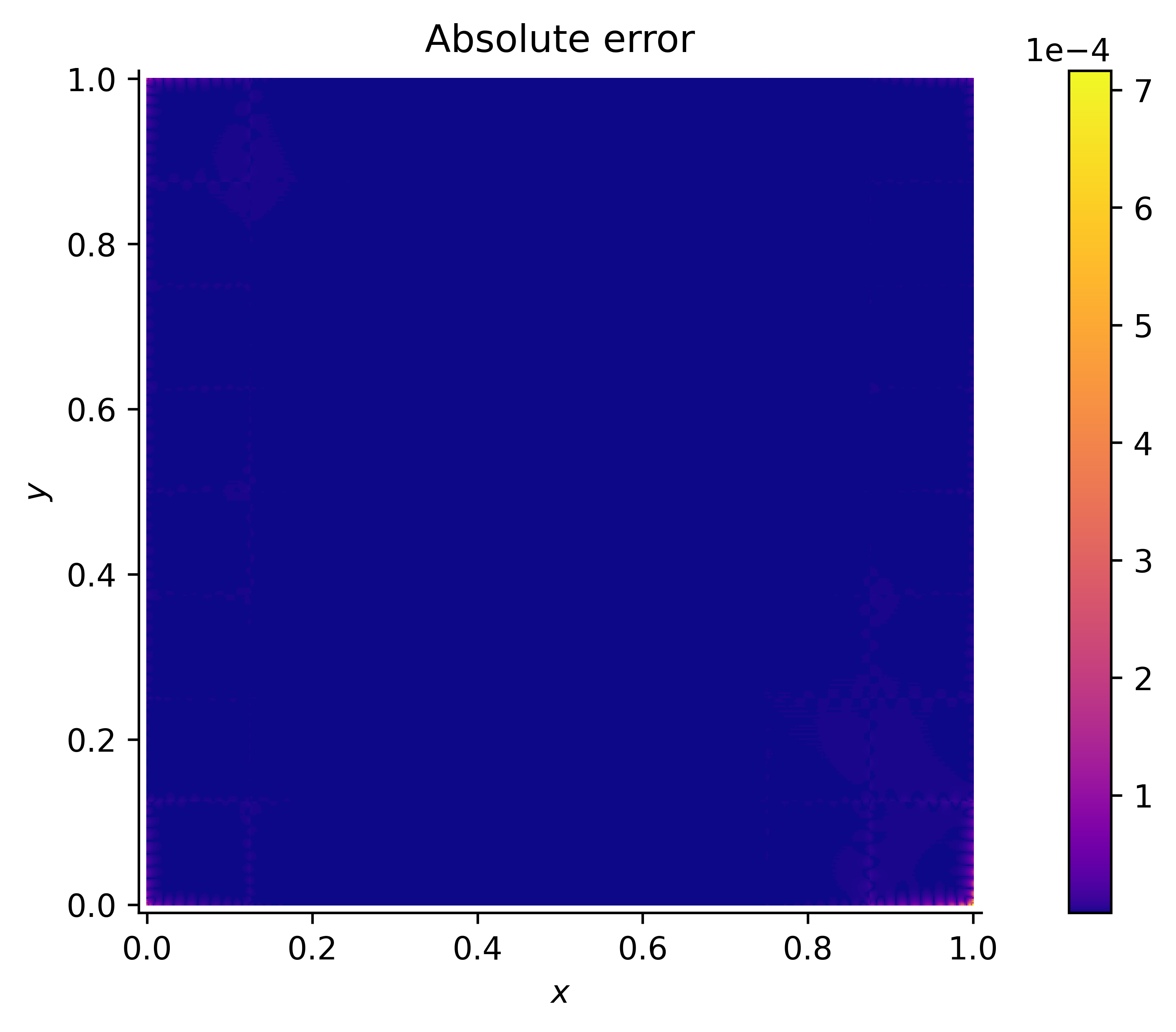}
    \caption*{(1c)}
  \end{subfigure}

  \begin{subfigure}{0.32\textwidth}
    \centering
    \includegraphics[width=\linewidth]{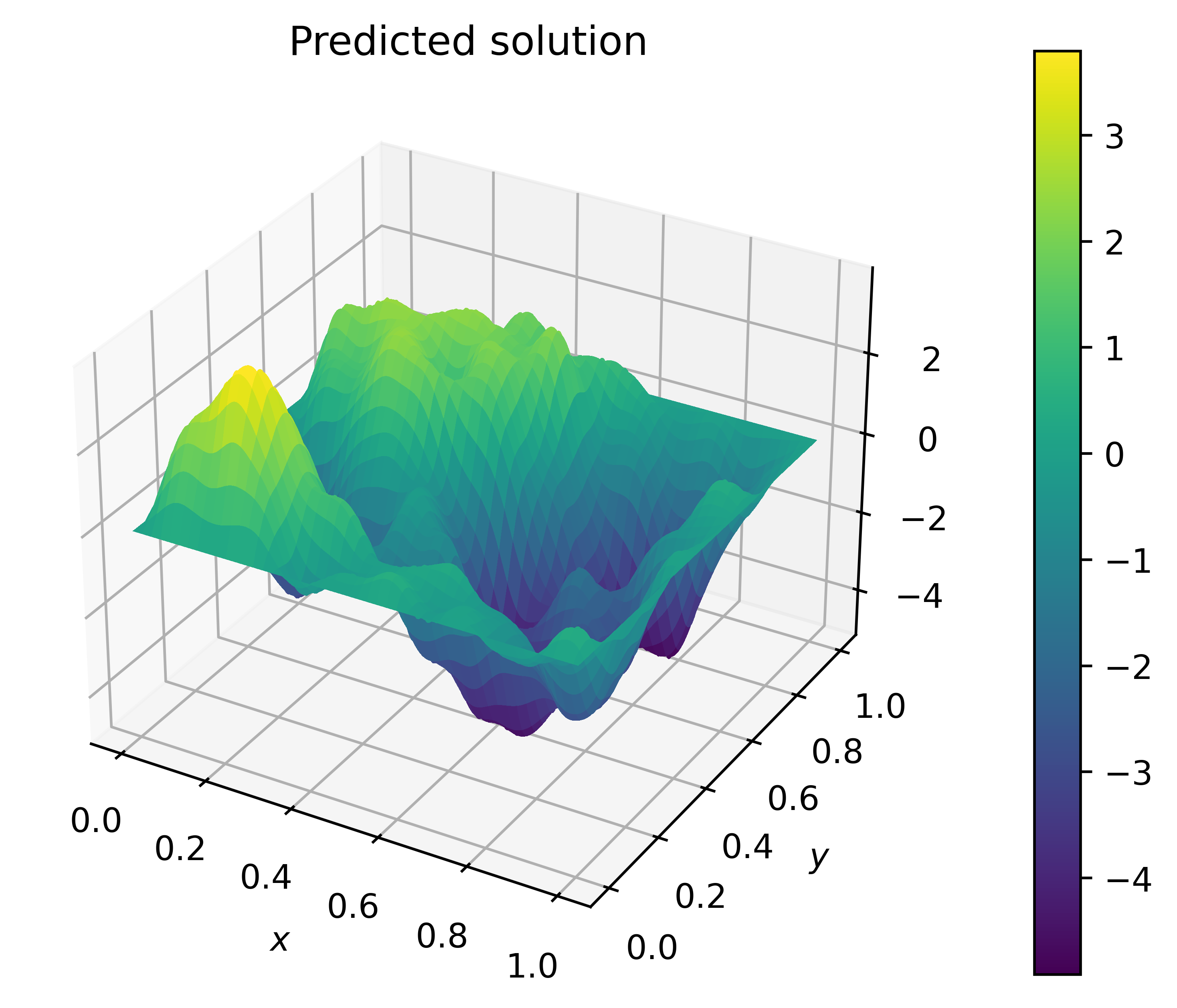}
    \caption*{(2a)}
  \end{subfigure}
  \hfill
  \begin{subfigure}{0.32\textwidth}
    \centering
    \includegraphics[width=.965\linewidth]{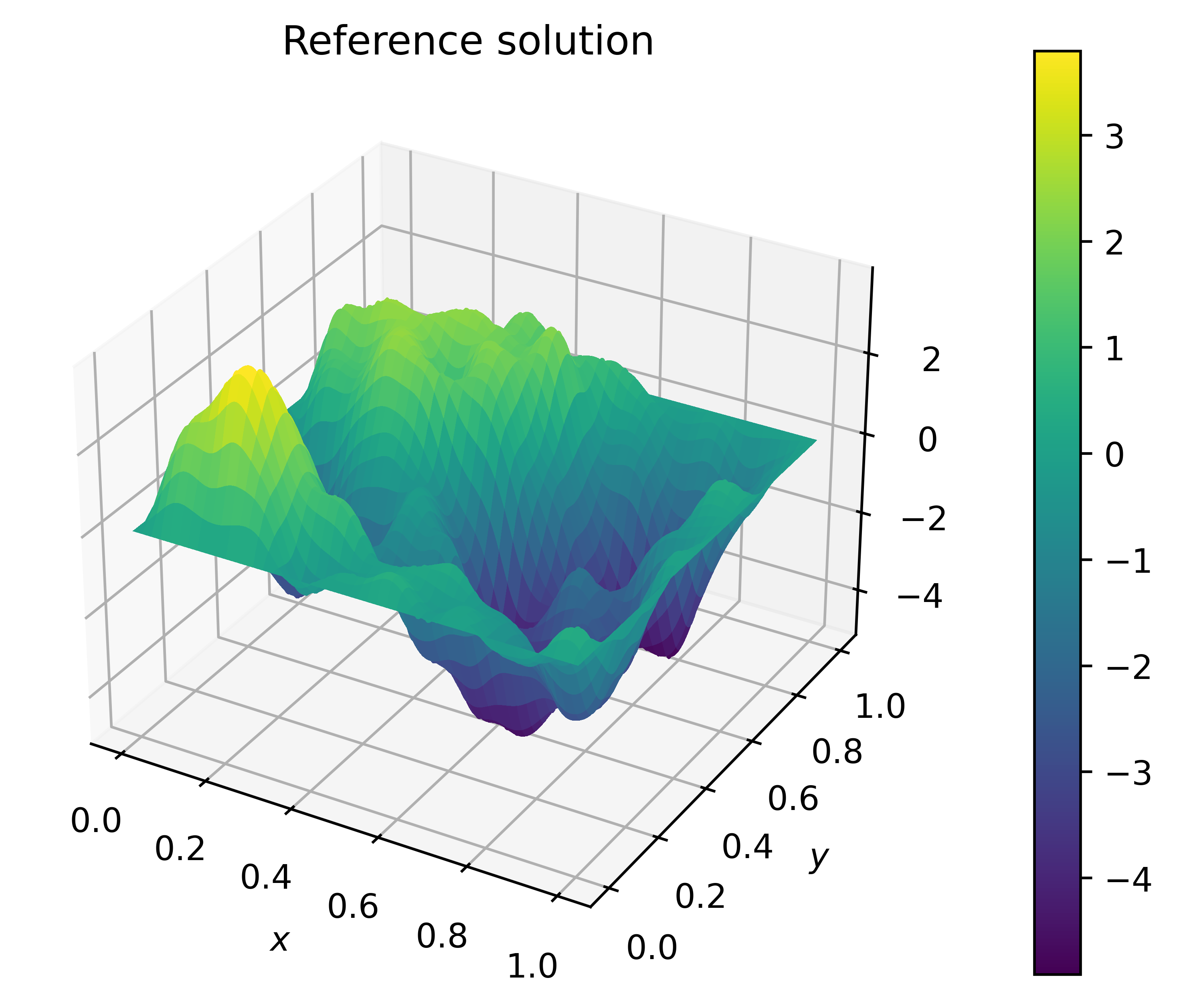}
    \caption*{(2b)}
  \end{subfigure}
  \hfill
  \begin{subfigure}{0.32\textwidth}
    \centering
    \includegraphics[width=.965\linewidth]{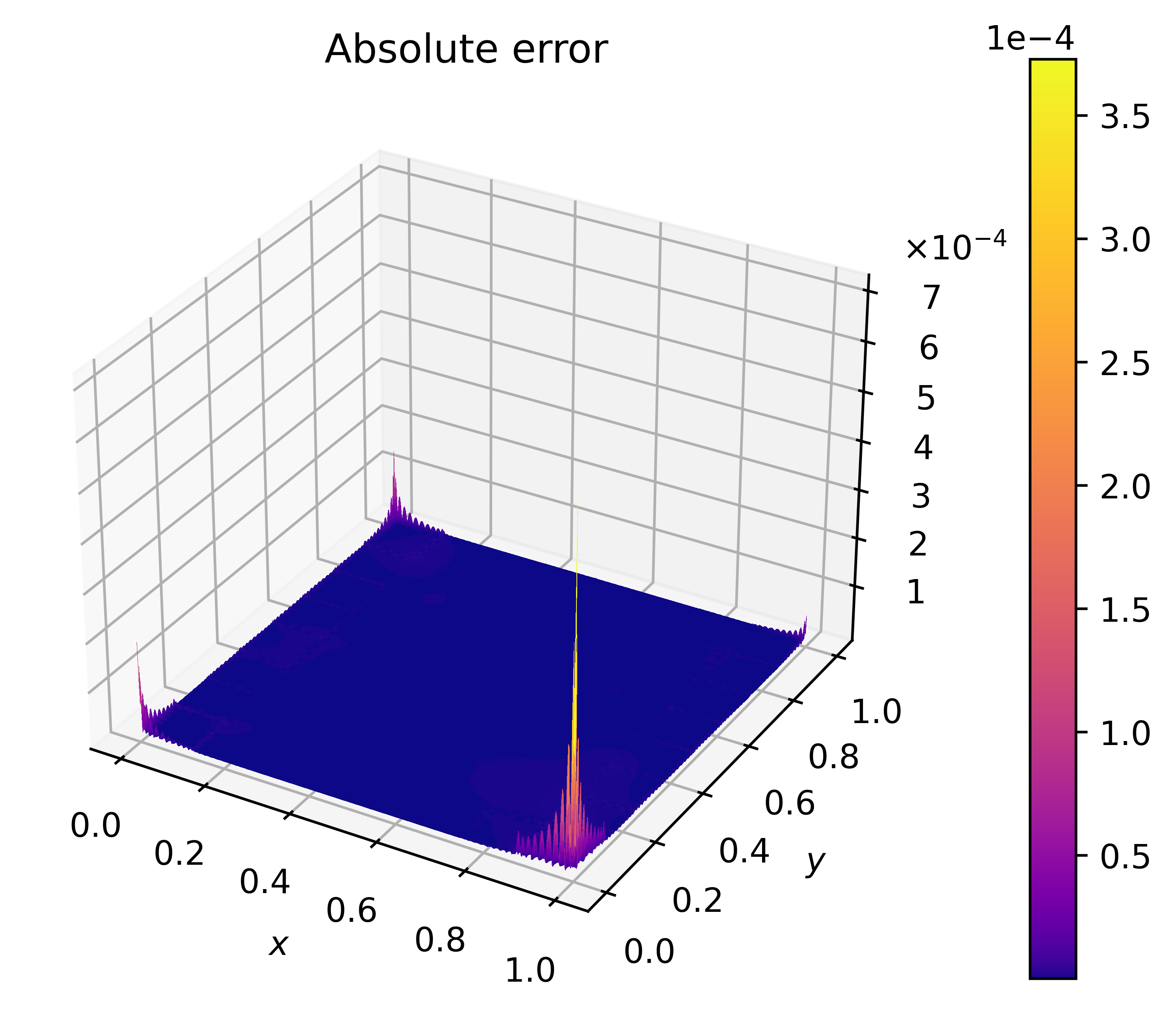}
    \caption*{(2c)}
  \end{subfigure}

  \begin{subfigure}{0.495\textwidth}
    \centering
    \includegraphics[width=\linewidth]{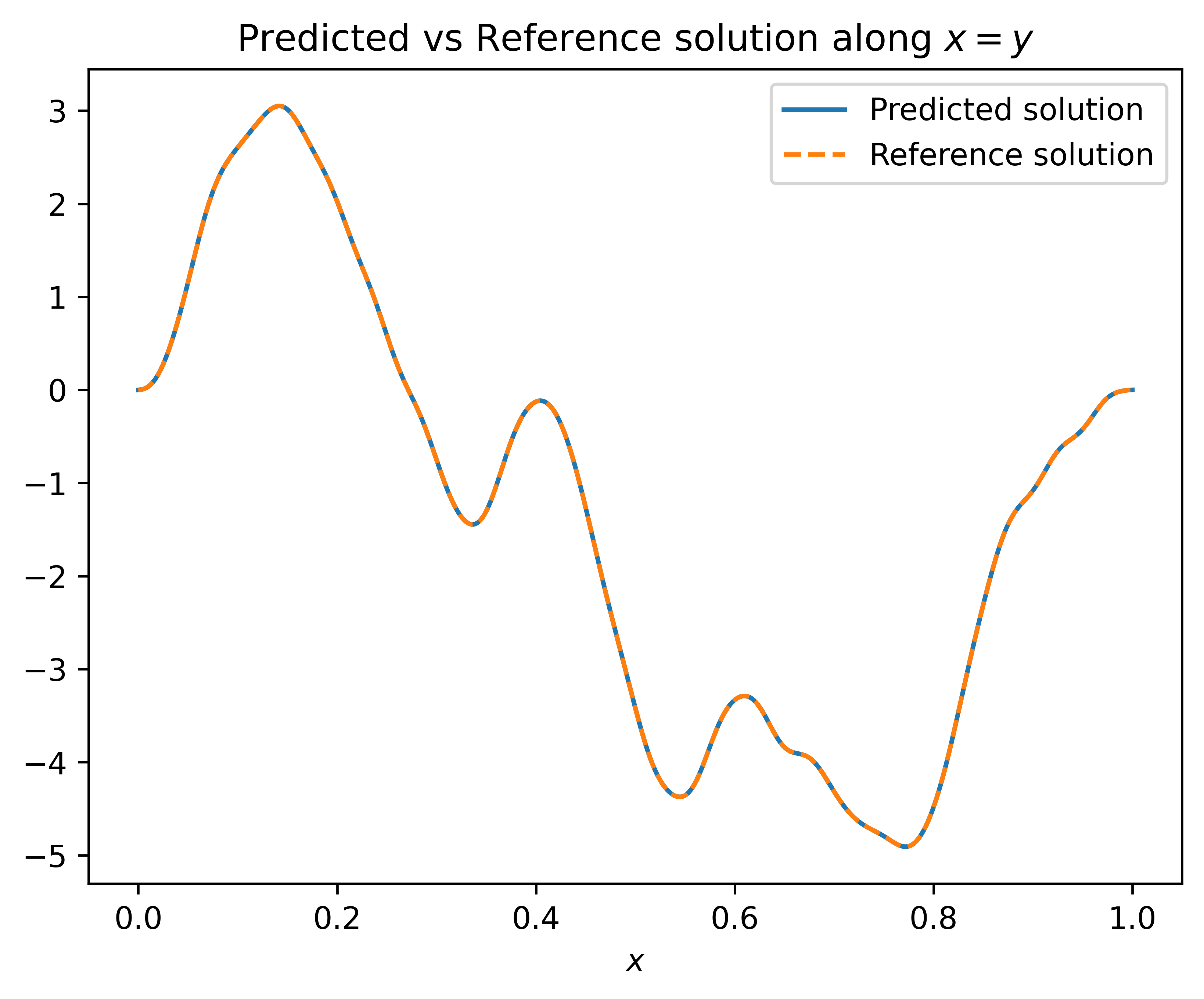}
    \caption*{(3a)}
  \end{subfigure}
  \hfill
  \begin{subfigure}{0.495\textwidth}
    \centering
    \includegraphics[width=.965\linewidth]{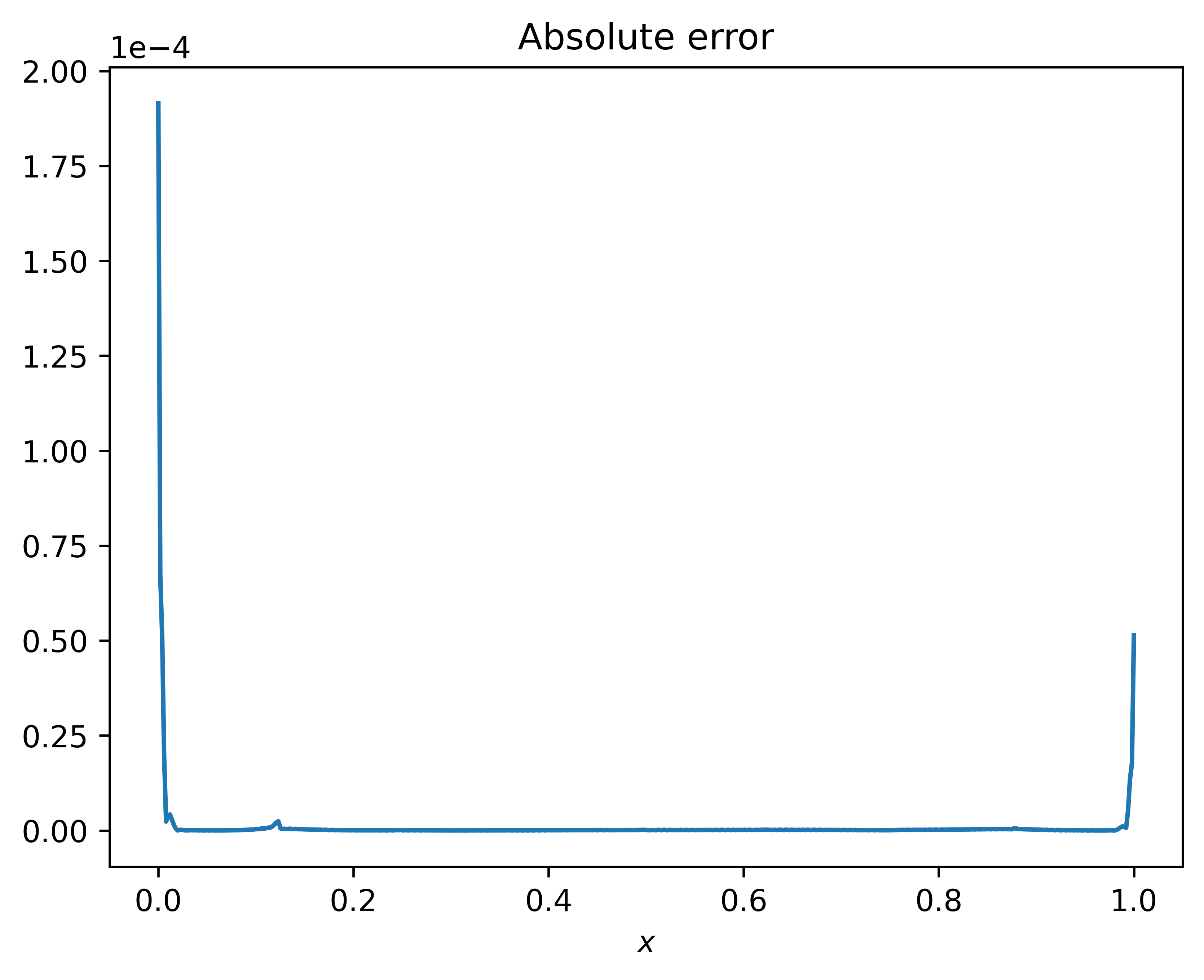}
    \caption*{(3b)}
  \end{subfigure}
  \caption{Solution and error of 2-D Poisson equation for the case $g=0$ and $f$ is an approximate GRF with $\alpha=32$.
  The DDELM-NN solution is computed on $8\times8$ subdomains.
  (1a) DDELM-NN solution 2D plot, (1b) Reference solution 2D plot, (1c) Absolute error of DDELM-NN solution 2D plot.
  (2a) DDELM-NN solution 3D plot, (2b) Reference solution 3D plot, (2c) Absolute error of DDELM-NN solution 3D plot.
  (3a) DDELM-NN and reference solution along the line $x=y$, (3b) Absolute error of DDELM-NN solution along the line $x=y$.}
  \label{fig:poisson_grf320b}
\end{figure}
Here, DDELM-NN shows greater improvement in both number of iterations and computation time.
The $8\times8$ subdomain DDELM-NN solution and error plots are shown in~\cref{fig:poisson_grf320b}.
We can observe sharp peaks in the error at the corners of the entire domain $\Omega$.
Fortunately, the error is still quite small and is not as large in the interior.

\subsubsection{Variable Coefficient Poisson Equation}
We consider a variable coefficient Poisson equation
\begin{equation*}
  \left\{\begin{aligned}
  -\nabla \cdot (\rho \nabla u) &= 1 &&\text{ in } \Omega=(0,1)^2, \\
  u &= 0 &&\text{ on } \partial \Omega,
  \end{aligned}\right.
\end{equation*}
where $\rho(\mathbf{x}) = \tanh(\mathrm{grf}(\mathbf{x})) + 1.1$.
We use an approximate GRF as in \eqref{eqn:grf} with $\alpha=32$ for the grf term in $\rho$.
Again, the exact solution is not known and we use the FEM solution as reference.
For this example, we also provide results with quadruple the number of neurons and training points: $2^{18}/N$ neurons and a $640/\sqrt{N}\times 640/\sqrt{N}$ grid for each local neural network.
We summarize the results in \cref{tab:varpoisson}.
\begin{table}
  \caption{
    Relative $L^2$ and $H^1$ errors, number of CG iterations, and wall-clock time solving the variable coefficient Poisson equation.
    $N$ is the number of subdomains.
    }
  \label{tab:varpoisson}
\centering
\begin{tabular}{ccccccc}
  \toprule
  Number of & \multirow{2}{*}{Method} & \multirow{2}{*}{$N$} & \multirow{2}{*}{$L^2$ error} & \multirow{2}{*}{$H^1$ error} & \multirow{2}{*}{Iterations} & Wall-clock\\
  neurons &&&&&& time (sec)\\
  \midrule \midrule
  \multirow{12}{*}{65,536}
  & \multirow{4}{*}{DDELM}
  & $2\times2$ & 9.10e-03 & 6.38e-02 & 1,580 & 460.87 \\
  & & $4\times4$ & 6.93e-03 & 6.33e-02 & 4,638 & 14.39 \\
  & & $8\times8$ & 1.27e-03 & 2.54e-02 & 10,402 & 13.13 \\
  & & $16\times16$ & 2.54e-04 & 8.27e-03 & 13,584 & 22.85 \\
  \cmidrule{2-7}
  & \multirow{4}{*}{DDELM-CS}
  & $2\times2$ & 1.06e-01 & 1.30e-01 & 967 & 461.18 \\
  & & $4\times4$ & 9.82e-02 & 1.33e-01 & 2,881 & 17.36 \\
  & & $8\times8$ & 1.39e-02 & 4.04e-02 & 5,637 & 12.08 \\
  & & $16\times16$ & 4.75e-04 & 7.23e-03 & 6,771 & 12.65 \\
  \cmidrule{2-7}
  & \multirow{4}{*}{DDELM-NN}
  & $2\times2$ & 1.11e-01 & 1.33e-01 & 123 & 922.90 \\
  & & $4\times4$ & 8.88e-02 & 1.31e-01 & 227 & 16.04 \\
  & & $8\times8$ & 1.30e-02 & 3.97e-02 & 435 & 1.87 \\
  & & $16\times16$ & 2.63e-04 & 6.84e-03 & 943 & 2.37 \\
  \midrule
  \multirow{9}{*}{262,144}
  & \multirow{3}{*}{DDELM}
  & $4\times4$ & 4.11e-07 & 2.95e-05 & 9,022 & 478.60 \\
  & & $8\times8$ & 2.24e-07 & 1.39e-05 & 18,718 & 54.78 \\
  & & $16\times16$ & 1.64e-07 & 9.06e-06 & 24,581 & 27.47 \\
  \cmidrule{2-7}
  & \multirow{3}{*}{DDELM-CS}
  & $4\times4$ & 5.50e-07 & 2.91e-05 & 5,299 & 477.14 \\
  & & $8\times8$ & 2.14e-07 & 1.28e-05 & 8,939 & 46.81 \\
  & & $16\times16$ & 2.74e-07 & 7.64e-06 & 10,494 & 21.16 \\
  \cmidrule{2-7}
  & \multirow{3}{*}{DDELM-NN}
  & $4\times4$ & 5.54e-07 & 2.96e-05 & 358 & 934.76 \\
  & & $8\times8$ & 3.10e-07 & 1.27e-05 & 542 & 55.68 \\
  & & $16\times16$ & 1.68e-07 & 8.29e-06 & 992 & 5.49 \\
  \bottomrule
\end{tabular}
\end{table}
\begin{figure}
  \begin{subfigure}{0.32\textwidth}
    \centering
    \includegraphics[width=\linewidth]{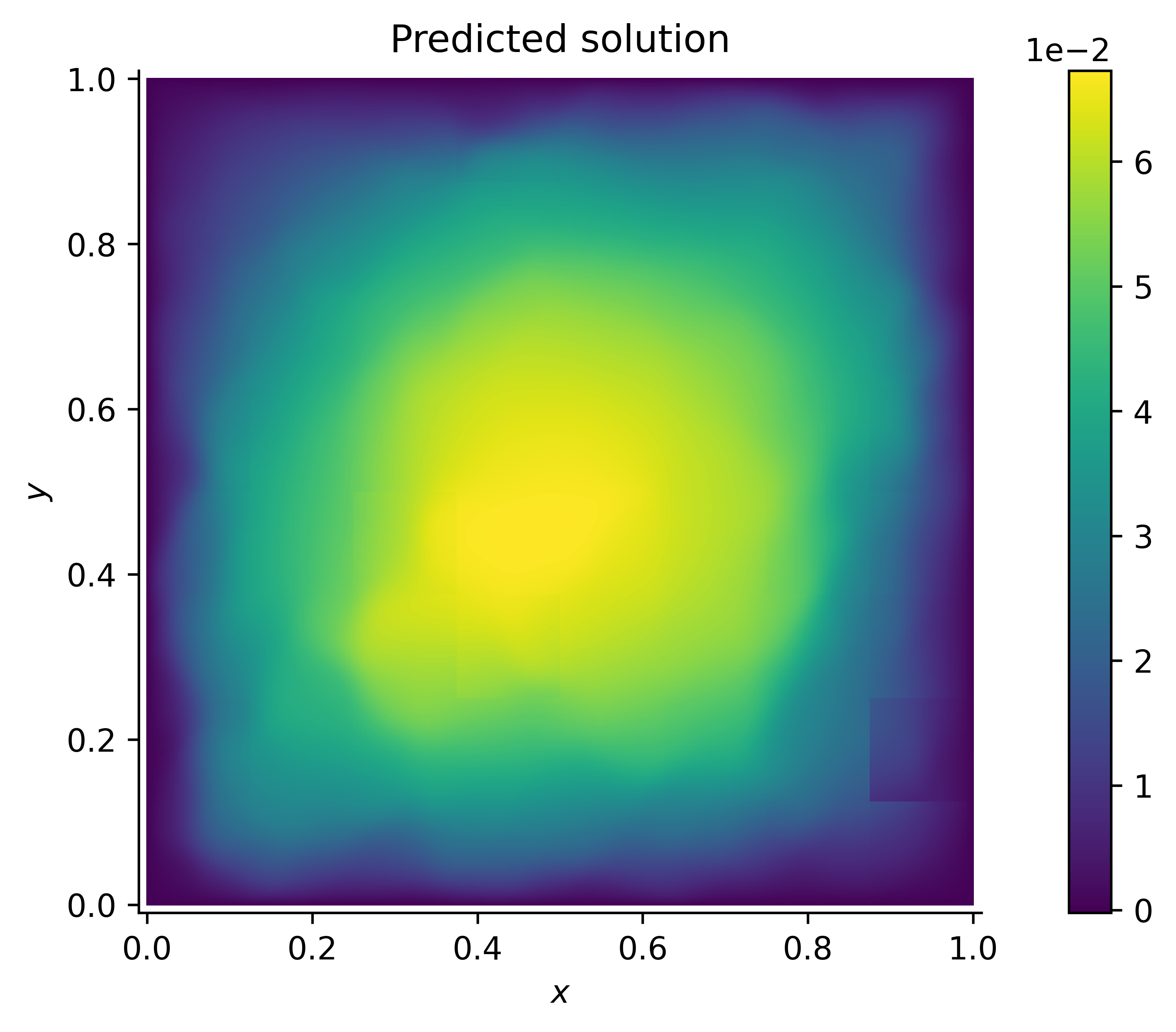}
    \caption*{(1a)}
  \end{subfigure}
  \hfill
  \begin{subfigure}{0.32\textwidth}
    \centering
    \includegraphics[width=.965\linewidth]{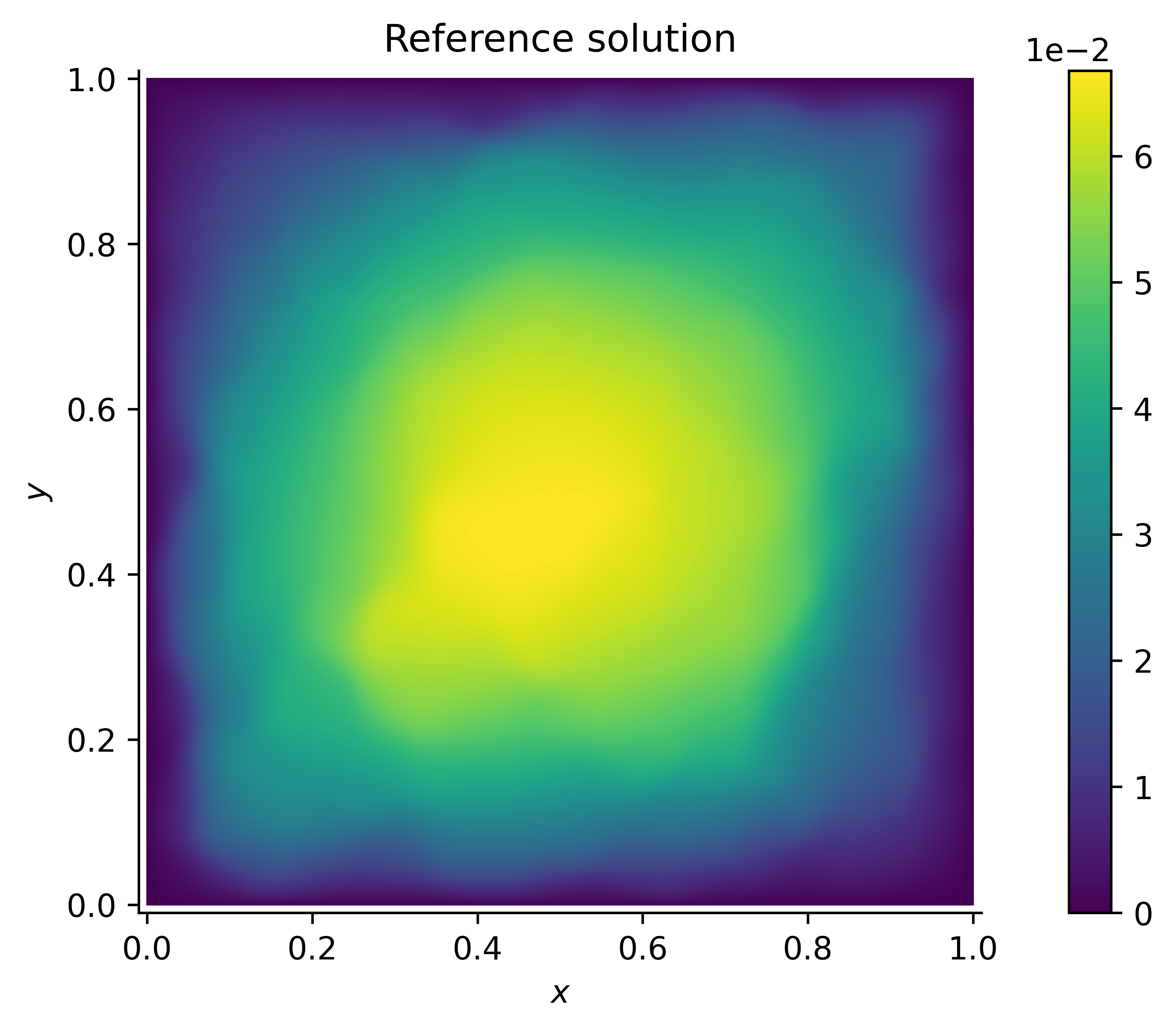}
    \caption*{(1b)}
  \end{subfigure}
  \begin{subfigure}{0.32\textwidth}
    \centering
    \includegraphics[width=.965\linewidth]{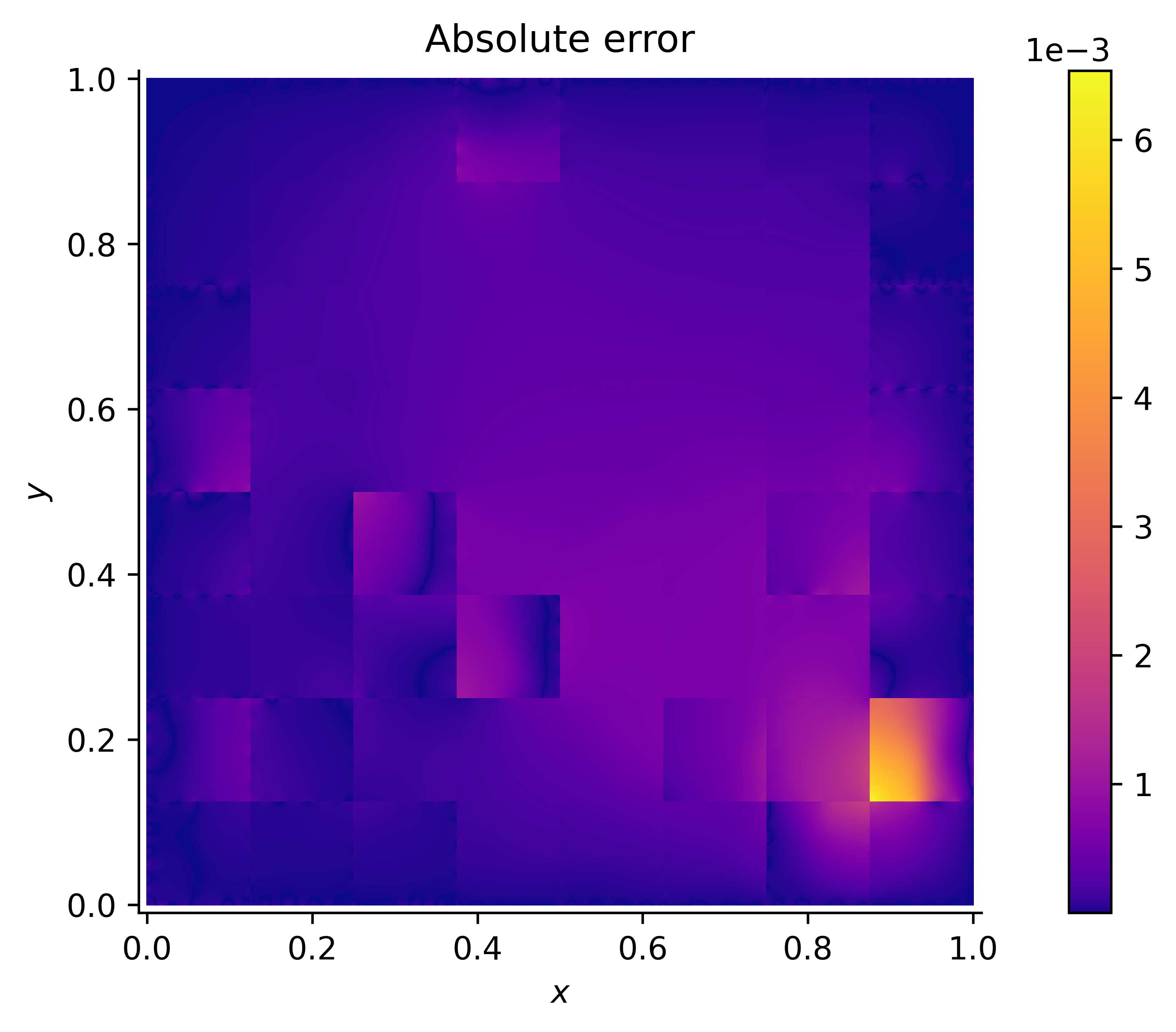}
    \caption*{(1c)}
  \end{subfigure}

  \begin{subfigure}{0.32\textwidth}
    \centering
    \includegraphics[width=\linewidth]{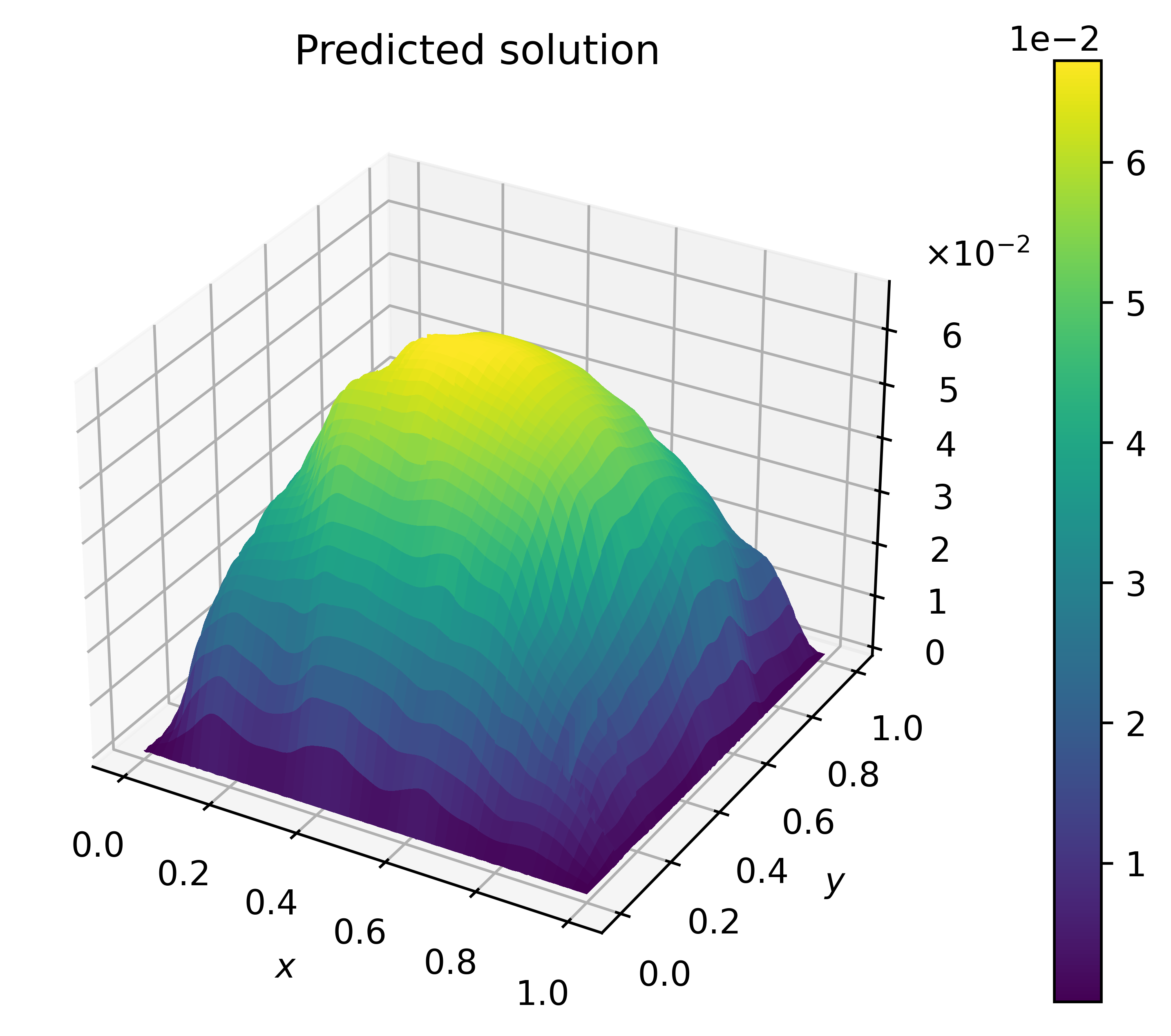}
    \caption*{(2a)}
  \end{subfigure}
  \hfill
  \begin{subfigure}{0.32\textwidth}
    \centering
    \includegraphics[width=.965\linewidth]{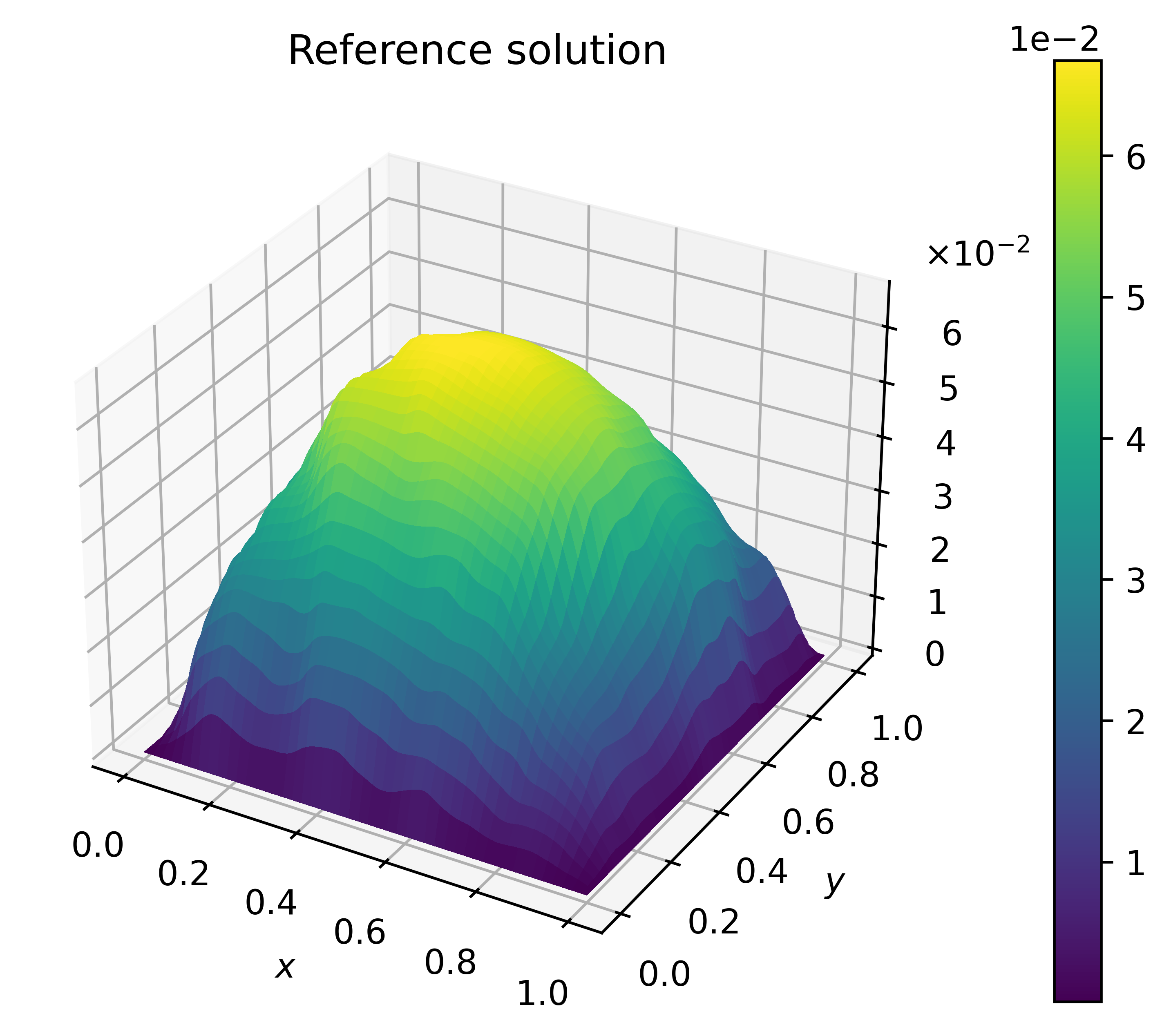}
    \caption*{(2b)}
  \end{subfigure}
  \hfill
  \begin{subfigure}{0.32\textwidth}
    \centering
    \includegraphics[width=.965\linewidth]{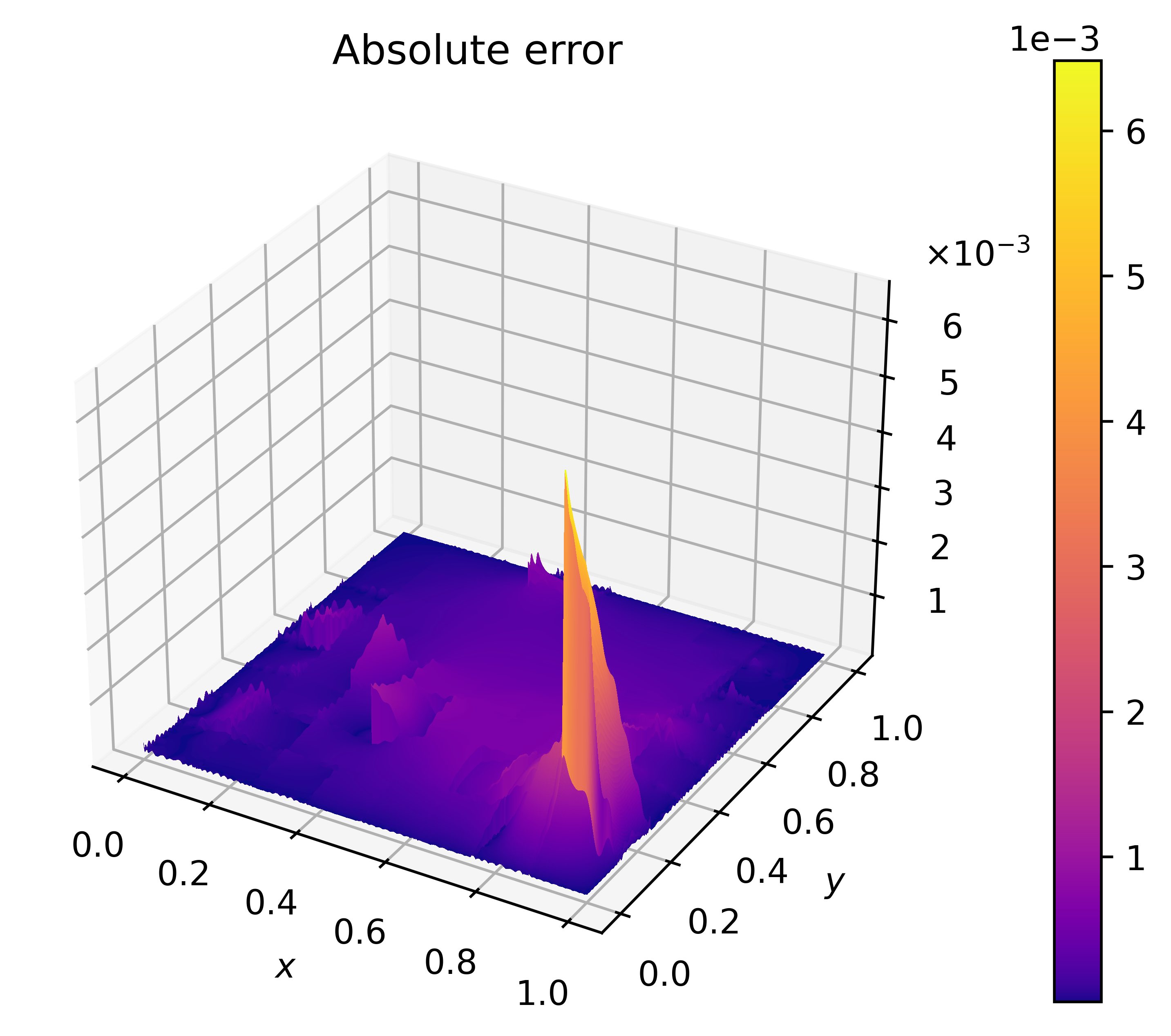}
    \caption*{(2c)}
  \end{subfigure}

  \begin{subfigure}{0.495\textwidth}
    \centering
    \includegraphics[width=\linewidth]{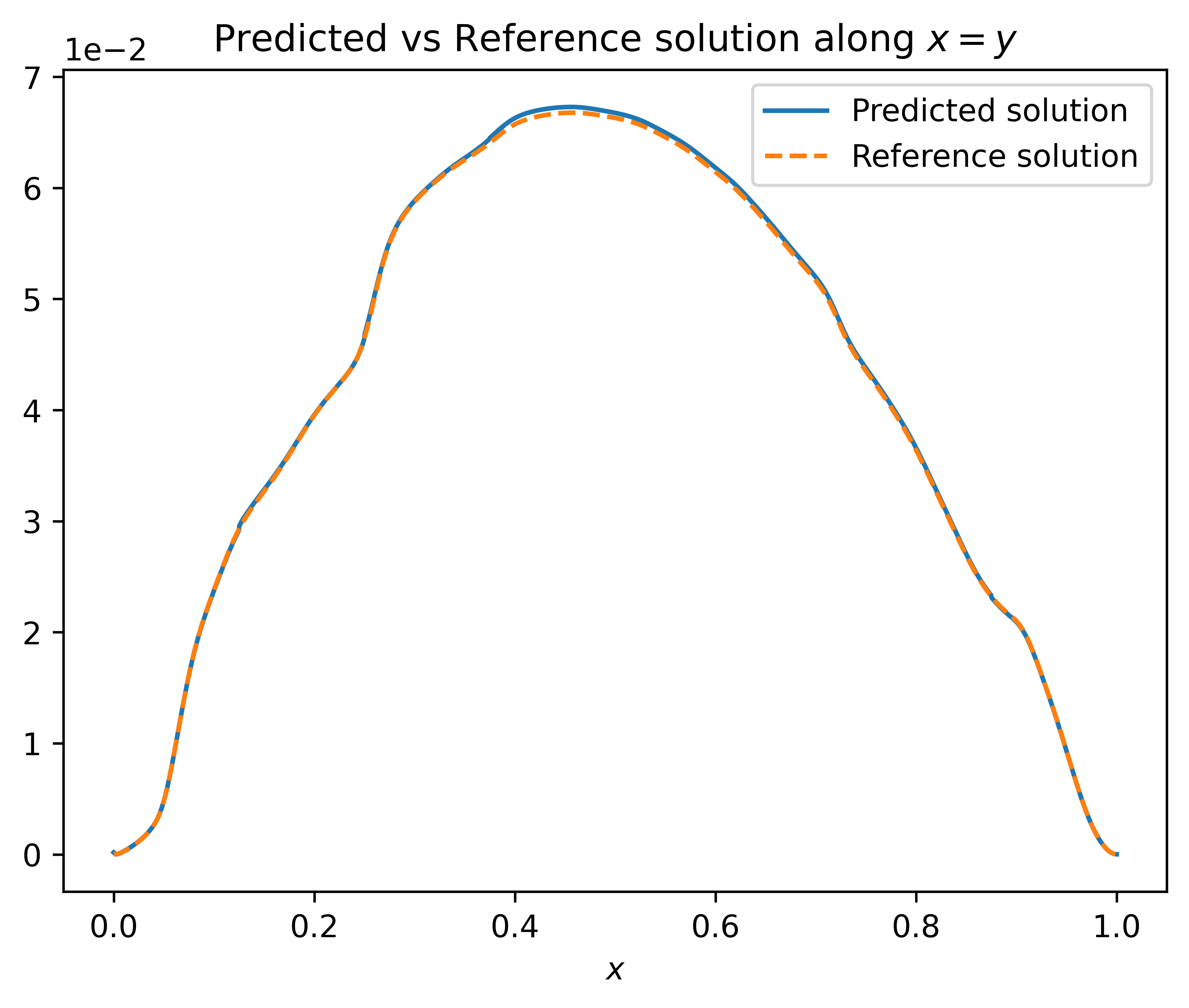}
    \caption*{(3a)}
  \end{subfigure}
  \hfill
  \begin{subfigure}{0.495\textwidth}
    \centering
    \includegraphics[width=.965\linewidth]{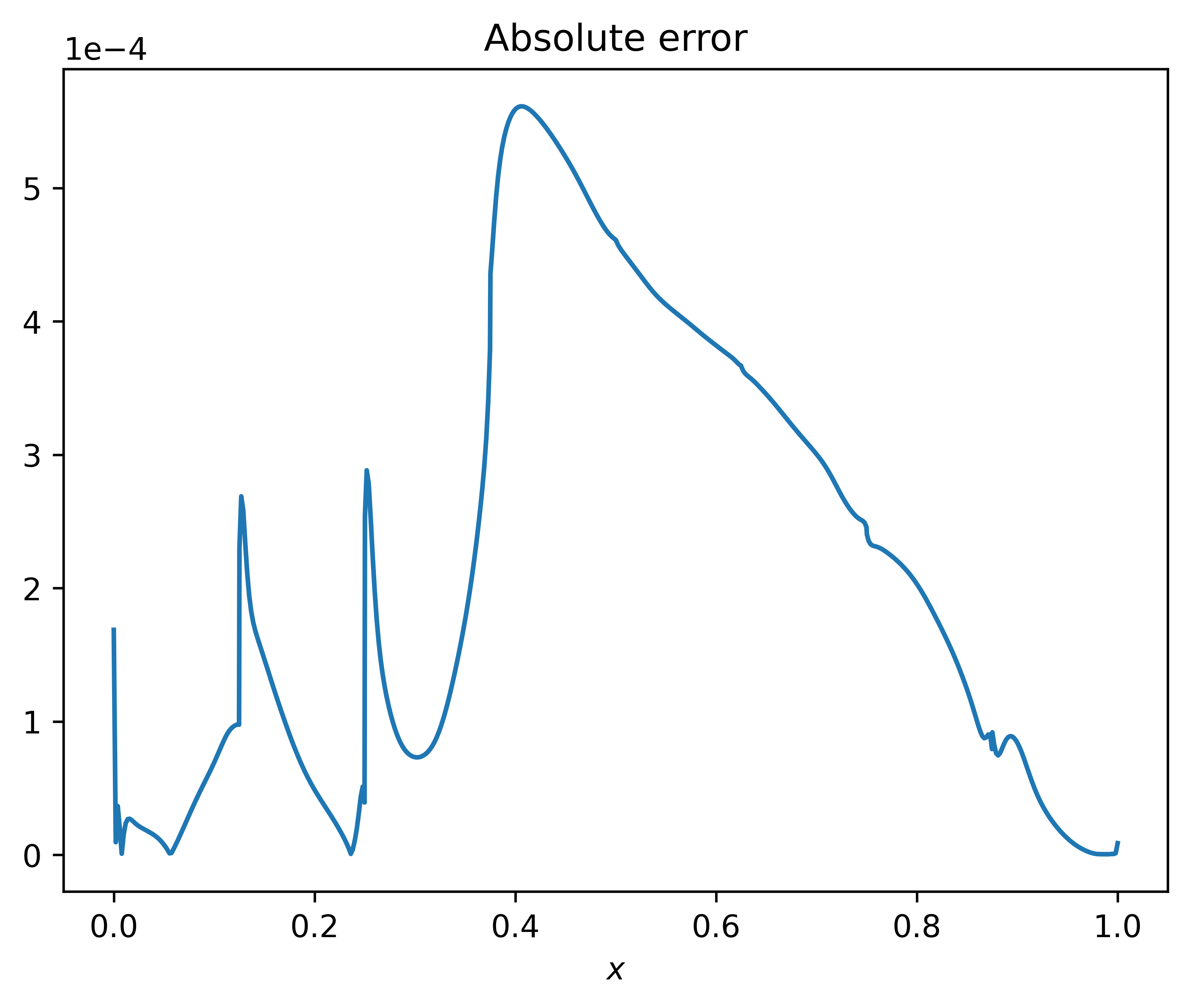}
    \caption*{(3b)}
  \end{subfigure}
  \caption{Solution and error of 2-D Poisson equation with variable coefficient when $\alpha=32$.
  The DDELM-NN solution is computed on $8\times8$ subdomains using 65,536 neurons.
  (1a) DDELM-NN solution 2D plot, (1b) Reference solution 2D plot, (1c) Absolute error of DDELM-NN solution 2D plot.
  (2a) DDELM-NN solution 3D plot, (2b) Reference solution 3D plot, (2c) Absolute error of DDELM-NN solution 3D plot.
  (3a) DDELM-NN and reference solution along the line $x=y$, (3b) Absolute error of DDELM-NN solution along the line $x=y$.}
  \label{fig:poivar_grf320b1f}
\end{figure}
\begin{figure}
  \begin{subfigure}{0.32\textwidth}
    \centering
    \includegraphics[width=\linewidth]{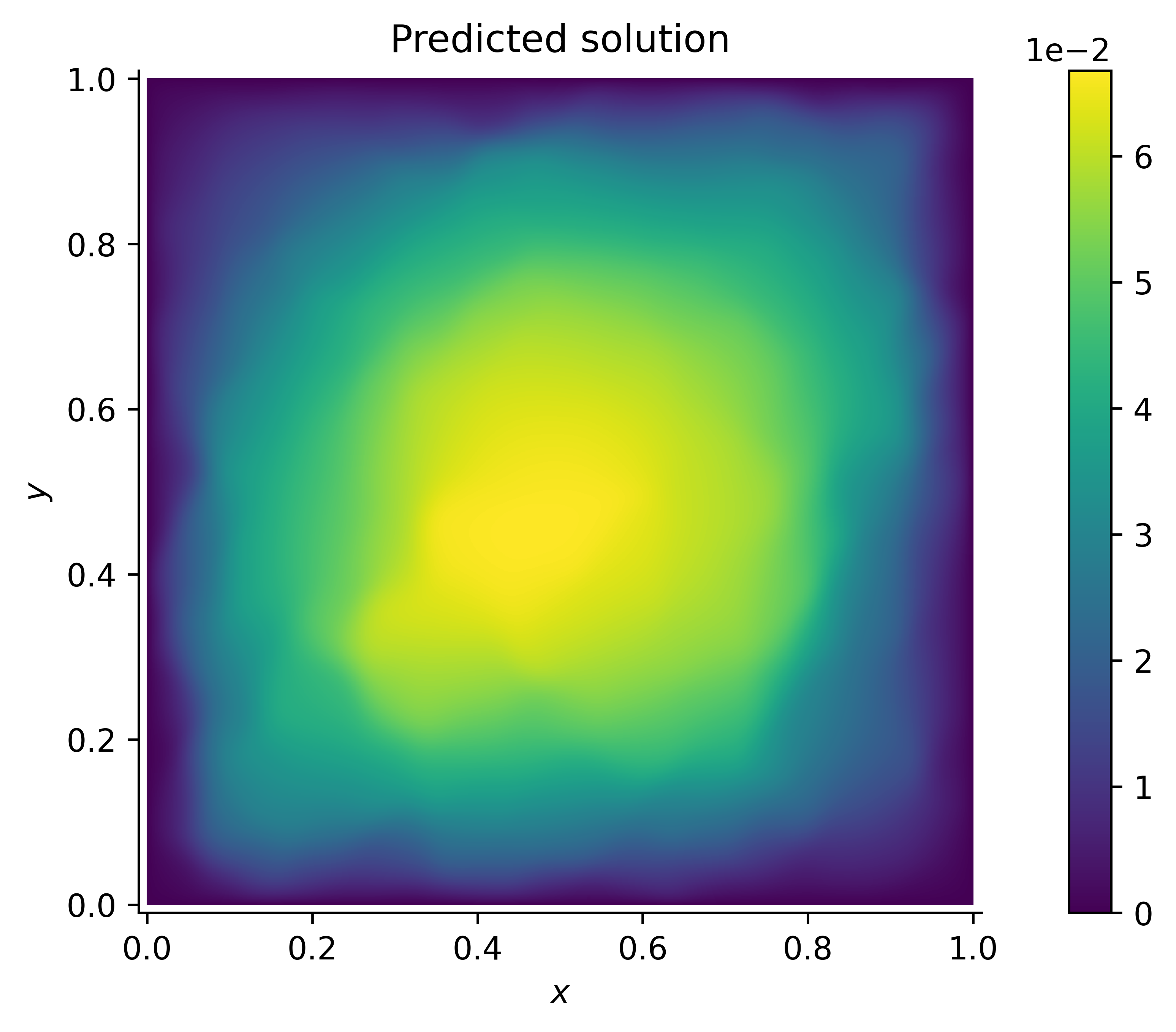}
    \caption*{(1a)}
  \end{subfigure}
  \hfill
  \begin{subfigure}{0.32\textwidth}
    \centering
    \includegraphics[width=.965\linewidth]{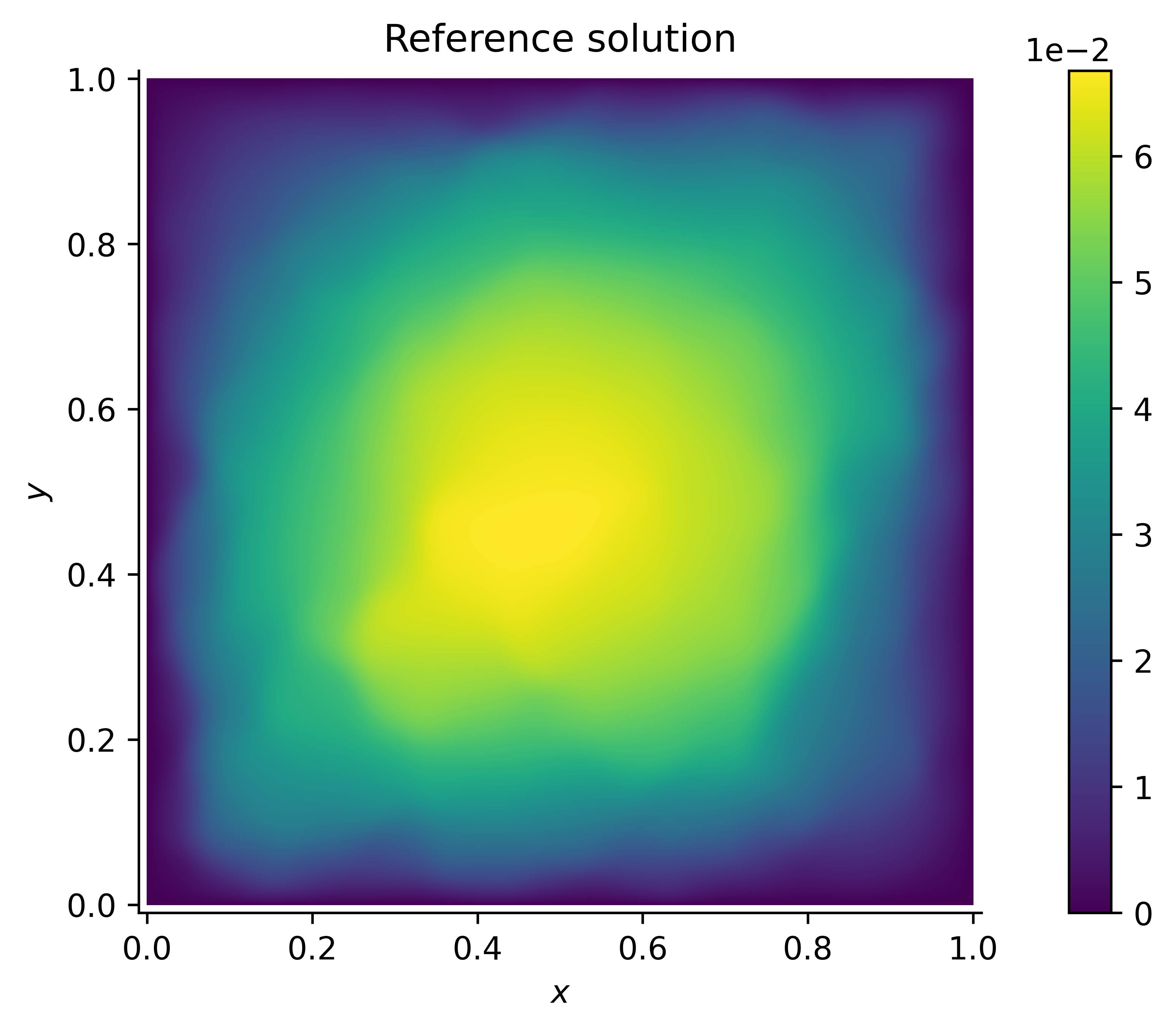}
    \caption*{(1b)}
  \end{subfigure}
  \begin{subfigure}{0.32\textwidth}
    \centering
    \includegraphics[width=.965\linewidth]{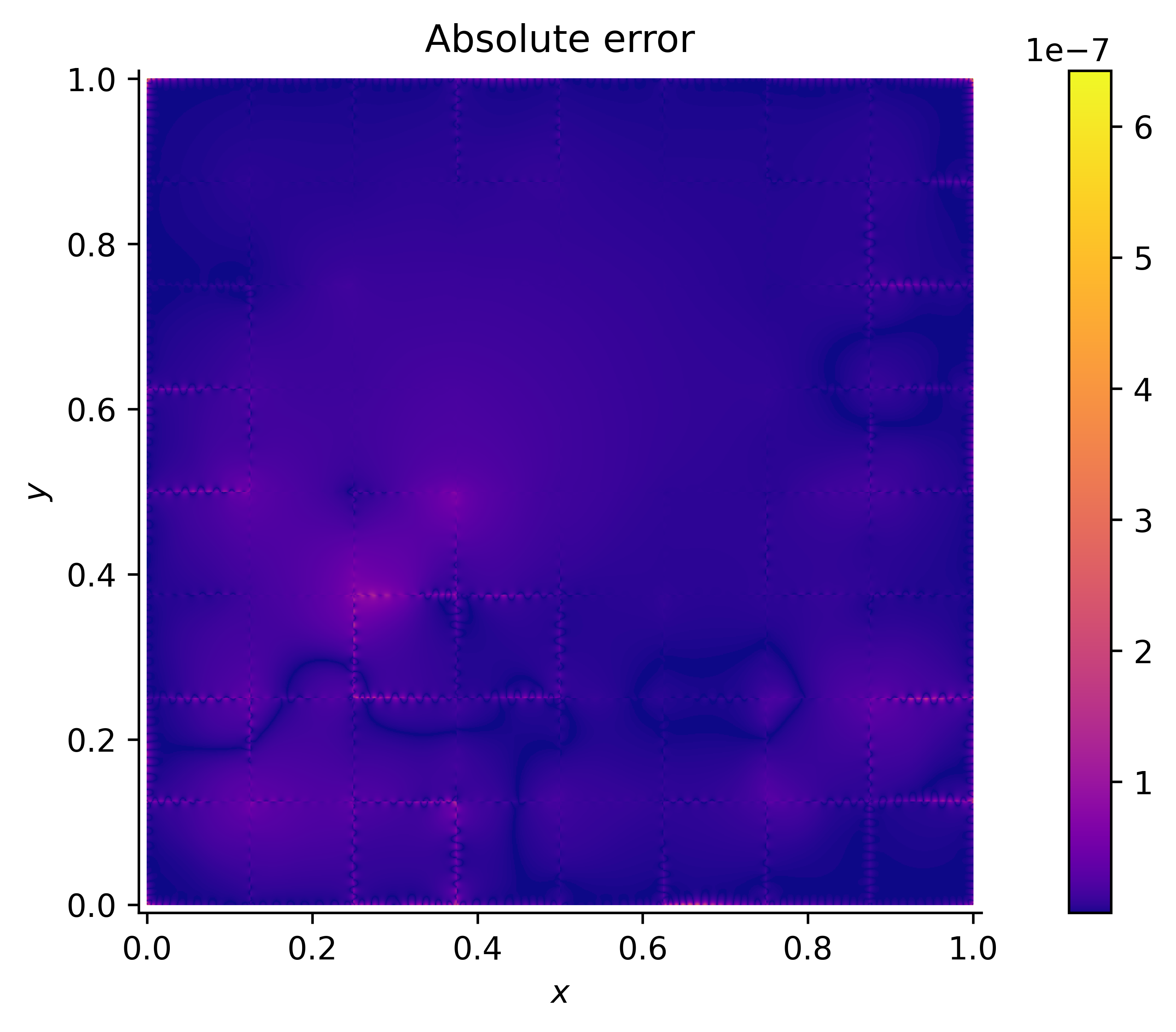}
    \caption*{(1c)}
  \end{subfigure}

  \begin{subfigure}{0.32\textwidth}
    \centering
    \includegraphics[width=\linewidth]{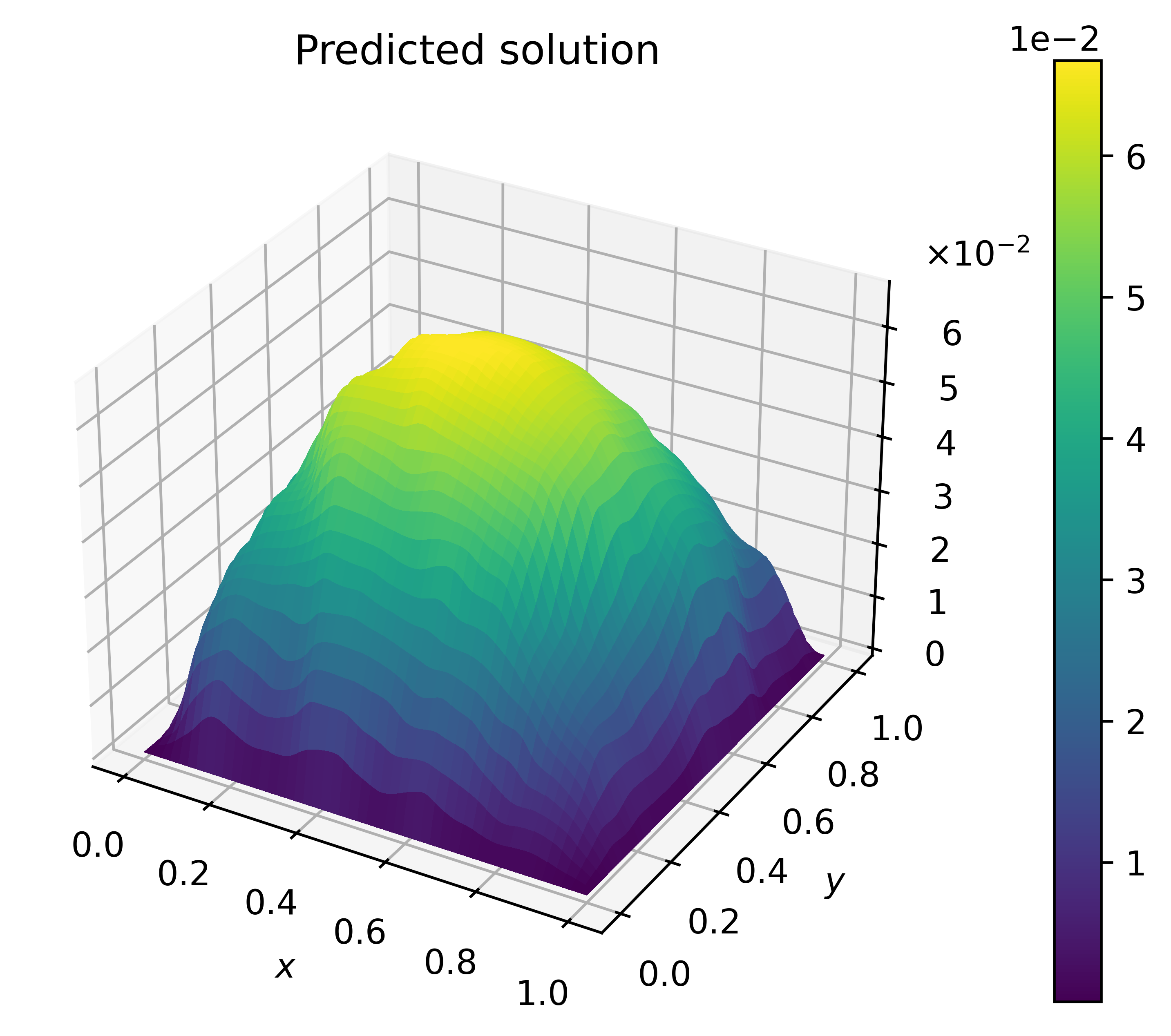}
    \caption*{(2a)}
  \end{subfigure}
  \hfill
  \begin{subfigure}{0.32\textwidth}
    \centering
    \includegraphics[width=.965\linewidth]{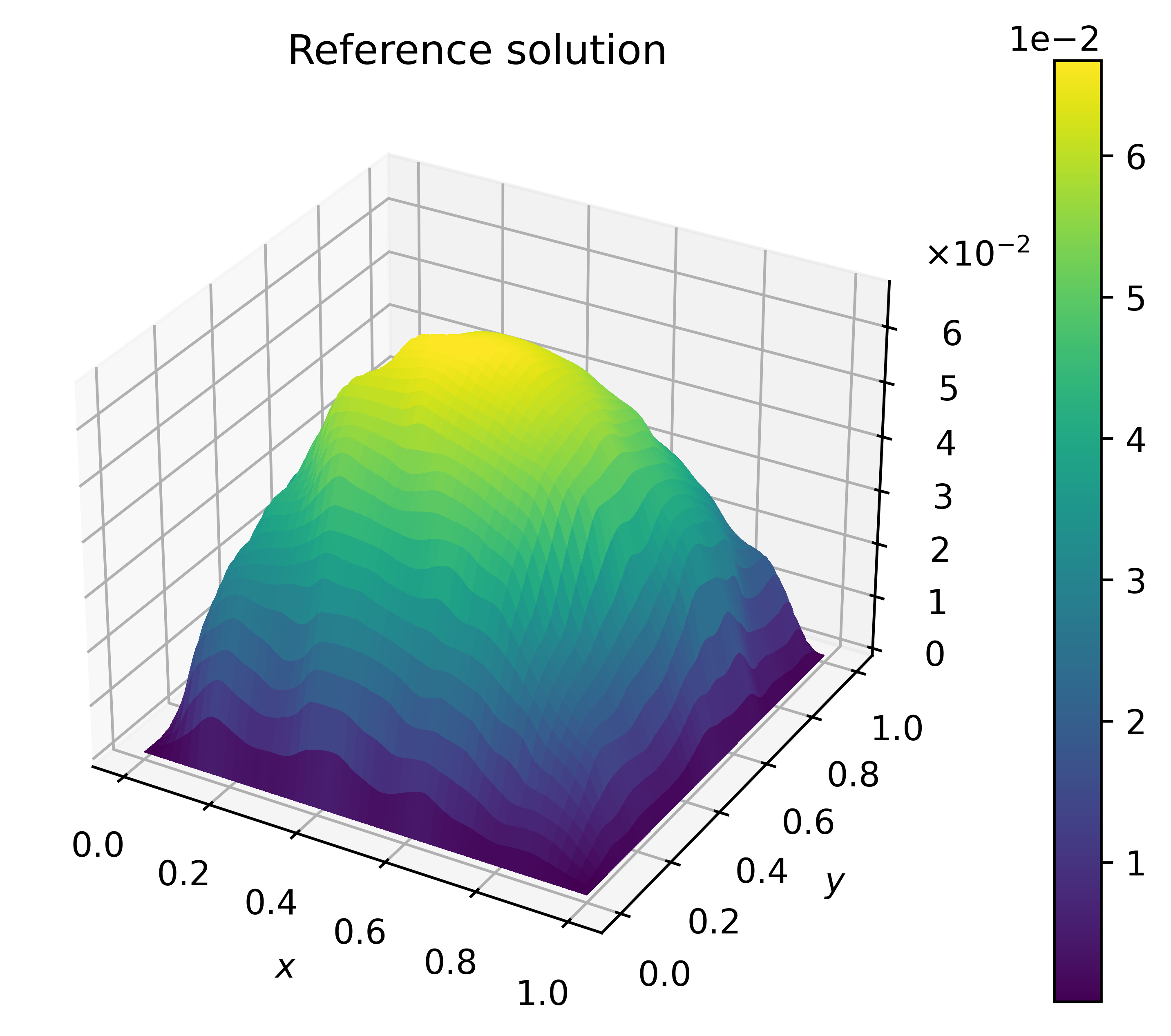}
    \caption*{(2b)}
  \end{subfigure}
  \hfill
  \begin{subfigure}{0.32\textwidth}
    \centering
    \includegraphics[width=.965\linewidth]{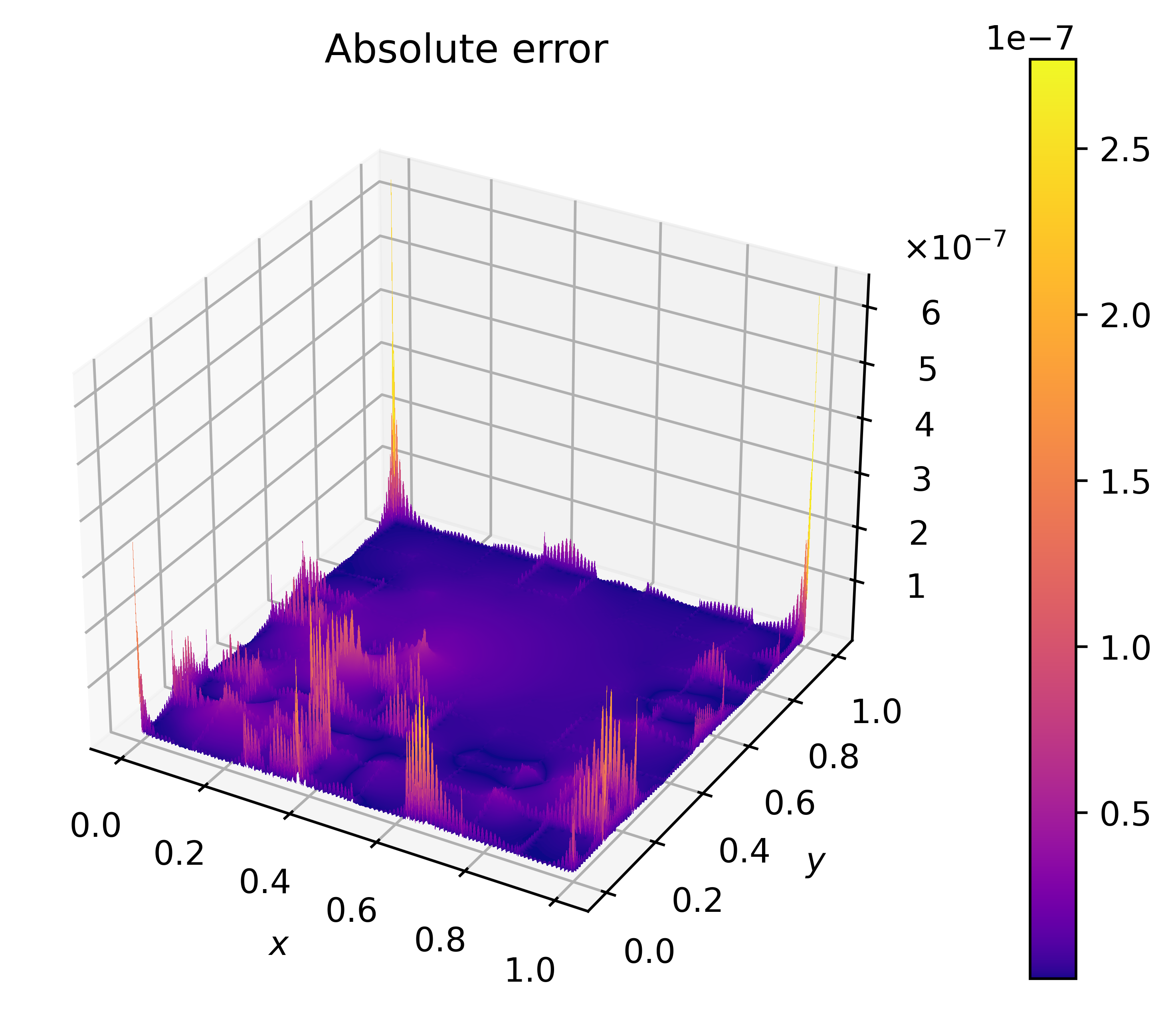}
    \caption*{(2c)}
  \end{subfigure}

  \begin{subfigure}{0.495\textwidth}
    \centering
    \includegraphics[width=\linewidth]{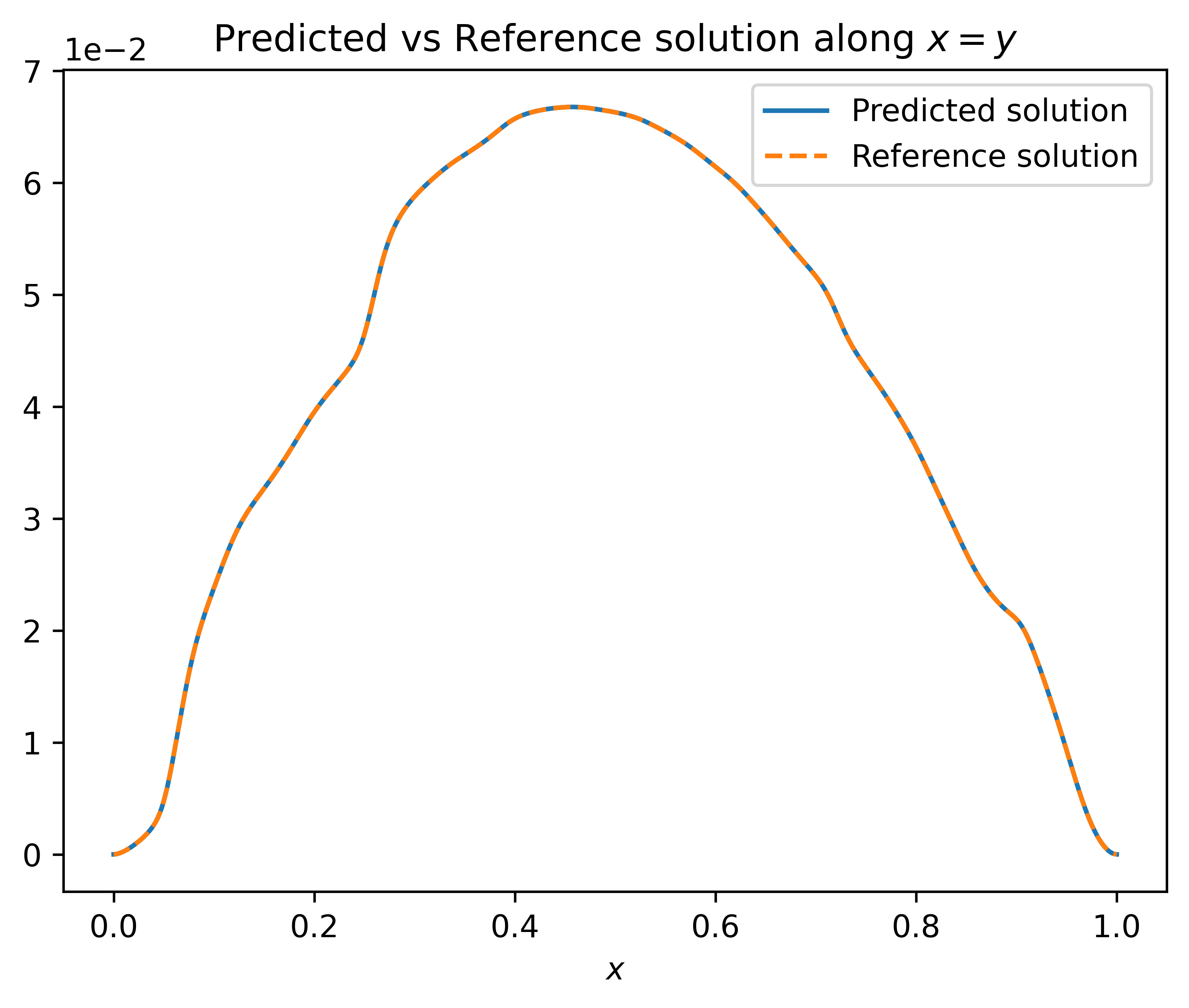}
    \caption*{(3a)}
  \end{subfigure}
  \hfill
  \begin{subfigure}{0.495\textwidth}
    \centering
    \includegraphics[width=.965\linewidth]{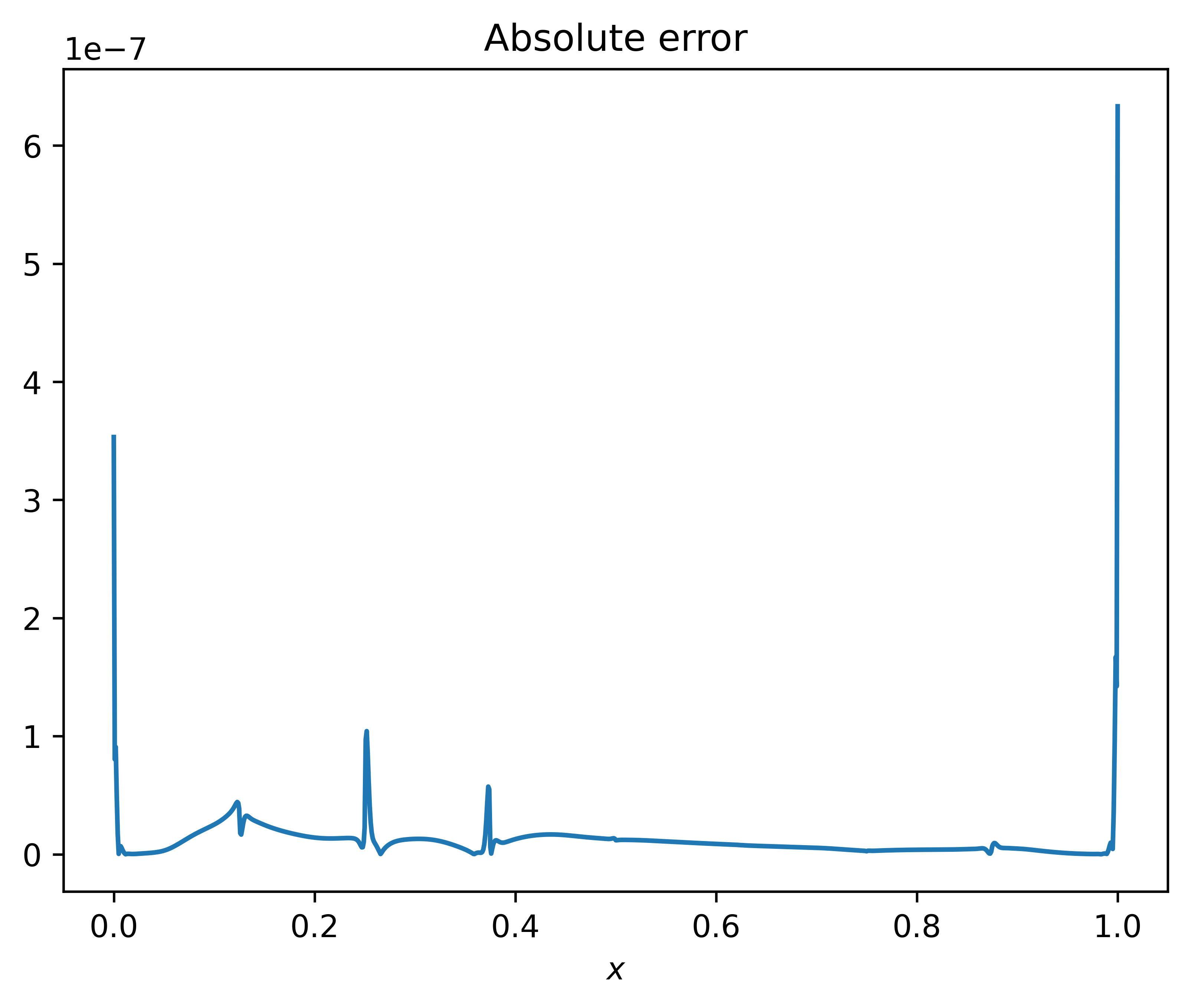}
    \caption*{(3b)}
  \end{subfigure}
  \caption{Solution and error of 2-D Poisson equation with variable coefficient when $\alpha=32$.
  The DDELM-NN solution is computed on $8\times8$ subdomains using 262,144 neurons.
  (1a) DDELM-NN solution 2D plot, (1b) Reference solution 2D plot, (1c) Absolute error of DDELM-NN solution 2D plot.
  (2a) DDELM-NN solution 3D plot, (2b) Reference solution 3D plot, (2c) Absolute error of DDELM-NN solution 3D plot.
  (3a) DDELM-NN and reference solution along the line $x=y$, (3b) Absolute error of DDELM-NN solution along the line $x=y$.}
  \label{fig:poivar_grf320b1f_quad}
\end{figure}
The performance of the methods is similar to the Poisson equation with oscillatory right hand side.
DDELM-NN reduces the number of iterations by a factor of nearly twenty and computation time by a factor of six.
Also, note that quadrupling the number of neurons and training points does not increase the number of iterations for DDELM-NN as much as it does for the other methods as $N$ increases.
This property is consistent with the optimality of the classical Neumann-Neumann method~\cite{toselli2005domain}.
Figures~\ref{fig:poivar_grf320b1f} and~\ref{fig:poivar_grf320b1f_quad} show the $8\times8$ subdomain DDELM-NN solution and error plots for the 65,536 and 262,144 neuron cases, respectively.
Although the 65,536 error plot shows particularly large error at one subdomain, the increased expressivity of the 262,144 neuron case shows a much smaller and more uniform error plot.

\subsubsection{Biharmonic Equation}
\begin{figure}
  \begin{subfigure}{0.32\textwidth}
    \centering
    \includegraphics[width=\linewidth]{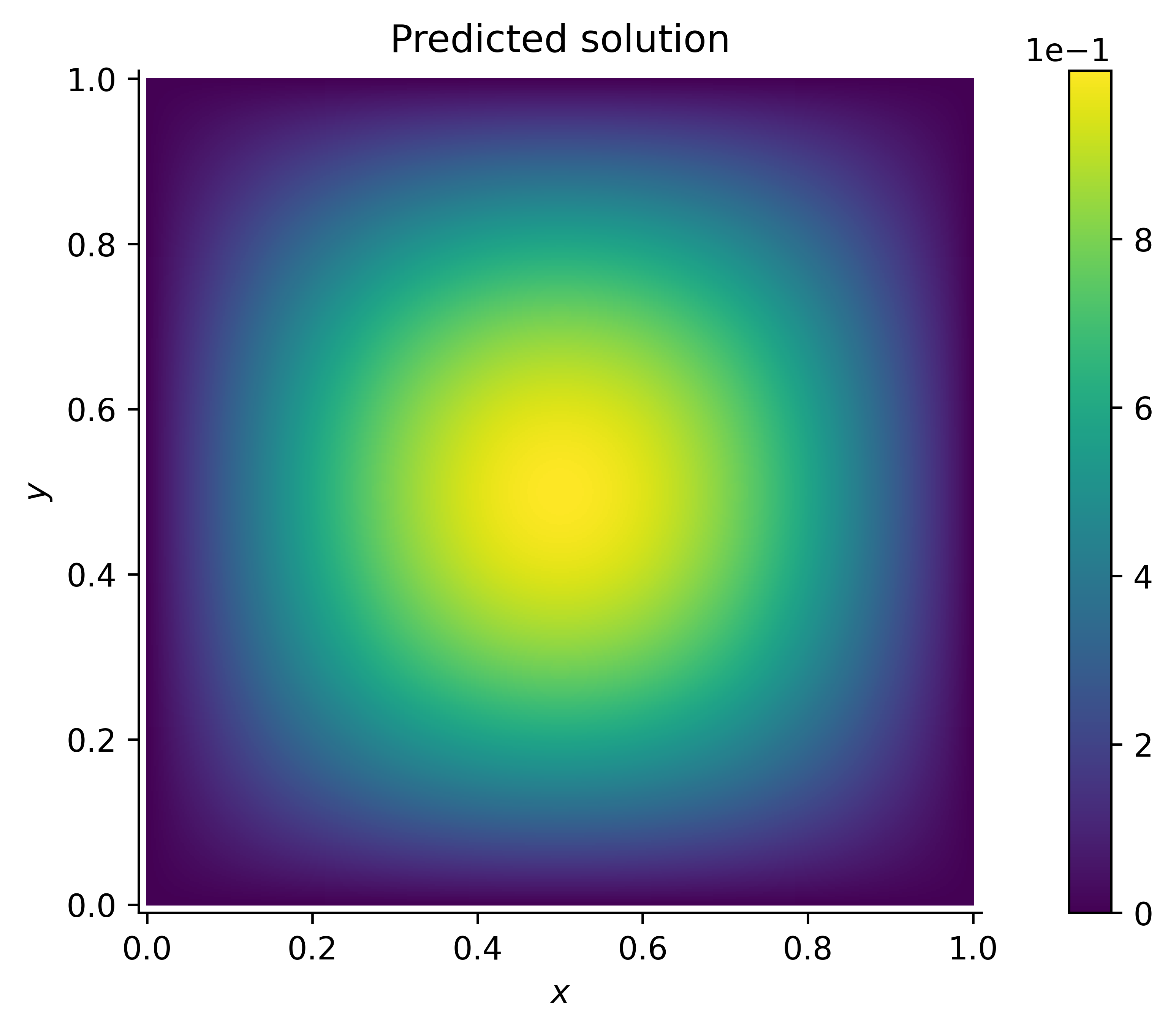}
    \caption*{(1a)}
  \end{subfigure}
  \hfill
  \begin{subfigure}{0.32\textwidth}
    \centering
    \includegraphics[width=.965\linewidth]{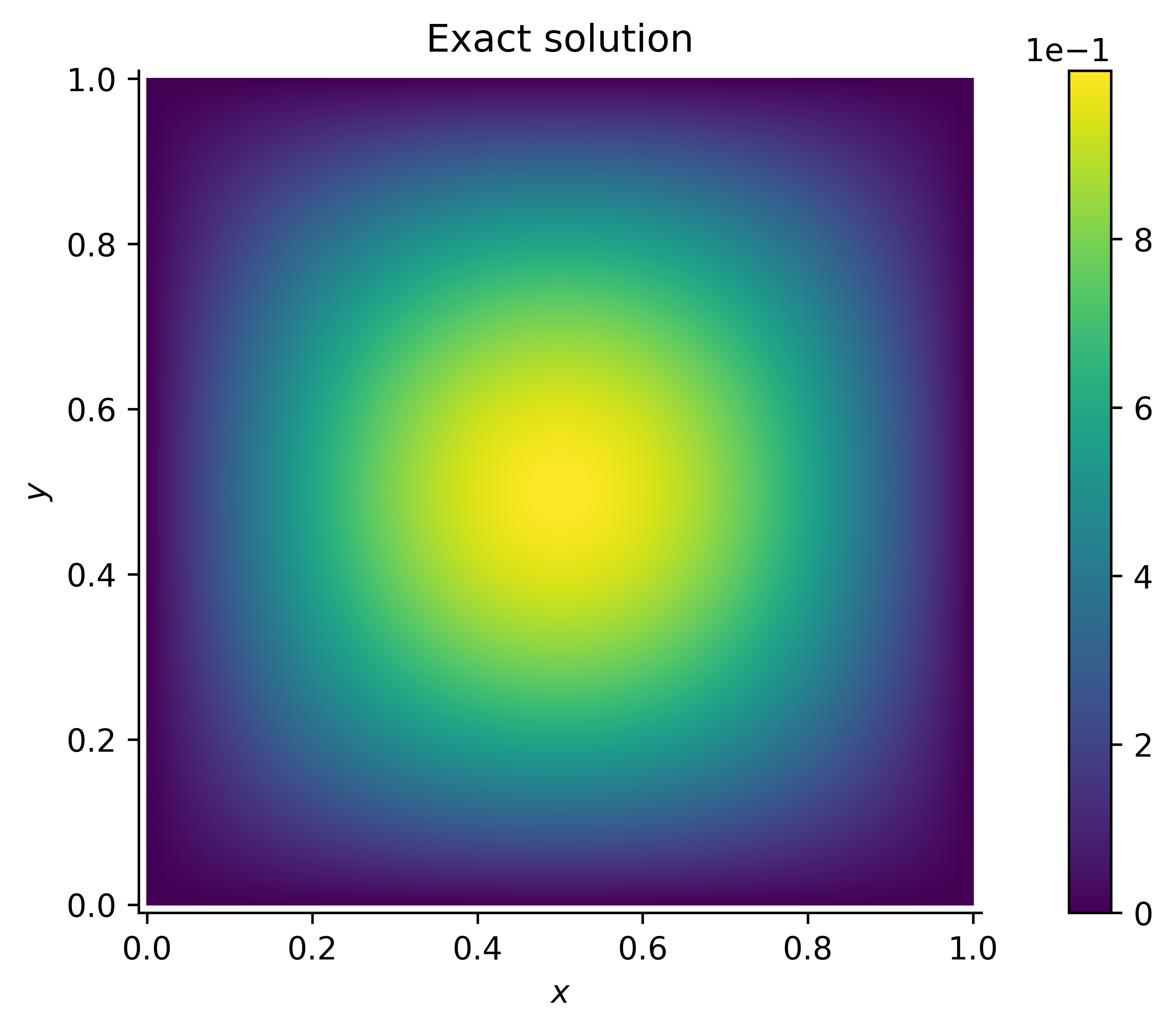}
    \caption*{(1b)}
  \end{subfigure}
  \begin{subfigure}{0.32\textwidth}
    \centering
    \includegraphics[width=.965\linewidth]{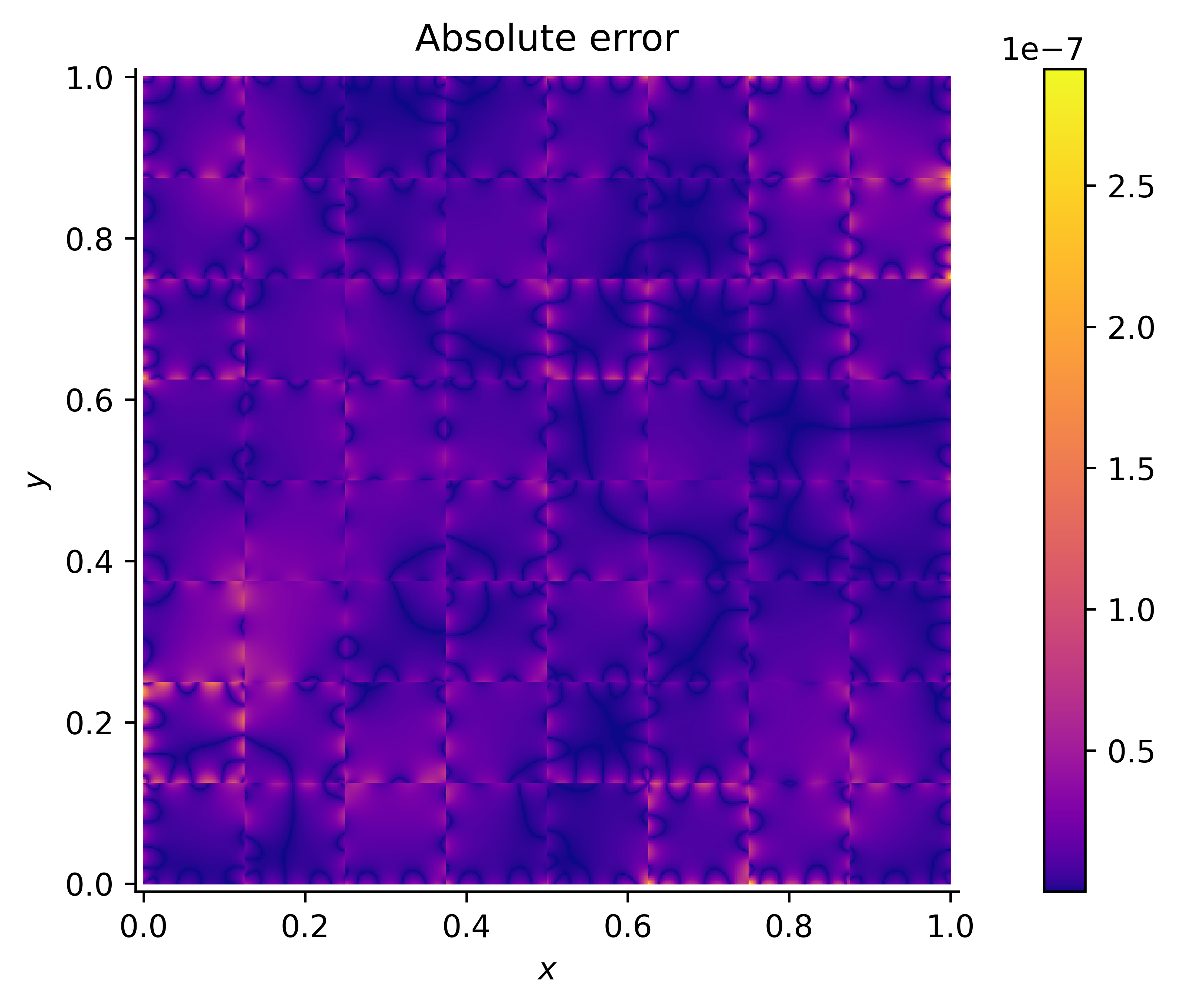}
    \caption*{(1c)}
  \end{subfigure}

  \begin{subfigure}{0.32\textwidth}
    \centering
    \includegraphics[width=\linewidth]{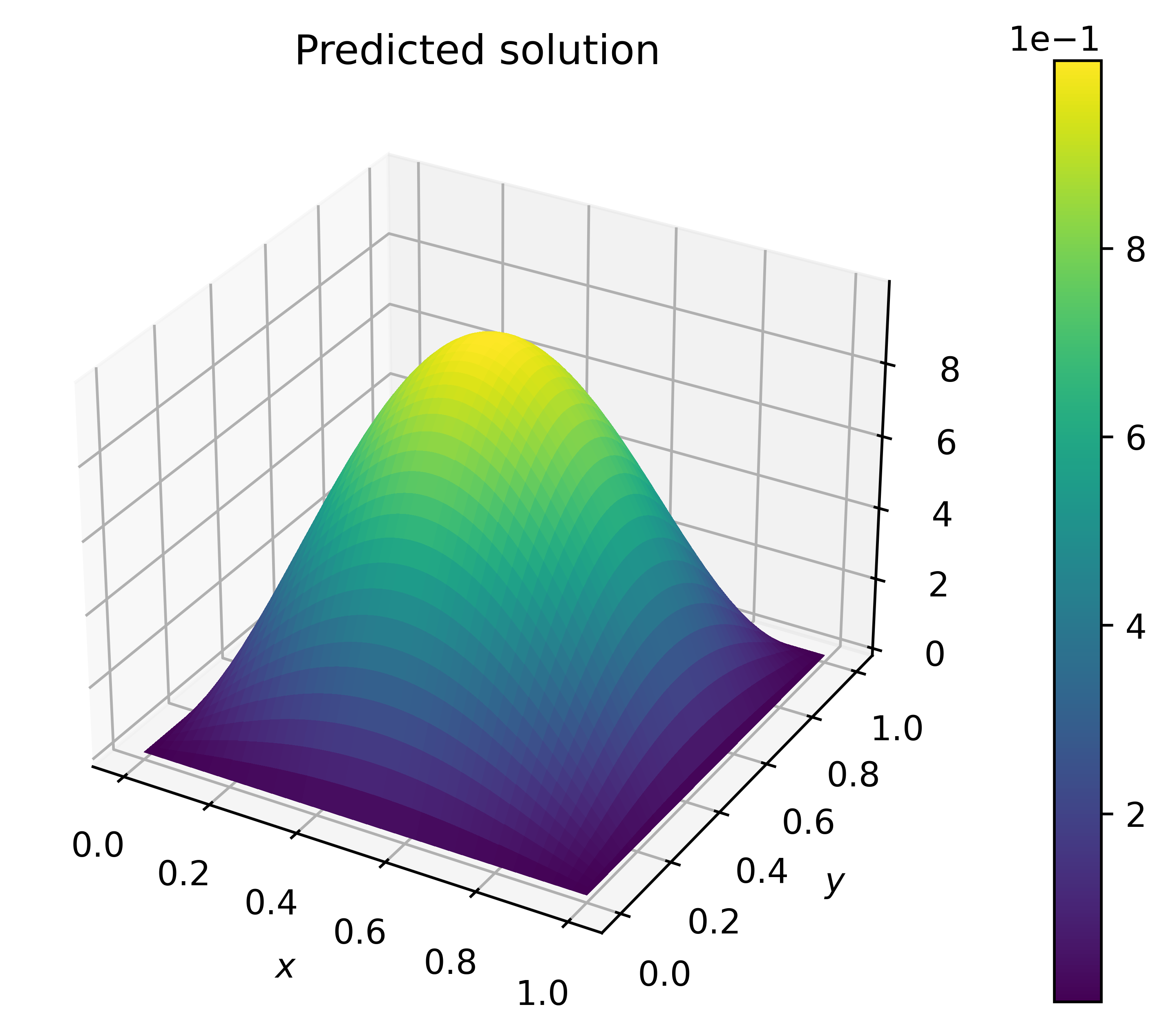}
    \caption*{(2a)}
  \end{subfigure}
  \hfill
  \begin{subfigure}{0.32\textwidth}
    \centering
    \includegraphics[width=.965\linewidth]{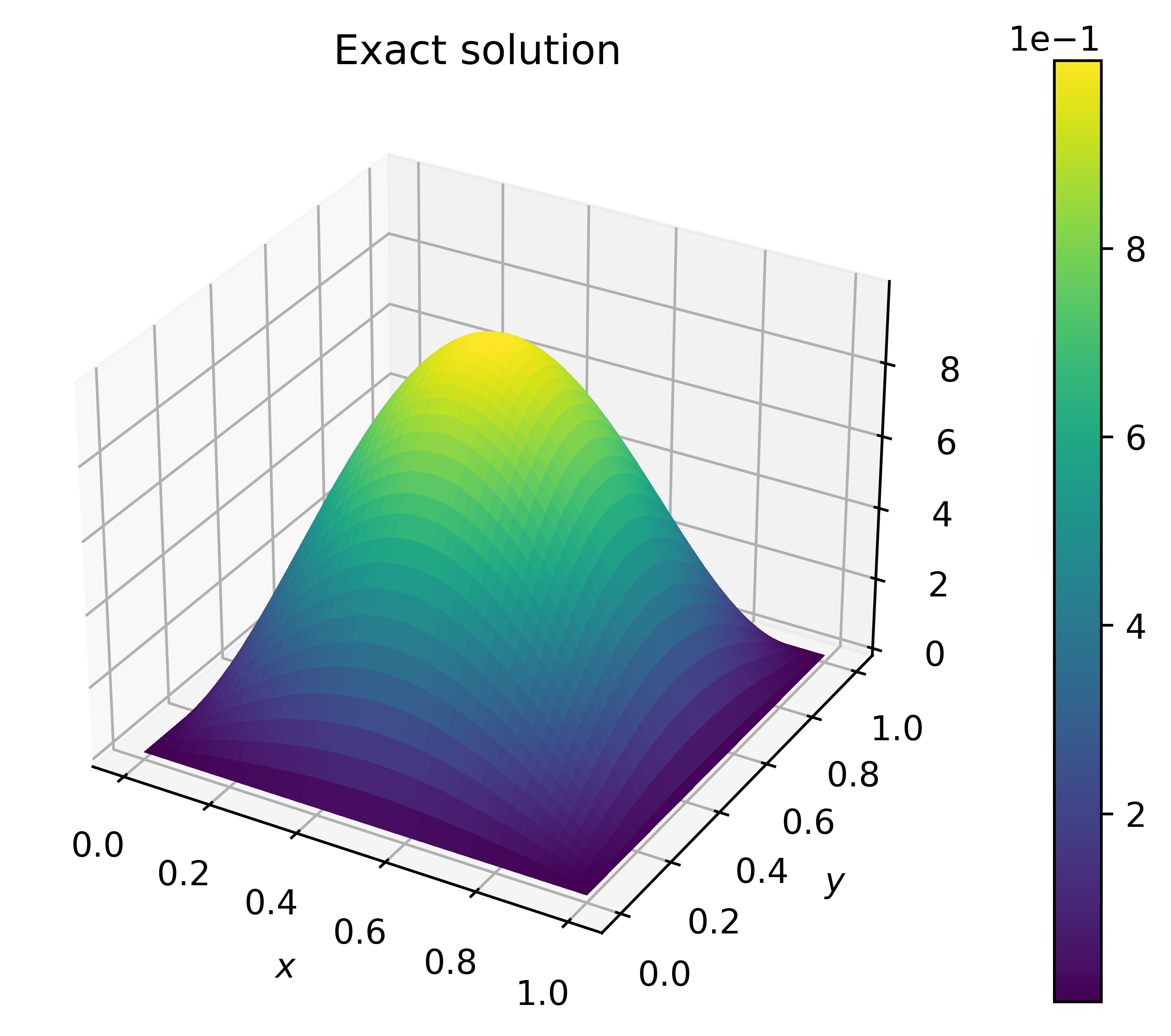}
    \caption*{(2b)}
  \end{subfigure}
  \hfill
  \begin{subfigure}{0.32\textwidth}
    \centering
    \includegraphics[width=.965\linewidth]{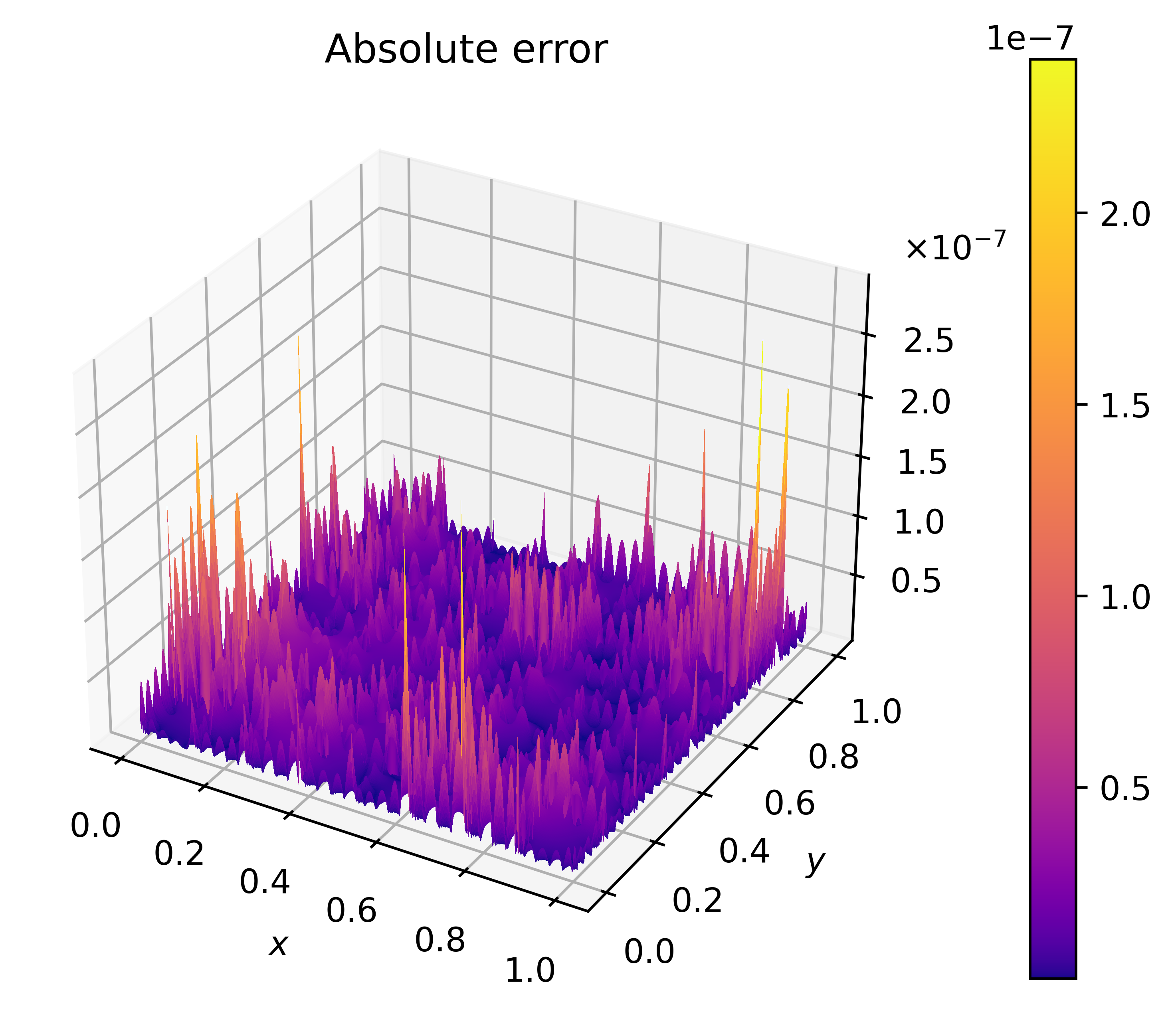}
    \caption*{(2c)}
  \end{subfigure}

  \begin{subfigure}{0.495\textwidth}
    \centering
    \includegraphics[width=\linewidth]{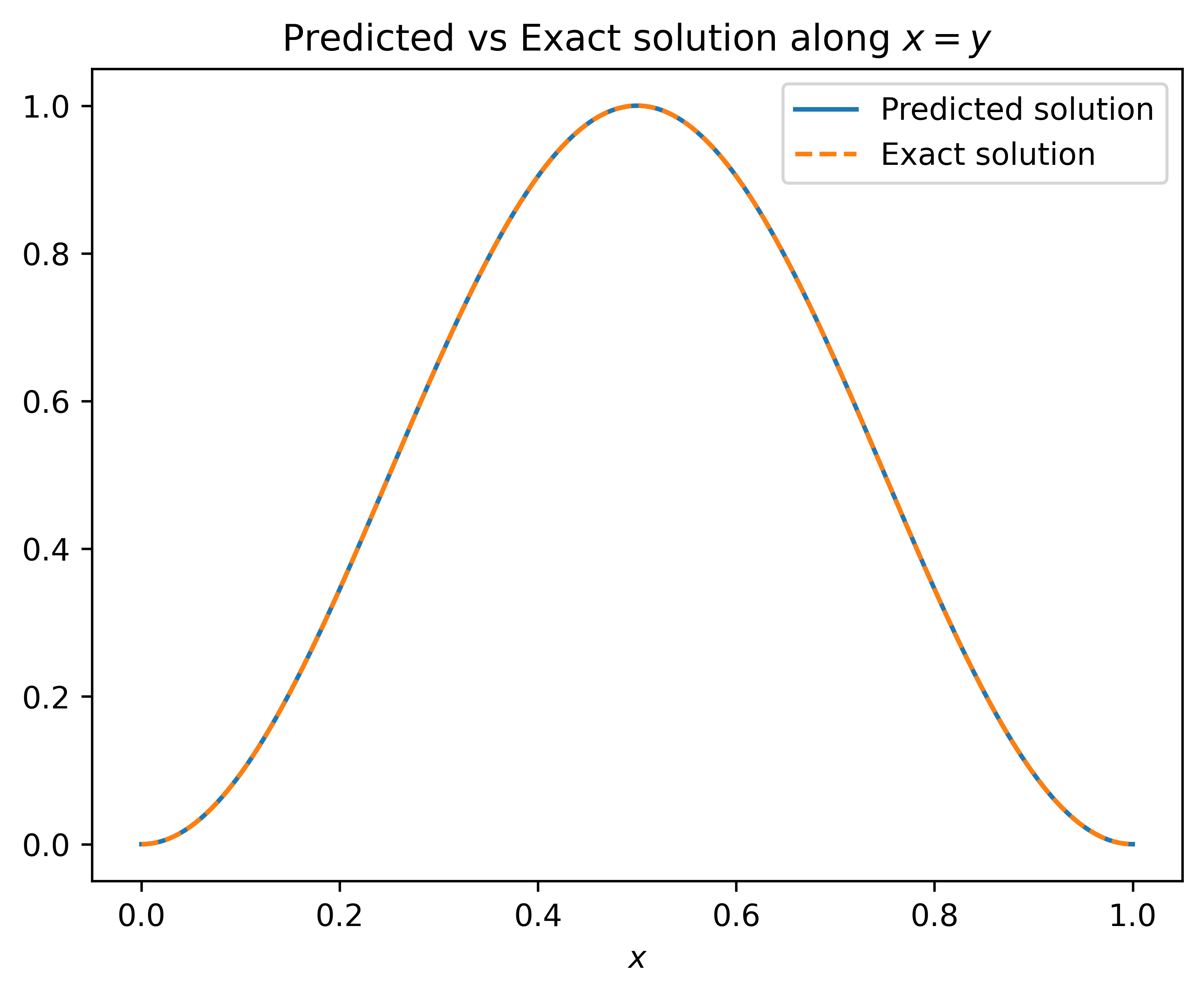}
    \caption*{(3a)}
  \end{subfigure}
  \hfill
  \begin{subfigure}{0.495\textwidth}
    \centering
    \includegraphics[width=.965\linewidth]{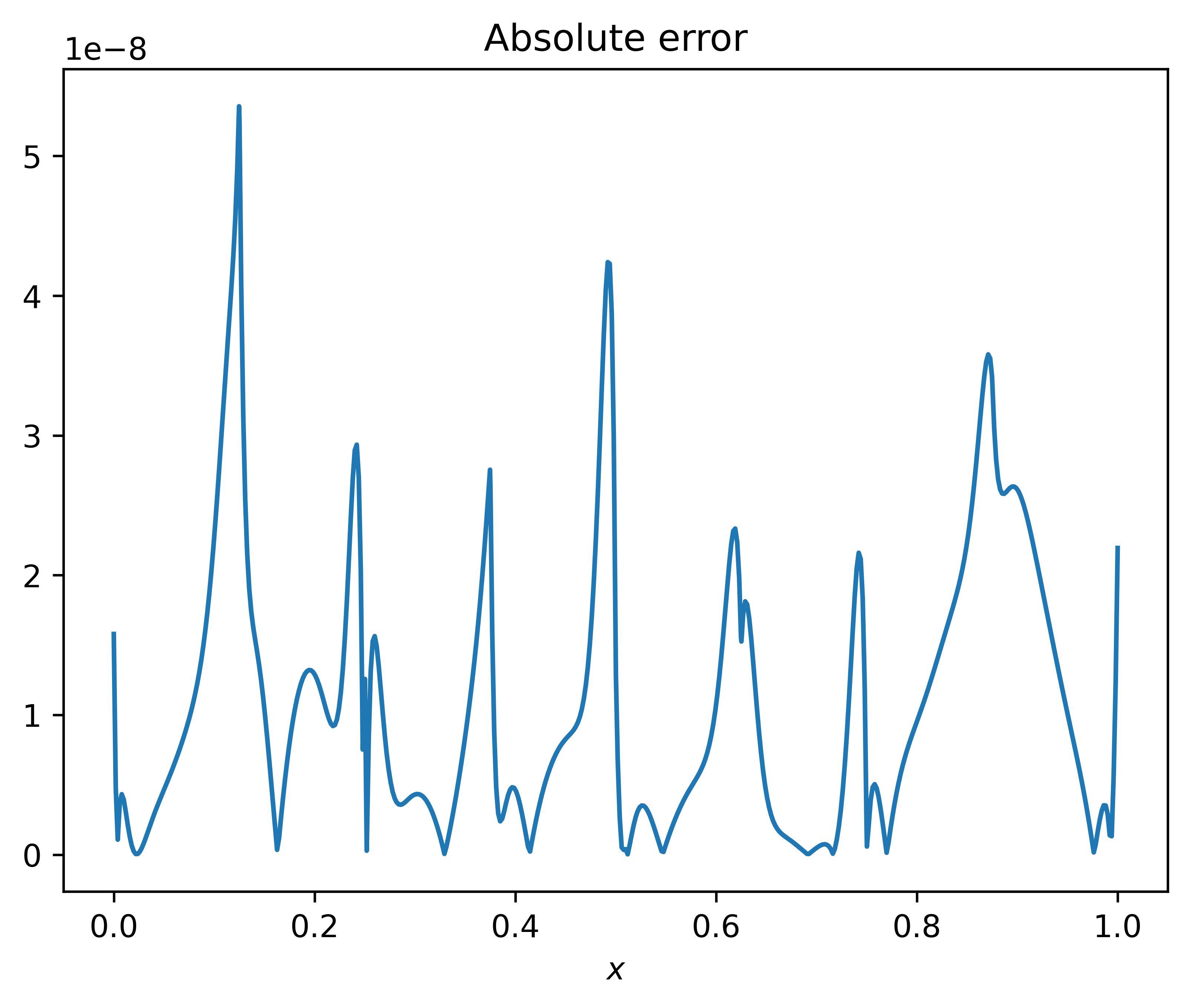}
    \caption*{(3b)}
  \end{subfigure}
  \caption{Solution and error of 2-D biharmonic equation with exact solution $u(x, y)=\sin(\pi x)\sin(\pi y)$.
  The DDELM-NN solution is computed on $8\times8$ subdomains.
  (1a) DDELM-NN solution 2D plot, (1b) Exact solution 2D plot, (1c) Absolute error of DDELM-NN solution 2D plot.
  (2a) DDELM-NN solution 3D plot, (2b) Exact solution 3D plot, (2c) Absolute error of DDELM-NN solution 3D plot.
  (3a) DDELM-NN and exact solution along the line $x=y$, (3b) Absolute error of DDELM-NN solution along the line $x=y$.}
  \label{fig:bih}
\end{figure}

Finally, we consider the 2-D biharmonic equation
\begin{equation*}
  \left\{\begin{aligned}
  \Delta^2 u &= f &&\text{ in } \Omega=(0,1)^2, \\
  u &= 0 &&\text{ on } \partial \Omega,\\
  \Delta u &= 0 &&\text{ on } \partial \Omega,
  \end{aligned}\right.
\end{equation*}
with the exact solution $u(x,y)=\sin(\pi x)\sin(\pi y)$.
Continuity and flux conditions are given in~\eqref{eqn:biharmonic_coupled}.
The results presented in \cref{tab:biharmonic}, are generally consistent with other experiments.
The $8\times8$ subdomain DDELM-NN solution and error plots are shown in~\cref{fig:bih}.
\begin{table}
  \caption{
    Relative $L^2$ and $H^1$ errors, number of CG iterations, and wall-clock time solving the biharmonic equation with exact solution $u(x,y)=\sin(\pi x)\sin(\pi y)$ using a total of 65,536 neurons.
    $N$ is the number of subdomains.
    }
  \label{tab:biharmonic}
\centering
\begin{tabular}{ccccccccc}
  \toprule
  Method & $N$ & $L^2$ & $H^1$ & Iterations & Wall-clock time (sec)\\
  \midrule \midrule
  \multirow{4}{*}{DDELM}
  & $2\times2$ & 3.52e-08 & 4.39e-07 & 826 & 460.74 \\
  & $4\times4$ & 2.63e-08 & 4.78e-07 & 2,694 & 12.27 \\
  & $8\times8$ & 9.02e-08 & 1.58e-06 & 3,957 & 6.02 \\
  & $16\times16$ & 7.25e-07 & 1.24e-05 & 5,809 & 6.21 \\
  \midrule
  \multirow{4}{*}{DDELM-CS}
  & $2\times2$ & 1.85e-08 & 2.89e-07 & 659 & 455.96 \\
  & $4\times4$ & 1.06e-08 & 1.85e-07 & 1,913 & 12.81 \\
  & $8\times8$ & 5.26e-08 & 1.07e-06 & 2,544 & 6.06 \\
  & $16\times16$ & 2.39e-07 & 5.56e-06 & 3,003 & 7.08 \\
  \midrule
  \multirow{4}{*}{DDELM-NN}
  & $2\times2$ & 1.50e-08 & 2.57e-07 & 86 & 909.09 \\
  & $4\times4$ & 9.27e-09 & 1.50e-07 & 169 & 16.41 \\
  & $8\times8$ & 3.47e-08 & 7.49e-07 & 294 & 1.52 \\
  & $16\times16$ & 2.16e-07 & 5.28e-06 & 736 & 2.29 \\
  \bottomrule
\end{tabular}
\end{table}

\section{Conclusions}
\label{Sec:Conc}
In this paper, we proposed a coarse space and a Neumann-Neumann acceleration for DDELM.
The coarse space formulation is based on partitioning of the auxiliary interface variables into coarse and non-coarse components.
A series of eliminations embeds a global problem into a system on the non-coarse variables.
We devised a change of variables that improves convergence, especially when used in conjunction with Neumann-Neumann acceleration.
By observing that the key term in DDELM is a Dirichlet to Neumann map, we proposed the Neumann-Neumann acceleration.
Numerical experiments validated two compoenents of the new method.
Experiments on several PDEs showed that the proposed method reduces both the number of CG iterations and wall-clock time.
An analysis of the spectral properties of the systems in question may enable further research on different coarse spaces and scalability.

\bibliographystyle{siam}
\bibliography{ddelm_ca}

\end{document}